\documentclass[11pt]{article}
\usepackage{graphicx}
\usepackage{amssymb,amsfonts,amsthm}
\usepackage{amsmath}
\usepackage{amscd}
\usepackage{color}

\usepackage{epsfig}

\newcommand{\atgs}{asymptotically tree-graded }
\newcommand{\atg}{asymptotically tree-graded}

\newcommand{\atgrs}{asymptotically tree-graded with respect to }
\newcommand{\wrt}{with respect to }
\newcommand{\uas}{$\omega$-almost surely}
\newcommand{\uass}{$\omega$-almost surely }

\newtheorem{theorem}{Theorem}[section]
\newtheorem{lemma}[theorem]{Lemma}

\newtheorem{cor}[theorem]{Corollary}
\newtheorem{proposition}[theorem]{Proposition}
\theoremstyle{definition}
\newtheorem{definition}[theorem]{Definition}

\newtheorem{rmks}[theorem]{Remarks}
\newtheorem{rmk}[theorem]{Remark}

\newtheorem{n}[theorem]{Notation}
\newtheorem{cvn}[theorem]{Convention}

\newcommand{\fq}{\mathfrak{Q}}
\newcommand{\cgs}{\mathrm{Cayley} (G,S)}
\newcommand{\cgsh}{\mathrm{Cayley} (G,\, S\cup \mathcal{EH})}
\newcommand{\dsh}{{\mathrm{dist}}_{S\cup \mathcal{EH}}}
\newcommand{\ds}{{\mathrm{dist}}_S}

\newcommand{\hh}{{\mathcal H}}
\def\mh{{\mathcal{EH}}} 

\newcommand{\Con}{{\mathrm{Con}}}

\newcommand{\dist}{{\mathrm{dist}}}
\newcommand{\oo}{{\mathcal O}}

\newcommand{\pp}{{\mathcal P}}

\newcommand{\nn}{{\mathcal N}}
\newcommand{\onn}{\overline{{\mathcal N}}}
\newcommand{\aaa}{{\mathcal A}}
\newcommand{\aau}{{\mathcal A}_\omega}
\newcommand{\diam}{{\mathrm{diam}}}

\def\sltwor{\mathrm{SL}(2,{\mathbb{R}})}   
\def\calt{\mathcal{T}}   
\def\calb{\mathcal{B}}   
\def\call{\mathcal{L}}   
\def\cali{\mathbf{A}}   
\def\calk{\mathcal{K}}   
\def\calf{\mathcal{F}}   

\newcommand {\N}{\mathbb{N}} 
\newcommand {\free}{\mathbb{F}} 
\newcommand {\hip}{\mathbb{H}} 
\newcommand {\sph}{\mathbb{S}} 
\newcommand {\q}{\mathfrak q} 
\newcommand {\g}{\mathfrak g} 
\newcommand {\pgot}{\mathfrak p}

\newcommand {\bc}{{\textbf{(C)}}}
\newcommand {\bcs}{{\textbf{(C)}} }
\newcommand {\bp}{{\textbf{(P)}}}
\newcommand {\bps}{{\textbf{(P)}} }
\newcommand {\me}{\medskip}

\newcommand {\iv}{^{-1}}
\newcommand{\lio}[1]{\lim_\omega\left( #1 \right)}
\newcommand{\co}[1]{\Con_\omega\left( #1 \right)}

\newcommand{\nk}[1]{{\mathcal{N}}_\varkappa \left( #1\right)}

\newcommand{\stab}{{\mathrm{Stab}}}
\newcommand {\fn}{\footnote}
\newcommand{\bB}{\bar{B}}
\begin{document}

\title{Relatively hyperbolic groups: geometry and quasi-isometric invariance}
\author{Cornelia Dru\c{t}u}
\date{}

\maketitle

\begin{abstract}
In this paper it is proved that relative hyperbolicity is an
invariant of quasi-isometry. As a byproduct of the arguments, simplified
definitions of relative hyperbolicity are obtained.
In particular we obtain a new definition
very similar
to the one of hyperbolicity, relying on the existence for every
quasi-geodesic triangle of a central left coset of peripheral
subgroup.
\end{abstract}

\tableofcontents


\section{Introduction}

\subsection{Rigidity result}

M. Gromov asked
(\cite{Gromov:hyperbolic},\cite{Gromov:Asymptotic}) what
properties of infinite finitely generated groups are invariant by
quasi-isometry. Such properties are sometimes called
\emph{geometric}, while a class of groups defined by a geometric
property is called \emph{rigid}.

\newpage

\textit{Examples of  rigid/non-rigid classes of groups:}
\begin{enumerate}
    \item the class of virtually nilpotent groups is rigid
    \cite{Gromov:PolynomialGrowth};
    \item  the class of virtually solvable groups is not rigid
    \cite{Dyubina:solvable}; but smaller classes of virtually solvable
    groups are rigid (\cite{FarbMosher:BSOne}, \cite{FarbMosher:BSTwo},
    \cite{EskinFisherWhyte:Sol});
    \item amenability is a geometric property;
    \item property (T) is not geometric (see for instance \cite{Valette:bourbakiT});
    \item hyperbolicity is a geometric property \cite{Gromov:hyperbolic};
    \item different classes of lattices of semisimple groups are rigid
    (this statement includes many deep results of different
    authors; see \cite{Farb:Lattices} and \cite{Drutu:grenoble} for surveys of these results).
\end{enumerate}

Recall that a group is said to \emph{virtually} satisfy a property
(P) if a finite index subgroup of it has property (P).

The present paper focuses on the class of relatively hyperbolic
groups.\fn{By relatively hyperbolic group we mean what is
sometimes called in the literature \textit{strongly relatively
hyperbolic group}, in contrast with \textit{weakly relatively
hyperbolic group}.} This notion was introduced by M. Gromov in
\cite{Gromov:hyperbolic}. Other definitions, as well as
developments of the theory of relatively hyperbolic groups can be
found in \cite{Bowditch:RelHyp}, \cite{Farb:RelHyp},
\cite{Dahmani:thesis}, \cite{Yaman:RelHyp},
\cite{DrutuSapir:TreeGraded}, \cite{Osin:RelHyp}. In
$\S$~\ref{marh} and $\S$~\ref{snewdef} we discuss in more detail
different ways to define relative hyperbolicity.

\me

\textit{Examples of relatively hyperbolic groups:}
\begin{enumerate}
    \item a hyperbolic group is hyperbolic relative to $\{
    1\}$;
    \item an amalgamated product $A*_F B$, where $F$ is finite, is hyperbolic relative to $A$ and $B$; more
    generally, fundamental groups of finite graphs of groups with
    finite edge groups are hyperbolic relative to the vertex
    groups \cite{Bowditch:RelHyp};
    \item fundamental groups of complete finite volume manifolds of pinched
    negative sectional curvature are hyperbolic relative to the
    fundamental groups of their cusps (\cite{Bowditch:RelHyp}, \cite{Farb:RelHyp});
    \item fundamental groups of (non-geometric) Haken manifolds with at least one
    hyperbolic component are hyperbolic relative to fundamental
    groups of maximal graph-manifold components and to fundamental groups of
    tori and Klein bottles not contained in a graph-manifold component;
    \item  fully residually free groups, also known as limit
    groups, are hyperbolic relative to their maximal
    Abelian non-cyclic subgroups \cite{Dahmani:combination}. Moreover they are CAT(0) with isolated flats
    \cite{AlibegovicBestvina:limit}.
\end{enumerate}

Note that there are also some interesting examples of groups
displaying a sort of ``intermediate'' relative hyperbolicity: they
are weakly relatively hyperbolic, not (strongly) relatively
hyperbolic, but nevertheless they have some common features with
(strongly) relatively hyperbolic groups, for instance their
asymptotic cones have a similar metric structure. Such groups are
the mapping class groups of surfaces of complexity at least two
(\cite{Behrstock:asymptotic}, \cite{BehrstockDrutuMosher:thick1}),
fundamental groups of $3$-dimensional graph manifolds
(\cite{KapovichLeeb:3manifolds}, \cite{KKL:QI},
\cite{BehrstockDrutuMosher:thick1}), as well as many Artin groups
(\cite{KapovichSchupp:Artin}, \cite{BehrstockDrutuMosher:thick1}).

Recently, relatively hyperbolic groups have been used to construct
examples of infinite finitely generated groups with unusual
properties. Thus in \cite{Osin:twoconjugacy} it is proved that
there exist torsion-free two-generated groups with exactly two
conjugacy classes.

\begin{cvn}\label{conv1}
Throughout the paper all relatively hyperbolic groups are assumed
to be \emph{finitely generated} and hyperbolic relative to
\emph{finitely many proper subgroups of finite type}.
\end{cvn}

We also use the following terminology: if a group $G$ is
hyperbolic relative to subgroups $H_1,...H_m$ then the subgroups
$H_1,\dots , H_m$ are called \textit{peripheral subgroups}.


The present paper gives an affirmative answer to the question
whether relative hyperbolicity is a quasi-isometry invariant
(formulated also as Problem 1.15 in \cite{DrutuSapir:TreeGraded}).

\begin{theorem}[relative hyperbolicity is geometric, Theorem \ref{qir}]\label{iqir}
Let $G$ be a group hyperbolic relative to a family of subgroups
$H_1,...,H_n$. If a group $G'$ is quasi-isometric to $G$ then $G'$
is hyperbolic relative to $H_1',...,H_m'$, where each $H_i'$ can
be embedded quasi-isometrically in $H_j$ for some $j=j(i)\in
\{1,2,...,n\}$.
\end{theorem}

Rigidity has previously been proved for some sub-classes of
relatively hyperbolic groups (with stronger versions of rigidity
theorems): non-uniform lattices in rank one semisimple groups
different from $\sltwor$ \cite{Schwartz:RankOne}, fundamental
groups of non-geometric Haken manifolds with at least one
hyperbolic component (\cite{KapovichLeeb:nonpc},
\cite{KapovichLeeb:Haken}), fundamental groups of graphs of groups
with finite edge groups \cite{PapasogluWhyte:ends}.

In the full generality assumed in Theorem \ref{iqir}, the stronger
statement that each subgroup $H_i'$ is quasi-isometric to some
subgroup $H_j$ cannot hold. This can be seen in the example when
$G=G'=A*B*C*D$, with $G$ hyperbolic relative to $\{A*B,C*D\}$ and
$G'$ hyperbolic relative to $\{A,B,C,D\}$. In
\cite{BehrstockDrutuMosher:thick1} it is shown that if in Theorem
\ref{iqir} it is moreover assumed that each peripheral subgroup
$H_i$ is \emph{not relatively hyperbolic} then the rigidity result
holds, moreover each $H_i'$ is quasi-isometric to some $H_j$. This
generalizes previous results from \cite{DrutuSapir:TreeGraded}.
The proof of Theorem \ref{iqir} is completely different from the
proofs in \cite{DrutuSapir:TreeGraded} and in
\cite{BehrstockDrutuMosher:thick1}. The main ingredient in
\cite{BehrstockDrutuMosher:thick1} is the following result, proved
using results from the present paper: given a group $G$ hyperbolic
relative to $H_1,...,H_n$, every quasi-isometric embedding into
$G$ of a group which is not relatively hyperbolic has its image in
a bounded radius tubular neighborhood of a left coset $gH_i\, $;
moreover the radius of the neighborhood depends only on $G$,
$H_1,...,H_n$ and on the constants of quasi-isometry, \emph{not}
on the domain of the quasi-isometry \cite[Theorem
4.1]{BehrstockDrutuMosher:thick1}.

The main steps in the proof of Theorem \ref{iqir} are explained in
what follows.

\subsection{Metric and algebraic relative hyperbolicity}\label{marh}

In order to study rigidity it is necessary to have a definition of
relative hyperbolicity of a group only in terms of its Cayley graphs.
Most definitions (except the ones in \cite{DrutuSapir:TreeGraded}
and in \cite{Osin:RelHyp}) use not only a Cayley graph of the
group but also a metric space obtained from this graph by gluing
to each left coset of a peripheral subgroup some geometric object
(a hyperbolic horoball \cite{Gromov:hyperbolic}, countably many
edges with one common endpoint \cite{Farb:RelHyp} etc).

In what follows, we recall definitions provided in
\cite{DrutuSapir:TreeGraded}.

\me

A complete geodesic metric space $\free$ is {\em tree-graded with
respect to a collection} $\pp$ of closed geodesic subsets (called
{\it{pieces}}), if the following two properties are satisfied:
\begin{enumerate}

\item[($T_1$)] two different pieces have at most one point in
common;

\item[($T_2$)] any simple non-trivial geodesic triangle is contained
in one piece.
\end{enumerate}

A similar, though not equivalent, notion has been introduced in
\cite{KKL:QI} under the name of space of type I.

\me

A metric space $X$ is \emph{asymptotically tree-graded with
respect to a collection of subsets }$\aaa $ if every asymptotic
cone of $X$ is tree-graded with respect to the collection of limit
sets of sequences in $\aaa$. A definition of asymptotic cones
of metric spaces, and of limit sets can be found in Section
\ref{prelac}.



\me

Equivalently, $X$ is \atgs \wrt $\aaa$ if the following three
geometric properties are satisfied (for details see Theorem 4.1 in
\cite{DrutuSapir:TreeGraded} or Theorem \ref{tgi} in this paper):
\begin{itemize}
\item[$(\alpha_1)$] finite radius tubular neighborhoods of distinct elements
in $\aaa$ are either disjoint or intersect in sets of uniformly
bounded diameter;
\item[$(\alpha_2)$] a geodesic with endpoints at distance at most
one third of its length from a set $A$ in $\aaa$ intersects a
tubular neighborhood of $A$ of uniformly bounded radius;
\item[$(\alpha_3)$] any fat geodesic polygon is contained in a
tubular neighborhood of a set $A$ in $\aaa$ of uniformly bounded
radius (here the meaning of ``fat'' is the contrary of ``thin'' in
its metric hyperbolic sense; see Definition \ref{fatpoly}).
\end{itemize}

The space $X$ is \emph{properly} \atgrs $\aaa$ if it is not
contained in any finite radius tubular neighborhood of a subset in
$\aaa$.

\begin{cvn}\label{conv2}
In what follows we assume that all \atgs metric spaces are
properly \atg.
\end{cvn}

The notion of \atgs metric space is a metric version for the
relative hyperbolicity of groups. Other similar notions can be
found in \cite{BrockFarb:curvature}, and in \cite{HruskaKleiner}
in the context of $\mathrm{CAT(0})$ metric spaces. The fact that
the metric definition is coherent with the definition for groups
is illustrated by the following result.

\begin{theorem}[\cite{DrutuSapir:TreeGraded}, Theorem
1.11 and Appendix]\label{dsrh} A finitely generated group $G$ is
hyperbolic  relative to $H_1,...,H_m $ if and only if $G$ is
asymptotically tree-graded with respect to the collection of left
cosets $\call=\{ gH_i \; ;\; g\in G/H_i\, ,\, i\in
\{1,2,...,m\}\}$.
\end{theorem}






The equivalence in Theorem \ref{dsrh} suggests the following
question, which appears as Problem 1.16 in
\cite{DrutuSapir:TreeGraded}: if a group is asymptotically
tree-graded in a metric sense, that is with respect to a
collection of subsets $\aaa$, does it follow that it is relatively
hyperbolic with respect to some finite family of subgroups ? The
implication was previously known to be true only under some
restrictive metric conditions on $\aaa$ (see \cite[Theorem
5.13]{DrutuSapir:TreeGraded} and
\cite{BehrstockDrutuMosher:thick1}).

We answer this question in the affirmative.

\begin{theorem}[Theorem \ref{atgrh}]\label{t3}

Let $G$ be an infinite finitely generated group asymptotically
tree-graded with respect to a collection of subsets $\aaa$. Then
$G$ is relatively hyperbolic with respect to some subgroups
$H_1,...,H_m$, such that every $H_i$ is contained in a bounded
radius tubular neighborhood of a set $A_i\in \aaa$.

\end{theorem}

Theorem \ref{t3} implies Theorem \ref{iqir}. Indeed, a group
quasi-isometric to a relatively hyperbolic group is \atgs as a
metric space with respect to the images by quasi-isometry of the
left cosets of peripheral subgroups \cite[Theorem
5.1]{DrutuSapir:TreeGraded}.

An outline of the proof of Theorem \ref{t3} will be given in the
following sections.

Theorem \ref{t3} is optimal  in the sense that if the group $G$
and the collection $\aaa$ satisfy less properties than those
required for \atgs metric spaces then the group $G$ may not be
relatively hyperbolic. This is shown by the examples of groups
constructed in \cite[$\S 7.1$]{BehrstockDrutuMosher:thick1} and in
\cite{OlshanskiiOsinSapir:cutpoints}. These groups are not
relatively hyperbolic, although they contain a collection of
subsets $\aaa$ such that all the asymptotic cones are tree-graded
with respect to {\emph{some}} limits of sequences in $\aaa$. But
in each cone, not all the limits of sequences in $\aaa$ are
considered as pieces: there are limits which are geodesic lines,
and different such lines intersect in more than one point. The
subsets in $\aaa$ do not satisfy property $(\alpha_1)$.

\subsection{New definitions of relative hyperbolicity}\label{snewdef}

If a group has an \atgs structure equivariant with respect to left
translations, then a standard argument shows that the group is
relatively hyperbolic (see Proposition \ref{equivatg}). Thus, the
main step in the proof of Theorem \ref{t3} is to construct an
equivariant \atgs structure on a group out of an arbitrary \atgs
structure. A natural idea is to consider all the translated \atgs
structures $g\aaa =\{gA \; ;\; A\in \aaa\}$ of a given \atgs
structure $\aaa$ on a group $G$, and to take non-empty intersections of the
form $\bigcap_{g\in G} gA_g$, with $A_g\in \aaa$. To make such an
argument work, it is necessary that the \atgs properties behave
well with respect to intersections. The following modification of
the list of three geometric properties defining an \atgs metric
space ensures this good behavior with respect to intersections.

\begin{theorem}[Theorem \ref{newdef}]\label{t1}
Let $(X,\dist)$ be a geodesic metric space and let $\aaa$ be a
collection of subsets of $X$. The space $X$ is \atgs \wrt $\aaa$
if and only if property $(\alpha_1)$ and the following two
properties are satisfied:
\begin{itemize}
     \item[$(\beta_2)$] a geodesic with endpoints at distance at most
$\frac{1}{k}$ of its length from a set $A$ in $\aaa$ (with $k$
large enough) has its
    middle third contained in a tubular
neighborhood of~$A$ of uniformly bounded radius;
  \item[$(\beta_3)$] any fat geodesic hexagon is contained in a
   tubular neighborhood of a set $A$ in $\aaa$ of uniformly bounded
radius.
\end{itemize}

\end{theorem}

It is not difficult to replace property $(\alpha_2)$ by
$(\beta_2)$, using results in \cite{DrutuSapir:TreeGraded}. But
replacing $(\alpha_3)$ by $(\beta_3)$ requires extra work.
Property $(\beta_3)$ implies that $(T_2)$ holds in any asymptotic
cone for simple triangles whose edges are limits of sequences of
geodesics (Proposition \ref{limtr}). But generically a geodesic in
an asymptotic cone of a group is not limit of a sequence of
geodesics (see the example in the end of $\S \, $\ref{prelac}). In
order to ensure $(T_2)$ for an arbitrary geodesic triangle the
argument in \cite{DrutuSapir:TreeGraded} was to prove that such a
triangle can be approximated by a geodesic triangle which is limit
of a sequence of fat polygons with the same number $m$ of edges
(see Lemma \ref{lime} in this paper). The number $m$ of edges must
increase when the constant of approximation decreases. This
approximation result and property $(\alpha_3)$ imply~$(T_2)$. In
this paper we show (Corollary \ref{cbeta3}) that if property
$(T_1)$ holds in every asymptotic cone, an inductive argument
allows to deduce $(T_2)$ from $(\beta_3)$.

\bigskip

Asymptotically tree-graded metric spaces have a property that
strongly reminds of hyperbolic metric spaces. A metric space is
hyperbolic if and only if the edges of every quasi-geodesic
triangle intersect a ball of uniformly bounded radius \cite[$\S
6$]{Gromov:hyperbolic}. A space $X$ that is \atgrs a collection of
subsets $\aaa$ has the following property
\cite{DrutuSapir:TreeGraded}:
\begin{itemize}
    \item[$(*)$] the edges of any quasi-geodesic triangle in $X$ either intersect a finite radius ball
or a finite radius tubular neighborhood of a subset in $\aaa$.
Moreover, in the latter case the distance between the entrance
points into the tubular neighborhood of two edges with common
 origin is uniformly bounded.
\end{itemize}

If $(X,\aaa)$ satisfy property $(*)$ then the space $X$ is called
$(*)$--\emph{asymptotically tree-graded with respect to }$\aaa$
\cite{DrutuSapir:RD}. This notion is weaker than the notion of
\atgs metric space (see Remark \ref{starrh}, (2)). Property $(*)$
was essential in the proof of the fact that the property of Rapid
Decay transfers from the peripheral subgroups $H_1,..., H_m$ of a
relatively hyperbolic group to the group itself
\cite{DrutuSapir:RD}. A version of property $(*)$ in the context
of CAT(0) spaces appears in \cite{Hruska:isolated}, where it is
called the Relatively Thin Triangle Property.

A natural question to ask is under what additional conditions is a
$(*)$--asymptotically tree-graded metric space also asymptotically
tree-graded. The arguments used to prove Theorem \ref{t1} can be
adapted to answer this question.

\begin{theorem}[Theorem \ref{tstar}]\label{T2}
Let $(X,\dist)$ and $\aaa$ be as in Theorem \ref{t1}. The metric
space $X$ is asymptotically tree-graded with respect to $\aaa$ if
and only if $(X,\aaa )$ satisfy properties $(\alpha_1)$ and
$(\alpha_2)$, and moreover $X$ is $(*)$--asymptotically
tree-graded with respect to $\aaa$.
\end{theorem}

\subsection{Organization of the paper}\label{org}

Section \ref{prel} contains preliminaries on asymptotic cones,
 as well as notation used throughout
the paper.

In Section \ref{stgr} are recalled some basic facts about
tree-graded spaces. Proposition \ref{braid} proved in the same
section is very useful in different arguments deducing the general
property $(T_2)$ from $(T_1)$, and $(T_2)$ restricted to some
particular cases.

Section \ref{satg} begins with a short overview of properties of
\atgs metric spaces. In $\S \, $\ref{pt2} an induction argument
and Proposition \ref{braid} are used to show the following central
result. Denote by $(\Pi_3)$ the property $(T_2)$ restricted to
triangles with edges limits of sequences of geodesics. If in an
asymptotic cone $\co{X}$ of a metric space $X$ a collection $\aau$
of closed subsets satisfies $(T_1)$ and $(\Pi_3)$ then $\aau$
satisfies $(T_2)$ in full generality (Corollary \ref{pi3}).

This statement is the main ingredient in the proof of Theorem
\ref{t1}, given in $\S \, $\ref{ndr}. It also plays a central part
in the proof of Theorem \ref{T2} given in $\S \, $\ref{ndh}.
Another difficult step in the proof of Theorem \ref{T2} is to
deduce from properties $(*),\, (\alpha_1)$ and $(\alpha_2)$ the
fact that fat quadrilaterals are contained in finite radius
tubular neighborhoods of subsets in $\aaa$ (Lemma \ref{4fat}).
Once this last statement proved, from it as well as from property
$(*)$ and Proposition \ref{bigon} can be deduced property
$(\Pi_3)$. Corollary \ref{pi3} allows to finish the argument.

Theorem \ref{t3} is proved in Section \ref{sqirrh}. The first and
most difficult step of the proof is to construct from a given
\atgs structure on a group an equivariant \atgs structure. The
subsets in the new \atgs structure are indexed by equivalence
classes of fat hexagons. A simple argument then shows that the
existence of an equivariant \atgs structure implies that the group
is relatively hyperbolic (Proposition \ref{equivatg}). This
completes the proof of Theorem \ref{t3} and thus of Theorem
\ref{iqir}.

\me

{\bf Acknowledgement.} The author wishes to thank Mark Sapir
 and Jason Behrstock
for comments that helped improving the presentation of the paper.

\section{Preliminaries}\label{prel}

\subsection{Definitions and notation}\label{dn}

Let $Y$ be a subset in a metric space $(X,\dist)$. We denote by
$\nn_\delta (Y)$ the set $\{ x\in X \mid \dist (x,Y) < \delta \}$,
which we call the $\delta$\emph{-tubular neighborhood of} $Y$. We
denote by $\overline{\nn}_\delta (Y)$ the set  $\{ x\in X\mid
\dist (x,Y) \leq \delta \}$, called the $\delta$\emph{-closed
tubular neighborhood of $Y$}.

When $Y$ is a singleton $y$, we also use the notation
$B(y,\delta)$ and respectively $\overline{B}(y,\delta)$.

\begin{definition}\label{tk}
An action of a group $G$ on a metric space $X$ is called
$\calk$\emph{-transitive}, where $\calk$ is a non-negative
constant, if for every $x\in X$ the closed tubular neighborhood
$\onn_\calk (Gx)$ of the orbit of $x$ coincides with $X$.
\end{definition}

\medskip

An $(L,C)$--\textit{quasi-isometric embedding} of a metric space
$(X, \dist_X)$ into a metric space $(Y,\dist_Y)$ is a map $\q : X
\to Y$ such that for every $x_1,x_2\in X$,
\begin{equation}\label{qi}
\frac{1}{L}\dist_X (x_1,x_2)-C \leq \dist_Y (\q(x_1), \q(x_2))\leq
L\dist_X (x_1,x_2)+C\, ,
\end{equation} for some constants $L\geq 1$ and $C\geq 0$.

If moreover $Y$ is contained in the $C$--tubular neighborhood of
$\q (X)$ then $\q$ is called an $(L,C)$--\textit{ quasi-isometry}.
In this case there exists an $(L,C)$--quasi-isometry $\bar{\q} :Y
\to X$ such that $\bar{\q} \circ \q$ and $\q \circ \bar{\q}$ are
at uniformly bounded distance from the respective identity maps \cite{GhysHarpe:hyperboliques}. The quasi-isometry $\bar{\q}$ is called \textit{quasi-converse of} $\q$.

If $\q :[a,b]\to X$ is an $(L,C)$--quasi-isometric embedding then
$\q$ is called an $(L,C)$-\textit{quasi-geodesic (segment)} in
$X$. The same name is used for the image of $\q$.

\begin{n}\label{+-}
For every quasi-geodesic segment $\q$ in a metric space $X$, we
denote the origin of $\q$ by $\q_-$ and the endpoint of $\q $ by
$\q_+$.
\end{n}


\me

If $\q_i :[0,\ell_i]\to X\, ,\, i=1,2,$ are two quasi-geodesic
segments with $\q_1(\ell_1)=\q_2(0)$, then we denote by
$\q_1\sqcup \q_2$ the map $\q :[0, \ell_1+\ell_2]\to X$ defined by
$\q(t)=\q_1(t)$ for $t\in [0, \ell_1]$ and $\q (t)=\q_2
(t-\ell_1)$ for $t\in [\ell_1, \ell_1+\ell_2]$.

\me

If an $(L,C)$--quasi-geodesic $\q$ is
$L$--Lipschitz then $\q$ is called an {\em $(L,C)$--almost
geodesic}.

\subsection{Asymptotic cones of a metric space}\label{prelac}

The notion of asymptotic cone of a metric space was used
implicitly in \cite{Gromov:PolynomialGrowth}, and it was defined
in full generality and studied in \cite{DriesWilkie} and
\cite{Gromov:Asymptotic}. For the definition, one needs the notion
of \textit{non-principal ultrafilter}. This is a finitely additive
measure $\omega$ defined on the set of all subsets of $\N$ (or,
more generally, of a countable set) and taking values in
$\{0,1\}$, such that $\omega (F)=0$ for every finite subset $F$ of
$\N$.


\begin{cvn}
Throughout the paper all ultrafilters are non-principal, therefore
we will omit mentioning it each time.
\end{cvn}

\begin{n}
Let $A_n$ and $B_n$ be two sequences of objects and let
$\mathcal{R}$ be a relation that can be established between $A_n$
and $B_n$ for every $n\in \N$. We write $A_n \,
\mathcal{R}_\omega\, B_n$ if and only of $A_n\, \mathcal{R}\, B_n$
\uas, that is
$$\omega \left( \{ n\in \N \mid A_n\, \mathcal{R}\, B_n \}
\right)=1\, .$$

\emph{Examples:} $=_\omega\, ,\, <_\omega \, ,\, \subset_\omega $.
\end{n}

Given an ultrafilter $\omega$, an $\omega$--\textit{limit}
$\lim_\omega x_n$ of a sequence $(x_n)$ in a topological space $X$
is an element $x\in X$ such that for every neighborhood $\nn$ of
$x$, $x_n\in_\omega \nn$. In a Hausdorff separable space if the
$\omega$--limit of a sequence exists then it is unique. If $(x_n)$
is contained in a compact space then it has an $\omega$--limit
\cite{Bourbaki}.

Given a space $X$ one can define its \emph{ultrapower} $X^\omega$
as the quotient $X^\N / \approx$, where $(x_n)\approx (y_n)$ if
$x_n=_\omega y_n$.

Let now $(X,\dist )$ be a metric space, $e$ a fixed element in its
ultrapower $X^\omega$, $(e_n)$ a representative of $e$, and
$d=(d_n)$ a sequence of numbers in $(0,+\infty)$ such that
$\lim_\omega {d_n}=+\infty$.

Consider
\begin{equation}\label{se}
\mathcal{S}_e(X) = \left\{ (x_n) \in X^\N \: ;\:  \mbox{ there exists }
M_x\mbox{ such that } \dist (x_n, e_n)\leq_{\omega} M_x\, d_n
\right\} \, .
\end{equation}

Define the equivalence relation
$$
(x_n)\sim (y_n) \Leftrightarrow \lim_\omega \frac{\dist (x_n,
y_n)}{d_n} =0\, .
$$

The quotient space $\mathcal{S}_e (X) / \sim $ is denoted by
$\co{X; e,d }$ and it is called \textit{the asymptotic cone of $X$
with respect to the ultrafilter $\omega$, the scaling sequence $d$
and the sequence of observation centers $e$}. It is endowed with
the natural metric $\dist_\omega$ defined by
$$
\dist_\omega \left( x,y \right)= \lim_\omega \frac{\dist (x_n,
y_n)}{d_n}\, .
$$

Every asymptotic cone is a complete metric space.

A sequence of subsets $(A_n)$ in $X$ gives rise to a \textit{limit
subset} in the cone, defined by
$$
\lio{A_n} = \left\{ \lio{a_n} \mid a_n \in_\omega A_n \right\}\, .
$$

If $\lim_\omega \frac{\dist (e_n,A_n)}{d_n}=+\infty $ then
$\lio{A_n}=\emptyset$. Every non-empty limit subset $\lio{A_n}$ is
closed.

\medskip

If each set $A_n$ is a geodesic $\g_n$ with length of order
$O(d_n)$ and $\lio{\g_n}$ is non-empty, then it is a geodesic in
$\co{X; e,d }$. Therefore if $X$ is a geodesic space then every
asymptotic cone of it is geodesic.

\begin{definition}\label{limg}
We call a geodesic in $\co{X; e,d }$ which appears as $\lio{\g_n}$
with $\g_n$ geodesics in $X$ a \emph{limit geodesic}.
\end{definition}

Not every geodesic in $\co{X; e,d}$ is a limit geodesic, not even
in the particular case when $X$ is a group of finite type with a
word metric.

\medskip

\textit{Example of group with continuously many non-limit
geodesics in an asymptotic cone:}

\medskip

On the two-dimensional unit sphere  $\sph^2$ consider a family of
horizontal circles, and a family of vertical circles in parallel
planes, such that two consecutive circles in each family are at
spherical distance $\frac{\pi}{2^k}$, and such that the North and
the South points are on one vertical circle, and are at distance
$\frac{\pi}{2^k}$ from two respective horizontal circles.

The two families of circles compose a spherical grid $\Gamma_k'$.
We have that $\Gamma_k' \subset \Gamma_{k+1}'$. Let $\Gamma_k''$
be the graph obtained from $\Gamma_k'$ by joining with spherical geodesics
 all pairs of vertices not on the same vertical or horizontal circle, and at
distance at most $\frac{\pi}{\sqrt{2}^k}$. Let $\Gamma_k$ be the graph obtained from $\Gamma_k''$
by deleting all the vertical edges of length $\frac{\pi}{2^k}$
above the Equator, except the one having the East point $(1,0,0)$
as an endpoint, and replacing each of them by a path of double
length $\frac{\pi}{2^{k-1}}$. Let $\dist_k$ be the shortest-path
metric on $\Gamma_k$.

Proposition 7.26 from \cite{DrutuSapir:TreeGraded} applied to the
sequence of graphs $(\Gamma_k ,\dist_k)$, and Lemma 7.5 from the
same paper imply that there exists a two-generated and recursively
presented group $G$ with one asymptotic cone tree-graded, with all
pieces isometric to $\sph^2$. Moreover, from the construction of
$G$ it follows that in each of the pieces, for an appropriate
choice of the North, South and East points, all geodesics joining
North and South and not containing East are not limit geodesics.

The same argument as in \cite[$\S 7$]{DrutuSapir:TreeGraded}
allows in fact to construct a two-generated and recursively
presented group with continuously many non-homeomorphic asymptotic
cones with the property that continuously many geodesics in each
of them are not limit geodesics.

\me

\section{Tree-graded metric spaces}\label{stgr}

\subsection{Definition and properties}

 The notion of tree-graded metric space has been introduced in
 \cite{DrutuSapir:TreeGraded}. In this paper we use the following version
  of this notion. Recall that a subset $A$ in a geodesic metric
  space $X$ is called \emph{geodesic} if every two points in $A$
  can be joined by a geodesic contained in $A$.

\begin{definition}\label{deftgr}
Let $\free$ be a complete geodesic metric space and let $\pp$ be a
collection of closed geodesic subsets, called {\it{pieces}}.
Suppose that the following two properties are satisfied:
\begin{enumerate}

\item[($T_1$)] Every two different pieces have at most one point in common.

\item[($T_2$)] Every simple non-trivial geodesic triangle in $\free$ is contained
in one piece.
\end{enumerate}

Then we say that the space $\free$ is {\em tree-graded with
respect to }$\pp$.

When there is no risk of confusion as to the set $\pp$, we simply
say that $\free$ is \emph{tree-graded}.
\end{definition}

\begin{rmks}[pieces need not cover the space]
\begin{itemize}
    \item[(1)] In \cite{DrutuSapir:TreeGraded} trivial geodesic triangles are
allowed in property ($T_2$). This is equivalent to asking that
$\free$ is covered by the pieces in $\pp$. In the present paper we
remove this convention. The reason is that a main purpose when
introducing the notion of tree-graded space is to produce a
convenient notion of relatively hyperbolic metric space
 (called asymptotically tree-graded metric space
 in \cite{DrutuSapir:TreeGraded} and in this paper, see Definition \ref{asco}).
 The condition that pieces cover $\free$ produces some unnatural restrictions for
 a space to be asymptotically tree-graded (i.e. relatively
 hyperbolic) with respect to a list of subsets.
 See Remark \ref{ratg} for details.
    \item[(2)] Possibly $\pp$ is empty, in which case $\free$ is a
    real tree.

    \item[(3)] When a group $G$ acts transitively on $\free$
    (for instance when $\free$ is an asymptotic cone
     of a group) and $G$ permutes the pieces,
      the condition that pieces cover $\free$ is
     automatically satisfied.
\end{itemize}
\end{rmks}

All properties of tree-graded spaces in \cite[$\S
2.1$]{DrutuSapir:TreeGraded} hold with the new definition
\ref{deftgr}, as none of the proofs uses the property that pieces
cover the space. In particular one has the following results.



\begin{lemma}[\cite{DrutuSapir:TreeGraded}, $\S 2.1$]\label{ptx}
Let $x$ be an arbitrary point in $\free$ and let $T_x$ be the set
of points $y\in \free$ which can be joined to $x$ by a topological
arc intersecting every piece in at most one point.

The subset $T_x$ is a real tree and a closed subset of $\free$, and every topological arc joining two points in $T_x$ is
    contained in $T_x$. Moreover, for every $y\in T_x$, $T_y=T_x$.
\end{lemma}

\begin{definition}\label{ttrees}
A subset $T_x$ as in Lemma \ref{ptx} is called a \emph{transversal
tree} in $\free$.

\end{definition}

In \cite{KKL:QI} is defined the notion of space of type I, which
is equivalent to that of a tree-graded space with the extra
property that for every $x$ the transversal tree $T_x$ is a
geodesically complete tree which branches everywhere.

\begin{rmk}\label{removetree}
 One can ensure that pieces in a tree-graded space cover it by adding to the list
    of pieces the transversal trees. Thus a tree-graded space $\free$ with set of pieces $\pp$ in the  sense of
Definition
 \ref{deftgr} can be seen as tree-graded in the sense of
 Definition 2.1 in \cite{DrutuSapir:TreeGraded} with respect to a
 set of pieces $\pp'$ such that $\pp'\setminus \pp$ is a
 collection of real trees.
\end{rmk}







\subsection{Topological bigons contained in pieces}\label{sbigon}

\begin{definition}\label{bigon}
We call \emph{topological bigon ($\calt$-bigon, in short) formed
by two topological arcs $\g_1$ and $\g_2$} a union of a sub-arc
$\g_1'$ of $\g_1$ with a sub-arc $\g_2'$ of $\g_2$ such that
$\g_1'$ and $\g_2'$ have common endpoints $x$ and $y$. The
\emph{endpoints of the $\calt$-bigon} are the points $x$ and $y$.
The \emph{interior of the $\calt$-bigon} is the set $\g_1'\cup
\g_2'\setminus \{x,y\}$.

If $\g_1'$ and $\g_2'$ intersect only in their endpoints then the
$\calt$-bigon is called \emph{simple} (in fact it is a simple loop
in this case).
\end{definition}

Note that a $\calt$-bigon with non-empty interior cannot be
trivial, i.e. reduced to a point.

The results in this section are useful in arguments aiming to
prove property $(T_2)$ for a collection of closed subsets of a
metric space. In several contexts it proves necessary to deduce
from $(T_1)$, and $(T_2)$ satisfied only for some special type of
geodesic bigons, the general property $(T_2)$.


\begin{lemma}\label{lbigon}
Let $\g_1$ and $\g_2$ be two topological arcs with common
endpoints. Then every point $z\in \g_1 \setminus \g_2$ is in the
interior of a simple $\calt$-bigon formed by $\g_1$ and $\g_2$.
\end{lemma}

\proof For $i=1,2$, $\g_i :[0,\ell_i]\to Y$ is a topological
embedding. Let $t\in [0,\ell_1]$ be such that $\g_1(t)=z$. The set
$K=\g_1\iv \left( \g_2 \left( [0,\ell_2]\right)\right)$ is a
compact set not containing $t$. Let $r$ be the maximal element of
the compact set $K\cap [0,t]$, and let $s$ be the minimal element
of the compact set $K\cap [t,\ell_1]$. Then $\g_1(r)=\g_2(r')$ and
$\g_1(s)=\g_2 (s')$ for some $r' , s' \in [0,\ell_2]$. The union
of $\g_1$ restricted to $[r,s]$ with $\g_2$ restricted to
$[r',s']$ is a simple $\calt$-bigon formed by $\g_1$ and $\g_2$,
containing $z$ in its interior.\endproof

\begin{lemma}\label{convbigon}
Let $Y$ be a metric space and let $\calb$ be a collection of
subsets of $Y$, $\calb$ satisfying property $(T_1)$.

Let $\g_1$ and $\g_2$ be two topological arcs with common
endpoints and with the property that any non-trivial simple
$\calt$-bigon formed by $\g_1$ and $\g_2$ is contained in a subset
in $\calb$.

If $\g_1$ is contained in $B\in \calb$ then $\g_2$ is contained in
$B$.
\end{lemma}

\proof  Take $z$ an arbitrary point in $\g_2 \setminus \g_1$. By
Lemma \ref{lbigon} the point $z$ is in the interior of a simple
$\calt$-bigon formed by $\g_1$ and $\g_2$, of endpoints $z_1,z_2$.
By hypothesis this $\calt$-bigon is contained in a subset $B_z\in
\calb$. As $\{z_1,z_2\}$ is in $B\cap B_z$ it follows by $(T_1)$
that $B_z=B$ and that $z\in B$.\endproof

\bigskip

\begin{figure}[!ht]
\centering
\unitlength 1mm 
\linethickness{0.4pt}
\ifx\plotpoint\undefined\newsavebox{\plotpoint}\fi 
\begin{picture}(106.5,90.5)(0,0)
\thicklines
\multiput(52.5,89.5)(.0337209302,-.0392026578){1505}{\line(0,-1){.0392026578}}
\thinlines \put(52.75,89.25){\line(0,1){0}}
\multiput(52.5,89.75)(-.0337301587,-.0806878307){756}{\line(0,-1){.0806878307}}
\put(27,28.75){\line(0,1){0}}
\qbezier(27,28.5)(33.25,31.38)(31.5,39.75)
\qbezier(31.5,39.25)(53.13,57)(103.25,30.75)
\qbezier(54.5,46.25)(49.63,52.75)(55.25,58.25)
\qbezier(55.25,58.25)(58.25,52.13)(54.25,46.5)
\qbezier(55,58)(33.25,64.75)(52.5,77.5)
\qbezier(52.5,77.25)(65,70.75)(55.5,58.25)
\qbezier(52.5,89.75)(57.5,82.63)(52.5,77)
\qbezier(52.5,77.5)(49,76.88)(52.5,89.75)
\put(24.5,25.25){\makebox(0,0)[cc]{$u$}}
\put(106.5,28.5){\makebox(0,0)[cc]{$w$}}
\put(54.5,90.5){\makebox(0,0)[cc]{$v$}}
\put(54.5,43){\makebox(0,0)[cc]{$a$}}
\put(29.25,40.75){\makebox(0,0)[cc]{$x$}}
\put(40.25,72){\makebox(0,0)[cc]{$y$}}
\put(85,61.5){\makebox(0,0)[cc]{${\mathcal{L}}_2$}}
\put(25.75,57.25){\makebox(0,0)[cc]{${\mathcal{L}}_1$}}
\put(43.5,53){\makebox(0,0)[cc]{$B_1$}}
\put(59.75,57){\makebox(0,0)[cc]{$y_1$}}
\put(56.75,78){\makebox(0,0)[cc]{$y_2$}}
\put(50.75,69.25){\makebox(0,0)[cc]{$B_2$}}
\put(40.5,33.75){\makebox(0,0)[cc]{${\mathfrak{g}}_1$}}
\put(70.75,34.5){\makebox(0,0)[cc]{${\mathfrak{g}}_2$}}
\end{picture}
\caption{Step 1.} \label{figbraid}
\end{figure}

\begin{proposition}\label{braid}
Let $Y$ be a metric space and let $\calb$ be a collection of
closed subsets of $Y$, $\calb$ with property $(T_1)$.

Let $\call_1$ and $\g_1$ be two topological arcs with common
endpoints $u,v$. Let $\call_2$ and $\g_2$ be two, possibly
identical, topological arcs with common endpoints $v,w$. Assume
that:
\begin{itemize}
    \item[(1)]  $\call_1 \cap \call_2 =\{v\}$;
    \item[(2)]  $\g_1 \cap \g_2$ contains a point $a\neq v$;
    \item[(3)] all non-trivial simple $\calt$-bigons formed either by $\g_1$ and
    $\g_2$, or by  $\g_i$ and $\call_i$, $i=1,2$, are contained in
    a subset in $\calb$.
\end{itemize}

Then the $\calt$-bigon formed by $\g_1$ and $\g_2$ with endpoints
$a$ and $v$ is contained in a subset in $\calb$.
\end{proposition}

\proof \textsc{Step 1.}\quad  Let $\g_i'$ denote the sub-arc of
$\g_i$ of endpoints $a$ and $v$, $i=1,2$.

We prove that there exists $b\in \g_1' \cap \g_2' \setminus
\{a\}$, such that the $\calt$-bigon formed by $\g_1'$ and $\g_2'$
of endpoints $a,b$ is contained in some $B\in \calb$.

Hypothesis (1) implies that either $a\not\in \call_1$ or $a\not\in
\call_2$. Without loss of generality we may assume that $a\not\in
\call_1$. Then $a$ is in the interior of a simple $\calt$-bigon
formed by $\call_1$ and $\g_1$, of endpoints $x$ and $y$, with $y$
on $\g_1'$. Property (3) implies that this $\calt$-bigon is
contained in a set $B_1\in \calb$.

If $y\in \g_2'$ then take $b=y$.

Assume that $y\not\in \g_2'$. Then $y$ is in the interior of a
simple $\calt$-bigon formed by $\g_1'$ and $\g_2'$, of endpoints
$y_1,y_2$ (with $y_2$ closer to $v$ than $y_1$ on $\g_1'$). By (3)
this $\calt$-bigon is contained in some $B_2\in \calb$. The
intersection $B_1\cap B_2$ contains $\{y,y_1\}$ hence by $(T_1)$
we have that $B_1=B_2=B$. Take $b=y_2$.

The sub-arc of $\g_1'$ with endpoints $a$ and $b$ is contained in
$B$. By property (3) we can apply Lemma \ref{convbigon} and obtain
that the sub-arc of $\g_2'$ in between $a$ and $b$ is also
contained in $B$.

\me

\textsc{Step 2.}\quad Let $\mathcal{E}$ be the set of points $b\in
\g_1' \cap \g_2' \setminus \{a\}$, such that the $\calt$-bigon
formed by $\g_1'$ and $\g_2'$ of endpoints $a,b$ is contained in
some $B\in \calb$. We prove that there exists $c\in \mathcal{E}$
such that $\g_1$ between $c$ and $v$ contains no other point from
$\mathcal{E}$.

Note that by property $(T_1)$ of $\calb$ all $\calt$-bigons of
endpoints $a$ and $b$, for some $b\in \mathcal{E}$, are contained
in the same $B_0\in \calb$.

Let $\varphi : [0,\ell]\to Y$ be a parametrization of $\g_1$,
$\varphi(\ell)=v$, and let $r$ be $\varphi^{-1} (a)$. The
pre-image $\mathcal{E}'=\varphi\iv \left(\mathcal{ E}\right)$ is
contained in $(r,\ell]$. Let $T$ be the supremum of
$\mathcal{E}'$. Then $T=\lim t_n$ for some increasing sequence
$(t_n)$ in $\mathcal{E}'$, hence $c=\varphi(T)$ is the limit of
the sequence of points $b_n=\varphi (t_n)\in \mathcal{E}$.
Obviously $c\in \g_1' \cap \g_2' \setminus \{a\}$. Since $B_0\in
\calb$ is closed and $b_n\in B_0$, it follows that $c\in B_0$.
Thus the sub-arc of $\g_1'$ between $a$ and $c$ is completely
contained in $B_0$. By Lemma \ref{convbigon} and property (3), the
$\calt$-bigon formed by $\g_1'$ and $\g_2'$ of endpoints $a$ and
$c$ is in $B_0$.

\me

\textsc{Step 3.}\quad We prove that the point $c$ obtained in Step
2 coincides with $v$.

Assume that $c\neq v$. Step 1 applied to the point $c$ instead of
$a$ implies that there exists $d\in \g_1' \cap \g_2' \setminus
\{c\}$, $d$ between $c$ and $v$ on both $\g_1'$ and $\g_2'$, such
that the $\calt$-bigon formed by $\g_1'$ and $\g_2'$ of endpoints
$c,d$ is contained in some $B'\in \calb$.

Since $c\neq v$ it cannot be contained simultaneously in $\call_1$
and in $\call_2$. Assume that $c\not\in \call_1$. Then $c$ is in
the interior of a simple $\calt$-bigon formed by $\g_1$ and
$\call_1$. According to (3) this $\calt$-bigon is contained in
some $B''\in \calb$. The intersections $B_0\cap B''$ and $B'\cap
B''$ both contain non-trivial sub-arcs of $\g_1'$, therefore
$B_0=B''=B'$. Thus the point $d$ is in the set $\mathcal{E}$ and
it is strictly between $c$ and $v$ on $\g_1'$. This contradicts
the choice of $c$.

We conclude that $c=v$. \endproof

\section{Asymptotically tree-graded metric spaces}\label{satg}

\subsection{Definition and properties}

Let $(X,\dist)$ be a geodesic metric space and let $\mathcal{A}=\{
A_i \mid i\in I \}$ be a collection of subsets of $X$. In every
asymptotic cone $\co{X;e,d}$, we consider the collection
$\mathcal{A}_\omega$ of limit subsets
$$
\left\{ {\lio{ A_{i_n}}}\; ;\; i=(i_n)^\omega\in I^\omega \hbox{
such that }\exists M_i \mbox{ with the property } \dist\left(
e_{n},A_{i_n} \right)\leq _\omega M_i\, d_{n} \right\}\, .
$$


\begin{definition}\label{asco}
The metric space $X$ is \textit{asymptotically tree-graded (ATG)
with respect to} $\aaa$ if every asymptotic cone $\co{X; e,d}$ is
tree-graded with respect to $\aaa_\omega $.
\end{definition}

Following Convention \ref{conv2}, in the rest of the paper we
shall assume that all ATG metric spaces are proper, that is no
subset $A\in \aaa$ contains $X$ in a tubular neighborhood of it.

The ATG property is meant as an extension of the property of
(strong) relative hyperbolicity from groups to metric spaces.
Theorem \ref{dsrh} emphasizes that it is the correct property to
work with.

\begin{rmk}\label{tna}
Let $X$ be ATG \wrt $\aaa=\{ A_i \; ;\; i\in I\}$.
\begin{itemize}
    \item[(1)] It is easy to see that for every $\tau >0$,
    the space $X$ is ATG \wrt $\{\nn_\tau (A_i)\; ;\;
i\in I\}$.

    \medskip

    \item[(2)] More generally, let $\calb$ be a collection of subsets
of $X$ such that there exists a constant $K\geq 0$ and a bijection
$\phi :\aaa \to \calb$ verifying $\dist_H(A, \phi(A))\leq K$. Then
$X$ is ATG \wrt $\calb$.
\end{itemize}
\end{rmk}








The notion of ATG metric space can also be defined by a list of
geometric conditions, without involving asymptotic cones. First we introduce
some notation and terminology.

\medskip

\begin{n}\label{breve}
Given $\pgot$ a quasi-geodesic and $r>0$ we denote by
$\breve{\pgot}_r$ the set $\pgot \setminus \nn_r \left( \left\{
\pgot_-\, ,\, \pgot_+ \right\} \right)$.
\end{n}

\medskip

We say that a metric space $P$ is a \emph{geodesic
(quasi-geodesic) $k$-gonal line} if it is a union of $k$ geodesics
(quasi-geodesics) $\q_1,...,\q_k$ such that
$(\q_i)_+=(\q_{i+1})_-$ for $i=1,...,k-1$. If moreover
$(\q_k)_+=(\q_{1})_-$ then we say that $P$ is a \emph{geodesic
(quasi-geodesic) $k$-gon}.

Let $P$ be a quasi-geodesic polygon, with set of vertices
$\mathcal{V}$. Points in $P\setminus \mathcal{V}$ are called
\textit{interior points of }$P$.

\medskip

\begin{n}\label{ox}
Given a vertex $x\in \mathcal{V}$ and $\q , \q'$ the consecutive
edges of $P$ such that $x=\q_+ =\q'_-$, we denote the polygonal
line $P\setminus ( \q\cup \q')$ by $\oo_{x}(P)$. When there is no
possibility of confusion we simply denote it by $\oo_{x}$.
\end{n}

\medskip

\begin{figure}[!ht]
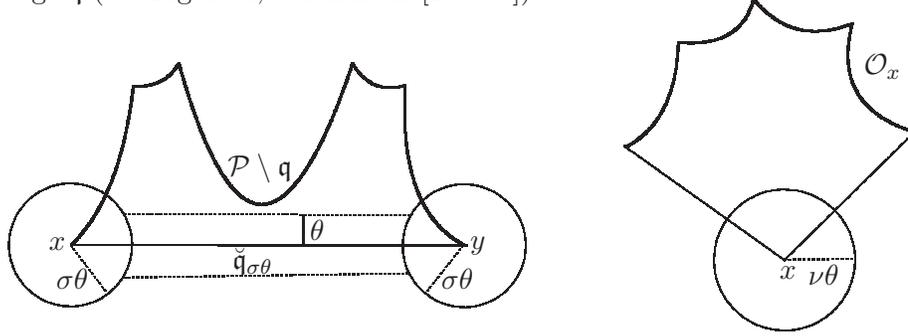

\centering \caption{Properties ($F_1$) and ($F_2$).} \label{fig1}
\end{figure}

Let $p\in P$. The \textit{inscribed radius in }$p$ with respect to
$P$ is either the distance from $p$ to the set $\mathcal{O}_p$, if
$p$ is a vertex, or the distance from $p$ to the set $P\setminus
\q$ if $p$ is an interior point contained in the edge $\q$ (see
Figure \ref{fig1}, taken from \cite{DrutuSapir:TreeGraded}).

\bigskip

\unitlength .3mm 
\linethickness{0.4pt}
\ifx\plotpoint\undefined\newsavebox{\plotpoint}\fi 
\begin{picture}(248.81,123.52)(0,0)
\multiput(47.32,42.92)(21.75,-.03125){8}{\line(1,0){21.75}}
\put(74.33,42.92){\line(0,1){1.138}}
\put(74.31,44.06){\line(0,1){1.136}}
\put(74.24,45.19){\line(0,1){1.132}}
\multiput(74.12,46.32)(-.03313,.22514){5}{\line(0,1){.22514}}
\multiput(73.95,47.45)(-.03035,.15969){7}{\line(0,1){.15969}}
\multiput(73.74,48.57)(-.03237,.1385){8}{\line(0,1){.1385}}
\multiput(73.48,49.68)(-.0305,.10962){10}{\line(0,1){.10962}}
\multiput(73.17,50.77)(-.03186,.09841){11}{\line(0,1){.09841}}
\multiput(72.82,51.86)(-.03294,.08891){12}{\line(0,1){.08891}}
\multiput(72.43,52.92)(-.031392,.074967){14}{\line(0,1){.074967}}
\multiput(71.99,53.97)(-.032193,.068685){15}{\line(0,1){.068685}}
\multiput(71.51,55)(-.032842,.063077){16}{\line(0,1){.063077}}
\multiput(70.98,56.01)(-.033361,.058025){17}{\line(0,1){.058025}}
\multiput(70.41,57)(-.03199,.050626){19}{\line(0,1){.050626}}
\multiput(69.81,57.96)(-.032371,.046785){20}{\line(0,1){.046785}}
\multiput(69.16,58.9)(-.032662,.043232){21}{\line(0,1){.043232}}
\multiput(68.47,59.8)(-.032872,.039929){22}{\line(0,1){.039929}}
\multiput(67.75,60.68)(-.033009,.036848){23}{\line(0,1){.036848}}
\multiput(66.99,61.53)(-.03308,.033962){24}{\line(0,1){.033962}}
\multiput(66.2,62.34)(-.034469,.032552){24}{\line(-1,0){.034469}}
\multiput(65.37,63.13)(-.037353,.032436){23}{\line(-1,0){.037353}}
\multiput(64.51,63.87)(-.040432,.032252){22}{\line(-1,0){.040432}}
\multiput(63.62,64.58)(-.045917,.03359){20}{\line(-1,0){.045917}}
\multiput(62.7,65.25)(-.049767,.03331){19}{\line(-1,0){.049767}}
\multiput(61.76,65.89)(-.053954,.032938){18}{\line(-1,0){.053954}}
\multiput(60.79,66.48)(-.058533,.032461){17}{\line(-1,0){.058533}}
\multiput(59.79,67.03)(-.063577,.031864){16}{\line(-1,0){.063577}}
\multiput(58.77,67.54)(-.074115,.033352){14}{\line(-1,0){.074115}}
\multiput(57.74,68.01)(-.081246,.032556){13}{\line(-1,0){.081246}}
\multiput(56.68,68.43)(-.08941,.03157){12}{\line(-1,0){.08941}}
\multiput(55.61,68.81)(-.10878,.03337){10}{\line(-1,0){.10878}}
\multiput(54.52,69.14)(-.12231,.032){9}{\line(-1,0){.12231}}
\multiput(53.42,69.43)(-.13898,.03023){8}{\line(-1,0){.13898}}
\multiput(52.31,69.67)(-.18683,.03253){6}{\line(-1,0){.18683}}
\multiput(51.18,69.87)(-.22563,.02965){5}{\line(-1,0){.22563}}
\put(50.06,70.02){\line(-1,0){1.133}}
\put(48.92,70.12){\line(-1,0){1.137}}
\put(47.79,70.17){\line(-1,0){1.138}}
\put(46.65,70.18){\line(-1,0){1.137}}
\put(45.51,70.14){\line(-1,0){1.134}}
\multiput(44.38,70.05)(-.22594,-.0272){5}{\line(-1,0){.22594}}
\multiput(43.25,69.91)(-.18717,-.0305){6}{\line(-1,0){.18717}}
\multiput(42.12,69.73)(-.1592,-.03282){7}{\line(-1,0){.1592}}
\multiput(41.01,69.5)(-.12265,-.03067){9}{\line(-1,0){.12265}}
\multiput(39.91,69.22)(-.10914,-.03218){10}{\line(-1,0){.10914}}
\multiput(38.82,68.9)(-.09791,-.03337){11}{\line(-1,0){.09791}}
\multiput(37.74,68.53)(-.081594,-.031672){13}{\line(-1,0){.081594}}
\multiput(36.68,68.12)(-.074473,-.032546){14}{\line(-1,0){.074473}}
\multiput(35.63,67.67)(-.06818,-.03325){15}{\line(-1,0){.06818}}
\multiput(34.61,67.17)(-.058882,-.031824){17}{\line(-1,0){.058882}}
\multiput(33.61,66.63)(-.054308,-.03235){18}{\line(-1,0){.054308}}
\multiput(32.63,66.04)(-.050126,-.032768){19}{\line(-1,0){.050126}}
\multiput(31.68,65.42)(-.046279,-.03309){20}{\line(-1,0){.046279}}
\multiput(30.76,64.76)(-.042722,-.033326){21}{\line(-1,0){.042722}}
\multiput(29.86,64.06)(-.039417,-.033485){22}{\line(-1,0){.039417}}
\multiput(28.99,63.32)(-.036334,-.033575){23}{\line(-1,0){.036334}}
\multiput(28.16,62.55)(-.033447,-.033601){24}{\line(0,-1){.033601}}
\multiput(27.35,61.74)(-.033407,-.036488){23}{\line(0,-1){.036488}}
\multiput(26.58,60.9)(-.033304,-.03957){22}{\line(0,-1){.03957}}
\multiput(25.85,60.03)(-.033129,-.042874){21}{\line(0,-1){.042874}}
\multiput(25.16,59.13)(-.032877,-.046431){20}{\line(0,-1){.046431}}
\multiput(24.5,58.21)(-.032537,-.050276){19}{\line(0,-1){.050276}}
\multiput(23.88,57.25)(-.0321,-.054456){18}{\line(0,-1){.054456}}
\multiput(23.3,56.27)(-.033525,-.062717){16}{\line(0,-1){.062717}}
\multiput(22.77,55.27)(-.032937,-.068332){15}{\line(0,-1){.068332}}
\multiput(22.27,54.24)(-.032204,-.074622){14}{\line(0,-1){.074622}}
\multiput(21.82,53.2)(-.031297,-.081739){13}{\line(0,-1){.081739}}
\multiput(21.41,52.13)(-.03292,-.09806){11}{\line(0,-1){.09806}}
\multiput(21.05,51.06)(-.03168,-.10928){10}{\line(0,-1){.10928}}
\multiput(20.74,49.96)(-.03011,-.12279){9}{\line(0,-1){.12279}}
\multiput(20.46,48.86)(-.03209,-.15935){7}{\line(0,-1){.15935}}
\multiput(20.24,47.74)(-.02964,-.18731){6}{\line(0,-1){.18731}}
\multiput(20.06,46.62)(-.0327,-.2826){4}{\line(0,-1){.2826}}
\put(19.93,45.49){\line(0,-1){1.135}}
\put(19.85,44.35){\line(0,-1){1.137}}
\put(19.81,43.22){\line(0,-1){1.138}}
\put(19.82,42.08){\line(0,-1){1.136}}
\put(19.88,40.94){\line(0,-1){1.133}}
\multiput(19.99,39.81)(.03068,-.22549){5}{\line(0,-1){.22549}}
\multiput(20.14,38.68)(.03339,-.18668){6}{\line(0,-1){.18668}}
\multiput(20.34,37.56)(.03086,-.13884){8}{\line(0,-1){.13884}}
\multiput(20.59,36.45)(.03256,-.12216){9}{\line(0,-1){.12216}}
\multiput(20.88,35.35)(.03079,-.09875){11}{\line(0,-1){.09875}}
\multiput(21.22,34.26)(.03198,-.08927){12}{\line(0,-1){.08927}}
\multiput(21.6,33.19)(.032929,-.081096){13}{\line(0,-1){.081096}}
\multiput(22.03,32.14)(.033692,-.073961){14}{\line(0,-1){.073961}}
\multiput(22.5,31.1)(.032156,-.06343){16}{\line(0,-1){.06343}}
\multiput(23.02,30.09)(.032729,-.058384){17}{\line(0,-1){.058384}}
\multiput(23.57,29.1)(.033185,-.053802){18}{\line(0,-1){.053802}}
\multiput(24.17,28.13)(.033538,-.049614){19}{\line(0,-1){.049614}}
\multiput(24.81,27.19)(.032191,-.043583){21}{\line(0,-1){.043583}}
\multiput(25.49,26.27)(.032437,-.040284){22}{\line(0,-1){.040284}}
\multiput(26.2,25.38)(.032608,-.037204){23}{\line(0,-1){.037204}}
\multiput(26.95,24.53)(.03271,-.034319){24}{\line(0,-1){.034319}}
\multiput(27.73,23.7)(.034113,-.032924){24}{\line(1,0){.034113}}
\multiput(28.55,22.91)(.036999,-.03284){23}{\line(1,0){.036999}}
\multiput(29.4,22.16)(.04008,-.032689){22}{\line(1,0){.04008}}
\multiput(30.29,21.44)(.043381,-.032463){21}{\line(1,0){.043381}}
\multiput(31.2,20.76)(.046933,-.032156){20}{\line(1,0){.046933}}
\multiput(32.14,20.12)(.053593,-.033521){18}{\line(1,0){.053593}}
\multiput(33.1,19.51)(.058177,-.033094){17}{\line(1,0){.058177}}
\multiput(34.09,18.95)(.063227,-.032552){16}{\line(1,0){.063227}}
\multiput(35.1,18.43)(.068832,-.031878){15}{\line(1,0){.068832}}
\multiput(36.13,17.95)(.080888,-.033436){13}{\line(1,0){.080888}}
\multiput(37.18,17.52)(.08906,-.03253){12}{\line(1,0){.08906}}
\multiput(38.25,17.13)(.09856,-.03141){11}{\line(1,0){.09856}}
\multiput(39.34,16.78)(.12196,-.03332){9}{\line(1,0){.12196}}
\multiput(40.43,16.48)(.13864,-.03173){8}{\line(1,0){.13864}}
\multiput(41.54,16.23)(.15983,-.02962){7}{\line(1,0){.15983}}
\multiput(42.66,16.02)(.22529,-.0321){5}{\line(1,0){.22529}}
\put(43.79,15.86){\line(1,0){1.132}}
\put(44.92,15.74){\line(1,0){1.136}}
\put(46.06,15.68){\line(1,0){1.138}}
\put(47.2,15.66){\line(1,0){1.137}}
\put(48.33,15.69){\line(1,0){1.135}}
\multiput(49.47,15.77)(.2828,.0309){4}{\line(1,0){.2828}}
\multiput(50.6,15.89)(.18749,.02847){6}{\line(1,0){.18749}}
\multiput(51.72,16.06)(.15955,.03109){7}{\line(1,0){.15955}}
\multiput(52.84,16.28)(.13835,.033){8}{\line(1,0){.13835}}
\multiput(53.95,16.54)(.10948,.031){10}{\line(1,0){.10948}}
\multiput(55.04,16.85)(.09826,.03231){11}{\line(1,0){.09826}}
\multiput(56.12,17.21)(.08876,.03335){12}{\line(1,0){.08876}}
\multiput(57.19,17.61)(.074822,.031736){14}{\line(1,0){.074822}}
\multiput(58.24,18.05)(.068537,.032509){15}{\line(1,0){.068537}}
\multiput(59.26,18.54)(.062925,.033132){16}{\line(1,0){.062925}}
\multiput(60.27,19.07)(.057871,.033627){17}{\line(1,0){.057871}}
\multiput(61.25,19.64)(.050479,.032222){19}{\line(1,0){.050479}}
\multiput(62.21,20.25)(.046636,.032585){20}{\line(1,0){.046636}}
\multiput(63.15,20.91)(.043081,.03286){21}{\line(1,0){.043081}}
\multiput(64.05,21.6)(.039778,.033055){22}{\line(1,0){.039778}}
\multiput(64.93,22.32)(.036696,.033178){23}{\line(1,0){.036696}}
\multiput(65.77,23.09)(.033809,.033236){24}{\line(1,0){.033809}}
\multiput(66.58,23.88)(.032393,.034618){24}{\line(0,1){.034618}}
\multiput(67.36,24.71)(.033731,.039207){22}{\line(0,1){.039207}}
\multiput(68.1,25.58)(.033593,.042512){21}{\line(0,1){.042512}}
\multiput(68.81,26.47)(.033379,.046071){20}{\line(0,1){.046071}}
\multiput(69.47,27.39)(.033081,.04992){19}{\line(0,1){.04992}}
\multiput(70.1,28.34)(.03269,.054105){18}{\line(0,1){.054105}}
\multiput(70.69,29.31)(.032192,.058682){17}{\line(0,1){.058682}}
\multiput(71.24,30.31)(.033677,.06797){15}{\line(0,1){.06797}}
\multiput(71.74,31.33)(.033012,.074268){14}{\line(0,1){.074268}}
\multiput(72.21,32.37)(.032182,.081395){13}{\line(0,1){.081395}}
\multiput(72.62,33.43)(.03115,.08956){12}{\line(0,1){.08956}}
\multiput(73,34.5)(.03287,.10893){10}{\line(0,1){.10893}}
\multiput(73.33,35.59)(.03144,.12246){9}{\line(0,1){.12246}}
\multiput(73.61,36.69)(.02959,.13912){8}{\line(0,1){.13912}}
\multiput(73.85,37.81)(.03167,.18698){6}{\line(0,1){.18698}}
\multiput(74.04,38.93)(.02861,.22576){5}{\line(0,1){.22576}}
\put(74.18,40.06){\line(0,1){1.134}}
\put(74.28,41.19){\line(0,1){1.728}}
\put(248.81,42.57){\line(0,1){1.138}}
\put(248.79,43.71){\line(0,1){1.136}}
\put(248.72,44.84){\line(0,1){1.132}}
\multiput(248.6,45.97)(-.03313,.22514){5}{\line(0,1){.22514}}
\multiput(248.43,47.1)(-.03035,.15969){7}{\line(0,1){.15969}}
\multiput(248.22,48.22)(-.03237,.1385){8}{\line(0,1){.1385}}
\multiput(247.96,49.33)(-.0305,.10962){10}{\line(0,1){.10962}}
\multiput(247.65,50.42)(-.03186,.09841){11}{\line(0,1){.09841}}
\multiput(247.3,51.51)(-.03294,.08891){12}{\line(0,1){.08891}}
\multiput(246.91,52.57)(-.031392,.074967){14}{\line(0,1){.074967}}
\multiput(246.47,53.62)(-.032193,.068685){15}{\line(0,1){.068685}}
\multiput(245.99,54.65)(-.032842,.063077){16}{\line(0,1){.063077}}
\multiput(245.46,55.66)(-.033361,.058025){17}{\line(0,1){.058025}}
\multiput(244.89,56.65)(-.03199,.050626){19}{\line(0,1){.050626}}
\multiput(244.29,57.61)(-.032371,.046785){20}{\line(0,1){.046785}}
\multiput(243.64,58.55)(-.032662,.043232){21}{\line(0,1){.043232}}
\multiput(242.95,59.45)(-.032872,.039929){22}{\line(0,1){.039929}}
\multiput(242.23,60.33)(-.033009,.036848){23}{\line(0,1){.036848}}
\multiput(241.47,61.18)(-.03308,.033962){24}{\line(0,1){.033962}}
\multiput(240.68,61.99)(-.034469,.032552){24}{\line(-1,0){.034469}}
\multiput(239.85,62.78)(-.037353,.032436){23}{\line(-1,0){.037353}}
\multiput(238.99,63.52)(-.040432,.032252){22}{\line(-1,0){.040432}}
\multiput(238.1,64.23)(-.045917,.03359){20}{\line(-1,0){.045917}}
\multiput(237.18,64.9)(-.049767,.03331){19}{\line(-1,0){.049767}}
\multiput(236.24,65.54)(-.053954,.032938){18}{\line(-1,0){.053954}}
\multiput(235.27,66.13)(-.058533,.032461){17}{\line(-1,0){.058533}}
\multiput(234.27,66.68)(-.063577,.031864){16}{\line(-1,0){.063577}}
\multiput(233.25,67.19)(-.074115,.033352){14}{\line(-1,0){.074115}}
\multiput(232.22,67.66)(-.081246,.032556){13}{\line(-1,0){.081246}}
\multiput(231.16,68.08)(-.08941,.03157){12}{\line(-1,0){.08941}}
\multiput(230.09,68.46)(-.10878,.03337){10}{\line(-1,0){.10878}}
\multiput(229,68.79)(-.12231,.032){9}{\line(-1,0){.12231}}
\multiput(227.9,69.08)(-.13898,.03023){8}{\line(-1,0){.13898}}
\multiput(226.79,69.32)(-.18683,.03253){6}{\line(-1,0){.18683}}
\multiput(225.66,69.52)(-.22563,.02965){5}{\line(-1,0){.22563}}
\put(224.54,69.67){\line(-1,0){1.133}}
\put(223.4,69.77){\line(-1,0){1.137}}
\put(222.27,69.82){\line(-1,0){1.138}}
\put(221.13,69.83){\line(-1,0){1.137}}
\put(219.99,69.79){\line(-1,0){1.134}}
\multiput(218.86,69.7)(-.22594,-.0272){5}{\line(-1,0){.22594}}
\multiput(217.73,69.56)(-.18717,-.0305){6}{\line(-1,0){.18717}}
\multiput(216.6,69.38)(-.1592,-.03282){7}{\line(-1,0){.1592}}
\multiput(215.49,69.15)(-.12265,-.03067){9}{\line(-1,0){.12265}}
\multiput(214.39,68.87)(-.10914,-.03218){10}{\line(-1,0){.10914}}
\multiput(213.3,68.55)(-.09791,-.03337){11}{\line(-1,0){.09791}}
\multiput(212.22,68.18)(-.081594,-.031672){13}{\line(-1,0){.081594}}
\multiput(211.16,67.77)(-.074473,-.032546){14}{\line(-1,0){.074473}}
\multiput(210.11,67.32)(-.06818,-.03325){15}{\line(-1,0){.06818}}
\multiput(209.09,66.82)(-.058882,-.031824){17}{\line(-1,0){.058882}}
\multiput(208.09,66.28)(-.054308,-.03235){18}{\line(-1,0){.054308}}
\multiput(207.11,65.69)(-.050126,-.032768){19}{\line(-1,0){.050126}}
\multiput(206.16,65.07)(-.046279,-.03309){20}{\line(-1,0){.046279}}
\multiput(205.24,64.41)(-.042722,-.033326){21}{\line(-1,0){.042722}}
\multiput(204.34,63.71)(-.039417,-.033485){22}{\line(-1,0){.039417}}
\multiput(203.47,62.97)(-.036334,-.033575){23}{\line(-1,0){.036334}}
\multiput(202.64,62.2)(-.033447,-.033601){24}{\line(0,-1){.033601}}
\multiput(201.83,61.39)(-.033407,-.036488){23}{\line(0,-1){.036488}}
\multiput(201.06,60.55)(-.033304,-.03957){22}{\line(0,-1){.03957}}
\multiput(200.33,59.68)(-.033129,-.042874){21}{\line(0,-1){.042874}}
\multiput(199.64,58.78)(-.032877,-.046431){20}{\line(0,-1){.046431}}
\multiput(198.98,57.86)(-.032537,-.050276){19}{\line(0,-1){.050276}}
\multiput(198.36,56.9)(-.0321,-.054456){18}{\line(0,-1){.054456}}
\multiput(197.78,55.92)(-.033525,-.062717){16}{\line(0,-1){.062717}}
\multiput(197.25,54.92)(-.032937,-.068332){15}{\line(0,-1){.068332}}
\multiput(196.75,53.89)(-.032204,-.074622){14}{\line(0,-1){.074622}}
\multiput(196.3,52.85)(-.031297,-.081739){13}{\line(0,-1){.081739}}
\multiput(195.89,51.78)(-.03292,-.09806){11}{\line(0,-1){.09806}}
\multiput(195.53,50.71)(-.03168,-.10928){10}{\line(0,-1){.10928}}
\multiput(195.22,49.61)(-.03011,-.12279){9}{\line(0,-1){.12279}}
\multiput(194.94,48.51)(-.03209,-.15935){7}{\line(0,-1){.15935}}
\multiput(194.72,47.39)(-.02964,-.18731){6}{\line(0,-1){.18731}}
\multiput(194.54,46.27)(-.0327,-.2826){4}{\line(0,-1){.2826}}
\put(194.41,45.14){\line(0,-1){1.135}}
\put(194.33,44){\line(0,-1){1.137}}
\put(194.29,42.87){\line(0,-1){1.138}}
\put(194.3,41.73){\line(0,-1){1.136}}
\put(194.36,40.59){\line(0,-1){1.133}}
\multiput(194.47,39.46)(.03068,-.22549){5}{\line(0,-1){.22549}}
\multiput(194.62,38.33)(.03339,-.18668){6}{\line(0,-1){.18668}}
\multiput(194.82,37.21)(.03086,-.13884){8}{\line(0,-1){.13884}}
\multiput(195.07,36.1)(.03256,-.12216){9}{\line(0,-1){.12216}}
\multiput(195.36,35)(.03079,-.09875){11}{\line(0,-1){.09875}}
\multiput(195.7,33.91)(.03198,-.08927){12}{\line(0,-1){.08927}}
\multiput(196.08,32.84)(.032929,-.081096){13}{\line(0,-1){.081096}}
\multiput(196.51,31.79)(.033692,-.073961){14}{\line(0,-1){.073961}}
\multiput(196.98,30.75)(.032156,-.06343){16}{\line(0,-1){.06343}}
\multiput(197.5,29.74)(.032729,-.058384){17}{\line(0,-1){.058384}}
\multiput(198.05,28.75)(.033185,-.053802){18}{\line(0,-1){.053802}}
\multiput(198.65,27.78)(.033538,-.049614){19}{\line(0,-1){.049614}}
\multiput(199.29,26.84)(.032191,-.043583){21}{\line(0,-1){.043583}}
\multiput(199.97,25.92)(.032437,-.040284){22}{\line(0,-1){.040284}}
\multiput(200.68,25.03)(.032608,-.037204){23}{\line(0,-1){.037204}}
\multiput(201.43,24.18)(.03271,-.034319){24}{\line(0,-1){.034319}}
\multiput(202.21,23.35)(.034113,-.032924){24}{\line(1,0){.034113}}
\multiput(203.03,22.56)(.036999,-.03284){23}{\line(1,0){.036999}}
\multiput(203.88,21.81)(.04008,-.032689){22}{\line(1,0){.04008}}
\multiput(204.77,21.09)(.043381,-.032463){21}{\line(1,0){.043381}}
\multiput(205.68,20.41)(.046933,-.032156){20}{\line(1,0){.046933}}
\multiput(206.62,19.77)(.053593,-.033521){18}{\line(1,0){.053593}}
\multiput(207.58,19.16)(.058177,-.033094){17}{\line(1,0){.058177}}
\multiput(208.57,18.6)(.063227,-.032552){16}{\line(1,0){.063227}}
\multiput(209.58,18.08)(.068832,-.031878){15}{\line(1,0){.068832}}
\multiput(210.61,17.6)(.080888,-.033436){13}{\line(1,0){.080888}}
\multiput(211.66,17.17)(.08906,-.03253){12}{\line(1,0){.08906}}
\multiput(212.73,16.78)(.09856,-.03141){11}{\line(1,0){.09856}}
\multiput(213.82,16.43)(.12196,-.03332){9}{\line(1,0){.12196}}
\multiput(214.91,16.13)(.13864,-.03173){8}{\line(1,0){.13864}}
\multiput(216.02,15.88)(.15983,-.02962){7}{\line(1,0){.15983}}
\multiput(217.14,15.67)(.22529,-.0321){5}{\line(1,0){.22529}}
\put(218.27,15.51){\line(1,0){1.132}}
\put(219.4,15.39){\line(1,0){1.136}}
\put(220.54,15.33){\line(1,0){1.138}}
\put(221.68,15.31){\line(1,0){1.137}}
\put(222.81,15.34){\line(1,0){1.135}}
\multiput(223.95,15.42)(.2828,.0309){4}{\line(1,0){.2828}}
\multiput(225.08,15.54)(.18749,.02847){6}{\line(1,0){.18749}}
\multiput(226.2,15.71)(.15955,.03109){7}{\line(1,0){.15955}}
\multiput(227.32,15.93)(.13835,.033){8}{\line(1,0){.13835}}
\multiput(228.43,16.19)(.10948,.031){10}{\line(1,0){.10948}}
\multiput(229.52,16.5)(.09826,.03231){11}{\line(1,0){.09826}}
\multiput(230.6,16.86)(.08876,.03335){12}{\line(1,0){.08876}}
\multiput(231.67,17.26)(.074822,.031736){14}{\line(1,0){.074822}}
\multiput(232.72,17.7)(.068537,.032509){15}{\line(1,0){.068537}}
\multiput(233.74,18.19)(.062925,.033132){16}{\line(1,0){.062925}}
\multiput(234.75,18.72)(.057871,.033627){17}{\line(1,0){.057871}}
\multiput(235.73,19.29)(.050479,.032222){19}{\line(1,0){.050479}}
\multiput(236.69,19.9)(.046636,.032585){20}{\line(1,0){.046636}}
\multiput(237.63,20.56)(.043081,.03286){21}{\line(1,0){.043081}}
\multiput(238.53,21.25)(.039778,.033055){22}{\line(1,0){.039778}}
\multiput(239.41,21.97)(.036696,.033178){23}{\line(1,0){.036696}}
\multiput(240.25,22.74)(.033809,.033236){24}{\line(1,0){.033809}}
\multiput(241.06,23.53)(.032393,.034618){24}{\line(0,1){.034618}}
\multiput(241.84,24.36)(.033731,.039207){22}{\line(0,1){.039207}}
\multiput(242.58,25.23)(.033593,.042512){21}{\line(0,1){.042512}}
\multiput(243.29,26.12)(.033379,.046071){20}{\line(0,1){.046071}}
\multiput(243.95,27.04)(.033081,.04992){19}{\line(0,1){.04992}}
\multiput(244.58,27.99)(.03269,.054105){18}{\line(0,1){.054105}}
\multiput(245.17,28.96)(.032192,.058682){17}{\line(0,1){.058682}}
\multiput(245.72,29.96)(.033677,.06797){15}{\line(0,1){.06797}}
\multiput(246.22,30.98)(.033012,.074268){14}{\line(0,1){.074268}}
\multiput(246.69,32.02)(.032182,.081395){13}{\line(0,1){.081395}}
\multiput(247.1,33.08)(.03115,.08956){12}{\line(0,1){.08956}}
\multiput(247.48,34.15)(.03287,.10893){10}{\line(0,1){.10893}}
\multiput(247.81,35.24)(.03144,.12246){9}{\line(0,1){.12246}}
\multiput(248.09,36.34)(.02959,.13912){8}{\line(0,1){.13912}}
\multiput(248.33,37.46)(.03167,.18698){6}{\line(0,1){.18698}}
\multiput(248.52,38.58)(.02861,.22576){5}{\line(0,1){.22576}}
\put(248.66,39.71){\line(0,1){1.134}}
\put(248.76,40.84){\line(0,1){1.728}}
\put(48,27.88){\makebox(0,0)[cc]{$\sigma\theta$}}
\put(218.5,28.06){\makebox(0,0)[cc]{$\sigma\theta$}}
\put(127.99,35.96){\makebox(0,0)[cc]{$\breve{\q}_{\sigma
\theta}$}}
\put(70.64,56.8){\line(1,0){.999}}
\put(72.64,56.79){\line(1,0){.999}}
\put(74.64,56.78){\line(1,0){.999}}
\put(76.64,56.77){\line(1,0){.999}}
\put(78.64,56.76){\line(1,0){.999}}
\put(80.63,56.74){\line(1,0){.999}}
\put(82.63,56.73){\line(1,0){.999}}
\put(84.63,56.72){\line(1,0){.999}}
\put(86.63,56.71){\line(1,0){.999}}
\put(88.63,56.7){\line(1,0){.999}}
\put(90.63,56.69){\line(1,0){.999}}
\put(92.63,56.68){\line(1,0){.999}}
\put(94.63,56.67){\line(1,0){.999}}
\put(96.63,56.66){\line(1,0){.999}}
\put(98.62,56.65){\line(1,0){.999}}
\put(100.62,56.63){\line(1,0){.999}}
\put(102.62,56.62){\line(1,0){.999}}
\put(104.62,56.61){\line(1,0){.999}}
\put(106.62,56.6){\line(1,0){.999}}
\put(108.62,56.59){\line(1,0){.999}}
\put(110.62,56.58){\line(1,0){.999}}
\put(112.62,56.57){\line(1,0){.999}}
\put(114.62,56.56){\line(1,0){.999}}
\put(116.61,56.55){\line(1,0){.999}}
\put(118.61,56.54){\line(1,0){.999}}
\put(120.61,56.52){\line(1,0){.999}}
\put(122.61,56.51){\line(1,0){.999}}
\put(124.61,56.5){\line(1,0){.999}}
\put(126.61,56.49){\line(1,0){.999}}
\put(128.61,56.48){\line(1,0){.999}}
\put(130.61,56.47){\line(1,0){.999}}
\put(132.61,56.46){\line(1,0){.999}}
\put(134.6,56.45){\line(1,0){.999}}
\put(136.6,56.44){\line(1,0){.999}}
\put(138.6,56.42){\line(1,0){.999}}
\put(140.6,56.41){\line(1,0){.999}}
\put(142.6,56.4){\line(1,0){.999}}
\put(144.6,56.39){\line(1,0){.999}}
\put(146.6,56.38){\line(1,0){.999}}
\put(148.6,56.37){\line(1,0){.999}}
\put(150.6,56.36){\line(1,0){.999}}
\put(152.59,56.35){\line(1,0){.999}}
\put(154.59,56.34){\line(1,0){.999}}
\put(156.59,56.33){\line(1,0){.999}}
\put(158.59,56.31){\line(1,0){.999}}
\put(160.59,56.3){\line(1,0){.999}}
\put(162.59,56.29){\line(1,0){.999}}
\put(164.59,56.28){\line(1,0){.999}}
\put(166.59,56.27){\line(1,0){.999}}
\put(168.59,56.26){\line(1,0){.999}}
\put(170.58,56.25){\line(1,0){.999}}
\put(172.58,56.24){\line(1,0){.999}}
\put(174.58,56.23){\line(1,0){.999}}
\put(176.58,56.22){\line(1,0){.999}}
\put(178.58,56.2){\line(1,0){.999}}
\put(180.58,56.19){\line(1,0){.999}}
\put(182.58,56.18){\line(1,0){.999}}
\put(184.58,56.17){\line(1,0){.999}}
\put(186.58,56.16){\line(1,0){.999}}
\put(188.57,56.15){\line(1,0){.999}}
\put(190.57,56.14){\line(1,0){.999}}
\put(192.57,56.13){\line(1,0){.999}}
\put(194.57,56.12){\line(1,0){.999}}
\put(196.57,56.11){\line(1,0){.999}}
\put(70.46,28.52){\line(1,0){.995}}
\put(72.45,28.55){\line(1,0){.995}}
\put(74.44,28.58){\line(1,0){.995}}
\put(76.43,28.6){\line(1,0){.995}}
\put(78.42,28.63){\line(1,0){.995}}
\put(80.41,28.66){\line(1,0){.995}}
\put(82.4,28.69){\line(1,0){.995}}
\put(84.39,28.71){\line(1,0){.995}}
\put(86.38,28.74){\line(1,0){.995}}
\put(88.37,28.77){\line(1,0){.995}}
\put(90.37,28.8){\line(1,0){.995}}
\put(92.36,28.83){\line(1,0){.995}}
\put(94.35,28.85){\line(1,0){.995}}
\put(96.34,28.88){\line(1,0){.995}}
\put(98.33,28.91){\line(1,0){.995}}
\put(100.32,28.94){\line(1,0){.995}}
\put(102.31,28.97){\line(1,0){.995}}
\put(104.3,28.99){\line(1,0){.995}}
\put(106.29,29.02){\line(1,0){.995}}
\put(108.28,29.05){\line(1,0){.995}}
\put(110.27,29.08){\line(1,0){.995}}
\put(112.26,29.11){\line(1,0){.995}}
\put(114.25,29.13){\line(1,0){.995}}
\put(116.24,29.16){\line(1,0){.995}}
\put(118.23,29.19){\line(1,0){.995}}
\put(120.22,29.22){\line(1,0){.995}}
\put(122.21,29.24){\line(1,0){.995}}
\put(124.2,29.27){\line(1,0){.995}}
\put(126.2,29.3){\line(1,0){.995}}
\put(128.19,29.33){\line(1,0){.995}}
\put(130.18,29.36){\line(1,0){.995}}
\put(132.17,29.38){\line(1,0){.995}}
\put(134.16,29.41){\line(1,0){.995}}
\put(136.15,29.44){\line(1,0){.995}}
\put(138.14,29.47){\line(1,0){.995}}
\put(140.13,29.5){\line(1,0){.995}}
\put(142.12,29.52){\line(1,0){.995}}
\put(144.11,29.55){\line(1,0){.995}}
\put(146.1,29.58){\line(1,0){.995}}
\put(148.09,29.61){\line(1,0){.995}}
\put(150.08,29.63){\line(1,0){.995}}
\put(152.07,29.66){\line(1,0){.995}}
\put(154.06,29.69){\line(1,0){.995}}
\put(156.05,29.72){\line(1,0){.995}}
\put(158.04,29.75){\line(1,0){.995}}
\put(160.03,29.77){\line(1,0){.995}}
\put(162.03,29.8){\line(1,0){.995}}
\put(164.02,29.83){\line(1,0){.995}}
\put(166.01,29.86){\line(1,0){.995}}
\put(168,29.89){\line(1,0){.995}}
\put(169.99,29.91){\line(1,0){.995}}
\put(171.98,29.94){\line(1,0){.995}}
\put(173.97,29.97){\line(1,0){.995}}
\put(175.96,30){\line(1,0){.995}}
\put(177.95,30.02){\line(1,0){.995}}
\put(179.94,30.05){\line(1,0){.995}}
\put(181.93,30.08){\line(1,0){.995}}
\put(183.92,30.11){\line(1,0){.995}}
\put(185.91,30.14){\line(1,0){.995}}
\put(187.9,30.16){\line(1,0){.995}}
\put(189.89,30.19){\line(1,0){.995}}
\put(191.88,30.22){\line(1,0){.995}}
\put(193.87,30.25){\line(1,0){.995}}
\put(195.86,30.28){\line(1,0){.995}}
\put(150.26,42.91){\line(0,1){13.43}}
\put(156.45,49.8){\makebox(0,0)[cc]{$\theta$}} \thicklines
\qbezier(47.2,43.08)(68.94,64.56)(74.42,113.27)
\qbezier(74.42,113.27)(88.74,112.74)(94.58,123.52)
\qbezier(94.58,123.52)(129.58,-1.55)(171.65,123.34)
\qbezier(171.65,123.34)(179.78,109.46)(194.98,113.62)
\qbezier(194.98,113.62)(191.8,57.23)(221.85,42.55)
\put(131,77.2){\makebox(0,0)[cc]{$\pp\setminus\q$}}
\put(41.32,43.7){\makebox(0,0)[cc]{$x$}}
\put(227.73,42.81){\makebox(0,0)[cc]{$y$}} \thinlines
\multiput(47.2,42.74)(.033727,-.044345){17}{\line(0,-1){.044345}}
\multiput(48.35,41.23)(.033727,-.044345){17}{\line(0,-1){.044345}}
\multiput(49.49,39.73)(.033727,-.044345){17}{\line(0,-1){.044345}}
\multiput(50.64,38.22)(.033727,-.044345){17}{\line(0,-1){.044345}}
\multiput(51.79,36.71)(.033727,-.044345){17}{\line(0,-1){.044345}}
\multiput(52.93,35.2)(.033727,-.044345){17}{\line(0,-1){.044345}}
\multiput(54.08,33.69)(.033727,-.044345){17}{\line(0,-1){.044345}}
\multiput(55.23,32.19)(.033727,-.044345){17}{\line(0,-1){.044345}}
\multiput(56.37,30.68)(.033727,-.044345){17}{\line(0,-1){.044345}}
\multiput(57.52,29.17)(.033727,-.044345){17}{\line(0,-1){.044345}}
\multiput(58.67,27.66)(.033727,-.044345){17}{\line(0,-1){.044345}}
\multiput(59.81,26.16)(.033727,-.044345){17}{\line(0,-1){.044345}}
\multiput(60.96,24.65)(.033727,-.044345){17}{\line(0,-1){.044345}}
\multiput(62.11,23.14)(.033727,-.044345){17}{\line(0,-1){.044345}}
\multiput(221.12,42.74)(-.033033,-.041882){18}{\line(0,-1){.041882}}
\multiput(219.93,41.23)(-.033033,-.041882){18}{\line(0,-1){.041882}}
\multiput(218.74,39.73)(-.033033,-.041882){18}{\line(0,-1){.041882}}
\multiput(217.55,38.22)(-.033033,-.041882){18}{\line(0,-1){.041882}}
\multiput(216.36,36.71)(-.033033,-.041882){18}{\line(0,-1){.041882}}
\multiput(215.18,35.2)(-.033033,-.041882){18}{\line(0,-1){.041882}}
\multiput(213.99,33.69)(-.033033,-.041882){18}{\line(0,-1){.041882}}
\multiput(212.8,32.19)(-.033033,-.041882){18}{\line(0,-1){.041882}}
\multiput(211.61,30.68)(-.033033,-.041882){18}{\line(0,-1){.041882}}
\multiput(210.42,29.17)(-.033033,-.041882){18}{\line(0,-1){.041882}}
\multiput(209.23,27.66)(-.033033,-.041882){18}{\line(0,-1){.041882}}
\multiput(208.04,26.16)(-.033033,-.041882){18}{\line(0,-1){.041882}}
\multiput(206.85,24.65)(-.033033,-.041882){18}{\line(0,-1){.041882}}
\multiput(205.66,23.14)(-.033033,-.041882){18}{\line(0,-1){.041882}}
\end{picture}
\unitlength .32mm 
\linethickness{0.4pt}
\ifx\plotpoint\undefined\newsavebox{\plotpoint}\fi 
\begin{picture}(100.75,125.5)(20,0)
\multiput(119.75,33.75)(-.0473163842,.0337217514){1416}{\line(-1,0){.0473163842}}
\multiput(119.5,34)(.0342055485,.0337326608){1586}{\line(1,0){.0342055485}}
\thicklines \qbezier(53.25,81)(73.88,91.75)(75,124.5)
\qbezier(75,124.5)(100.88,122)(106.25,142.5)
\qbezier(106.25,142.5)(122.25,122.63)(147.25,132.25)
\qbezier(147.25,132.25)(138.13,97.88)(173.5,87)
\put(160.25,114.5){\makebox(0,0)[cc]{$\oo_x$}}
\put(121,29){\makebox(0,0)[cc]{$x$}} \thinlines
\put(148.82,34){\line(0,1){1.191}}
\put(148.79,35.19){\line(0,1){1.189}}
\multiput(148.72,36.38)(-.0302,.2963){4}{\line(0,1){.2963}}
\multiput(148.6,37.56)(-.02814,.19653){6}{\line(0,1){.19653}}
\multiput(148.43,38.74)(-.03094,.16733){7}{\line(0,1){.16733}}
\multiput(148.22,39.92)(-.033,.14519){8}{\line(0,1){.14519}}
\multiput(147.95,41.08)(-.0311,.11499){10}{\line(0,1){.11499}}
\multiput(147.64,42.23)(-.03249,.1033){11}{\line(0,1){.1033}}
\multiput(147.28,43.36)(-.03361,.0934){12}{\line(0,1){.0934}}
\multiput(146.88,44.48)(-.032033,.078824){14}{\line(0,1){.078824}}
\multiput(146.43,45.59)(-.032861,.072294){15}{\line(0,1){.072294}}
\multiput(145.94,46.67)(-.033535,.066468){16}{\line(0,1){.066468}}
\multiput(145.4,47.74)(-.032184,.057823){18}{\line(0,1){.057823}}
\multiput(144.82,48.78)(-.032691,.053496){19}{\line(0,1){.053496}}
\multiput(144.2,49.79)(-.033095,.049518){20}{\line(0,1){.049518}}
\multiput(143.54,50.78)(-.033409,.04584){21}{\line(0,1){.04584}}
\multiput(142.84,51.75)(-.033641,.042425){22}{\line(0,1){.042425}}
\multiput(142.1,52.68)(-.032392,.037605){24}{\line(0,1){.037605}}
\multiput(141.32,53.58)(-.032537,.034808){25}{\line(0,1){.034808}}
\multiput(140.51,54.45)(-.033924,.033457){25}{\line(-1,0){.033924}}
\multiput(139.66,55.29)(-.036724,.033387){24}{\line(-1,0){.036724}}
\multiput(138.78,56.09)(-.039704,.033254){23}{\line(-1,0){.039704}}
\multiput(137.86,56.85)(-.042887,.03305){22}{\line(-1,0){.042887}}
\multiput(136.92,57.58)(-.046299,.03277){21}{\line(-1,0){.046299}}
\multiput(135.95,58.27)(-.049971,.032406){20}{\line(-1,0){.049971}}
\multiput(134.95,58.92)(-.05694,.033721){18}{\line(-1,0){.05694}}
\multiput(133.92,59.52)(-.06169,.033227){17}{\line(-1,0){.06169}}
\multiput(132.88,60.09)(-.066926,.032612){16}{\line(-1,0){.066926}}
\multiput(131.8,60.61)(-.072742,.031857){15}{\line(-1,0){.072742}}
\multiput(130.71,61.09)(-.085357,.033318){13}{\line(-1,0){.085357}}
\multiput(129.6,61.52)(-.09386,.03231){12}{\line(-1,0){.09386}}
\multiput(128.48,61.91)(-.10374,.03106){11}{\line(-1,0){.10374}}
\multiput(127.34,62.25)(-.12823,.03278){9}{\line(-1,0){.12823}}
\multiput(126.18,62.55)(-.14564,.03098){8}{\line(-1,0){.14564}}
\multiput(125.02,62.79)(-.1957,.03339){6}{\line(-1,0){.1957}}
\multiput(123.84,62.99)(-.23628,.0305){5}{\line(-1,0){.23628}}
\put(122.66,63.15){\line(-1,0){1.187}}
\put(121.47,63.25){\line(-1,0){1.19}}
\put(120.29,63.31){\line(-1,0){1.191}}
\put(119.09,63.32){\line(-1,0){1.19}}
\put(117.9,63.27){\line(-1,0){1.188}}
\multiput(116.72,63.19)(-.23665,-.02744){5}{\line(-1,0){.23665}}
\multiput(115.53,63.05)(-.19612,-.03086){6}{\line(-1,0){.19612}}
\multiput(114.36,62.86)(-.16689,-.03326){7}{\line(-1,0){.16689}}
\multiput(113.19,62.63)(-.12864,-.03112){9}{\line(-1,0){.12864}}
\multiput(112.03,62.35)(-.11455,-.03268){10}{\line(-1,0){.11455}}
\multiput(110.88,62.02)(-.09427,-.03109){12}{\line(-1,0){.09427}}
\multiput(109.75,61.65)(-.08578,-.032212){13}{\line(-1,0){.08578}}
\multiput(108.64,61.23)(-.078373,-.033122){14}{\line(-1,0){.078373}}
\multiput(107.54,60.77)(-.067342,-.031743){16}{\line(-1,0){.067342}}
\multiput(106.46,60.26)(-.062115,-.032426){17}{\line(-1,0){.062115}}
\multiput(105.41,59.71)(-.057372,-.032982){18}{\line(-1,0){.057372}}
\multiput(104.37,59.12)(-.053038,-.033428){19}{\line(-1,0){.053038}}
\multiput(103.37,58.48)(-.046719,-.032169){21}{\line(-1,0){.046719}}
\multiput(102.39,57.8)(-.043311,-.032493){22}{\line(-1,0){.043311}}
\multiput(101.43,57.09)(-.040131,-.032737){23}{\line(-1,0){.040131}}
\multiput(100.51,56.34)(-.037153,-.03291){24}{\line(-1,0){.037153}}
\multiput(99.62,55.55)(-.034354,-.033016){25}{\line(-1,0){.034354}}
\multiput(98.76,54.72)(-.032984,-.034384){25}{\line(0,-1){.034384}}
\multiput(97.93,53.86)(-.032876,-.037183){24}{\line(0,-1){.037183}}
\multiput(97.15,52.97)(-.032701,-.040161){23}{\line(0,-1){.040161}}
\multiput(96.39,52.05)(-.032453,-.043341){22}{\line(0,-1){.043341}}
\multiput(95.68,51.09)(-.033732,-.049085){20}{\line(0,-1){.049085}}
\multiput(95,50.11)(-.03338,-.053069){19}{\line(0,-1){.053069}}
\multiput(94.37,49.1)(-.03293,-.057402){18}{\line(0,-1){.057402}}
\multiput(93.78,48.07)(-.032369,-.062145){17}{\line(0,-1){.062145}}
\multiput(93.23,47.01)(-.031682,-.067371){16}{\line(0,-1){.067371}}
\multiput(92.72,45.93)(-.03305,-.078403){14}{\line(0,-1){.078403}}
\multiput(92.26,44.84)(-.032133,-.08581){13}{\line(0,-1){.08581}}
\multiput(91.84,43.72)(-.03101,-.0943){12}{\line(0,-1){.0943}}
\multiput(91.47,42.59)(-.03258,-.11458){10}{\line(0,-1){.11458}}
\multiput(91.14,41.44)(-.031,-.12867){9}{\line(0,-1){.12867}}
\multiput(90.86,40.29)(-.0331,-.16692){7}{\line(0,-1){.16692}}
\multiput(90.63,39.12)(-.03068,-.19615){6}{\line(0,-1){.19615}}
\multiput(90.45,37.94)(-.02722,-.23668){5}{\line(0,-1){.23668}}
\put(90.31,36.76){\line(0,-1){1.188}}
\put(90.22,35.57){\line(0,-1){2.382}}
\put(90.19,33.19){\line(0,-1){1.19}}
\put(90.25,32){\line(0,-1){1.187}}
\multiput(90.36,30.81)(.03071,-.23625){5}{\line(0,-1){.23625}}
\multiput(90.51,29.63)(.03357,-.19567){6}{\line(0,-1){.19567}}
\multiput(90.71,28.46)(.03112,-.14561){8}{\line(0,-1){.14561}}
\multiput(90.96,27.29)(.0329,-.1282){9}{\line(0,-1){.1282}}
\multiput(91.26,26.14)(.03115,-.10371){11}{\line(0,-1){.10371}}
\multiput(91.6,25)(.03239,-.09383){12}{\line(0,-1){.09383}}
\multiput(91.99,23.87)(.033397,-.085326){13}{\line(0,-1){.085326}}
\multiput(92.42,22.76)(.031924,-.072713){15}{\line(0,-1){.072713}}
\multiput(92.9,21.67)(.032673,-.066896){16}{\line(0,-1){.066896}}
\multiput(93.42,20.6)(.033283,-.06166){17}{\line(0,-1){.06166}}
\multiput(93.99,19.55)(.031996,-.053914){19}{\line(0,-1){.053914}}
\multiput(94.6,18.53)(.032452,-.049941){20}{\line(0,-1){.049941}}
\multiput(95.25,17.53)(.032813,-.046269){21}{\line(0,-1){.046269}}
\multiput(95.93,16.56)(.03309,-.042857){22}{\line(0,-1){.042857}}
\multiput(96.66,15.61)(.03329,-.039674){23}{\line(0,-1){.039674}}
\multiput(97.43,14.7)(.033421,-.036693){24}{\line(0,-1){.036693}}
\multiput(98.23,13.82)(.033489,-.033893){25}{\line(0,-1){.033893}}
\multiput(99.07,12.97)(.034838,-.032505){25}{\line(1,0){.034838}}
\multiput(99.94,12.16)(.037635,-.032357){24}{\line(1,0){.037635}}
\multiput(100.84,11.39)(.042456,-.033602){22}{\line(1,0){.042456}}
\multiput(101.78,10.65)(.045871,-.033366){21}{\line(1,0){.045871}}
\multiput(102.74,9.95)(.049548,-.033049){20}{\line(1,0){.049548}}
\multiput(103.73,9.28)(.053526,-.032642){19}{\line(1,0){.053526}}
\multiput(104.75,8.66)(.057852,-.032131){18}{\line(1,0){.057852}}
\multiput(105.79,8.09)(.066499,-.033474){16}{\line(1,0){.066499}}
\multiput(106.85,7.55)(.072324,-.032795){15}{\line(1,0){.072324}}
\multiput(107.94,7.06)(.078853,-.031961){14}{\line(1,0){.078853}}
\multiput(109.04,6.61)(.09343,-.03352){12}{\line(1,0){.09343}}
\multiput(110.16,6.21)(.10333,-.0324){11}{\line(1,0){.10333}}
\multiput(111.3,5.85)(.11502,-.03099){10}{\line(1,0){.11502}}
\multiput(112.45,5.54)(.14522,-.03287){8}{\line(1,0){.14522}}
\multiput(113.61,5.28)(.16736,-.03079){7}{\line(1,0){.16736}}
\multiput(114.78,5.06)(.23586,-.03355){5}{\line(1,0){.23586}}
\put(115.96,4.9){\line(1,0){1.185}}
\put(117.15,4.78){\line(1,0){1.189}}
\put(118.34,4.7){\line(1,0){2.382}}
\put(120.72,4.71){\line(1,0){1.189}}
\multiput(121.91,4.78)(.2962,.0305){4}{\line(1,0){.2962}}
\multiput(123.09,4.9)(.1965,.02832){6}{\line(1,0){.1965}}
\multiput(124.27,5.07)(.1673,.03109){7}{\line(1,0){.1673}}
\multiput(125.44,5.29)(.14516,.03313){8}{\line(1,0){.14516}}
\multiput(126.6,5.56)(.11496,.0312){10}{\line(1,0){.11496}}
\multiput(127.75,5.87)(.10327,.03259){11}{\line(1,0){.10327}}
\multiput(128.89,6.23)(.09337,.03369){12}{\line(1,0){.09337}}
\multiput(130.01,6.63)(.078794,.032105){14}{\line(1,0){.078794}}
\multiput(131.11,7.08)(.072263,.032928){15}{\line(1,0){.072263}}
\multiput(132.2,7.57)(.066437,.033596){16}{\line(1,0){.066437}}
\multiput(133.26,8.11)(.057793,.032237){18}{\line(1,0){.057793}}
\multiput(134.3,8.69)(.053466,.03274){19}{\line(1,0){.053466}}
\multiput(135.32,9.31)(.049487,.03314){20}{\line(1,0){.049487}}
\multiput(136.31,9.98)(.04581,.033451){21}{\line(1,0){.04581}}
\multiput(137.27,10.68)(.042394,.03368){22}{\line(1,0){.042394}}
\multiput(138.2,11.42)(.037575,.032426){24}{\line(1,0){.037575}}
\multiput(139.1,12.2)(.034778,.032569){25}{\line(1,0){.034778}}
\multiput(139.97,13.01)(.033426,.033955){25}{\line(0,1){.033955}}
\multiput(140.81,13.86)(.033354,.036755){24}{\line(0,1){.036755}}
\multiput(141.61,14.74)(.033217,.039735){23}{\line(0,1){.039735}}
\multiput(142.37,15.66)(.033011,.042917){22}{\line(0,1){.042917}}
\multiput(143.1,16.6)(.032728,.046329){21}{\line(0,1){.046329}}
\multiput(143.78,17.57)(.03236,.050001){20}{\line(0,1){.050001}}
\multiput(144.43,18.57)(.033669,.056971){18}{\line(0,1){.056971}}
\multiput(145.04,19.6)(.03317,.061721){17}{\line(0,1){.061721}}
\multiput(145.6,20.65)(.03255,.066956){16}{\line(0,1){.066956}}
\multiput(146.12,21.72)(.03179,.072771){15}{\line(0,1){.072771}}
\multiput(146.6,22.81)(.03324,.085387){13}{\line(0,1){.085387}}
\multiput(147.03,23.92)(.03222,.09389){12}{\line(0,1){.09389}}
\multiput(147.42,25.05)(.03096,.10377){11}{\line(0,1){.10377}}
\multiput(147.76,26.19)(.03266,.12826){9}{\line(0,1){.12826}}
\multiput(148.05,27.34)(.03085,.14567){8}{\line(0,1){.14567}}
\multiput(148.3,28.51)(.03321,.19573){6}{\line(0,1){.19573}}
\multiput(148.5,29.68)(.03028,.2363){5}{\line(0,1){.2363}}
\put(148.65,30.87){\line(0,1){1.187}}
\put(148.75,32.05){\line(0,1){1.948}}
\put(119.18,34.18){\line(1,0){.991}}
\put(121.16,34.2){\line(1,0){.991}}
\put(123.15,34.21){\line(1,0){.991}}
\put(125.13,34.23){\line(1,0){.991}}
\put(127.11,34.25){\line(1,0){.991}}
\put(129.09,34.27){\line(1,0){.991}}
\put(131.08,34.28){\line(1,0){.991}}
\put(133.06,34.3){\line(1,0){.991}}
\put(135.04,34.32){\line(1,0){.991}}
\put(137.02,34.33){\line(1,0){.991}}
\put(139.01,34.35){\line(1,0){.991}}
\put(140.99,34.37){\line(1,0){.991}}
\put(142.97,34.39){\line(1,0){.991}}
\put(144.96,34.4){\line(1,0){.991}}
\put(146.94,34.42){\line(1,0){.991}}
\put(135.75,28){\makebox(0,0)[cc]{$\nu\theta$}}
\end{picture}

\begin{definition}[fat polygon]\label{fatpoly}
Let $\theta>0$, $\sigma \geq 1$ and $\nu \geq 4\sigma$. We call a
$k$-gon $P$ with quasi-geodesic edges \textit{$(\theta , \sigma,
\nu )$--fat} if the following properties hold:
\begin{enumerate}
  \item[$(F_1)$](\textbf{large inscribed radii in interior points, large comparison
angles}) for every edge $\q $ we have, with the notation
\ref{breve}, that
$$\dist \left( \breve{\q}_{\sigma \theta}\, ,\, P \setminus
\q \right) \geq \theta;$$

  \item[$(F_2)$]
(\textbf{large inscribed radii in vertices, large edges}) for
every vertex $x$ we have that $$\dist(x, \oo_x)\geq \nu \theta.$$
\end{enumerate}

When $\sigma =2$ we say that $P$ is \textit{$(\theta , \nu
)$--fat}.
\end{definition}

\begin{lemma}\label{latop}
Let $P$ be a polygon $(\theta , \sigma , \nu )$--fat for some
$\theta >0$, $\sigma \geq 1$ and $\nu \geq 4 \sigma$. Then any two
edges of $P$ without a common vertex are at distance at least
$\theta$ from each other.
\end{lemma}

\proof Let $\q$ and $\q'$ be two edges without a common vertex.
Assume that there exists a point $a\in \q$ such that $\dist
(a,\q') <\theta$. Property $(F_1)$ implies that $a\in
\nn_{\sigma\theta } \left( \{ x,y\} \right)$, where $x,y$ are the
endpoints of $\q$. Property $(F_2)$ implies that $\dist \left( \{
x,y\} , \q' \right) \geq \nu \theta$. Therefore $\dist (a,
\q')\geq (\nu - \sigma ) \theta \geq 3\sigma \theta \geq 3\theta$.
This contradicts the assumption that $\dist (a,\q')
<\theta$.\endproof

\medskip

The following lemma describes a situation in which given two
consecutive edges of a geodesic polygon, any two points on each of
these edges which are at distance at least $2\theta$ from the
common vertex are at distance at least $\theta$ from one another.

\begin{lemma}\label{simplif}
Let $P$ be a geodesic polygon with two consecutive edges $[x,y]$
and $[y,z]$ such that $\dist (x, [y,z])=\dist (x,y)$. Then both
the distance from $[x,y]\setminus B(y,2\theta )$ to $[y,z]$, and
the distance from $[y,z]\setminus B(y,2\theta )$ to $[x,y]$ are at
least $\theta$.
\end{lemma}

\proof The distance from $[x,y]\setminus B(y,2\theta )$ to $[y,z]$
is $2\theta $ because of the hypothesis that $\dist
(x,[y,z])=\dist (x,y)$.

Assume that there exists $p\in [y,z]\setminus B(y,2\theta )$ and
$p'\in [x,y]$ such that $\dist (p,p')<\theta$. Then $\dist
(y,p')\geq \dist (y, p)-\dist (p,p') > \theta > \dist (p,p')$. It
follows that
 $\dist (x,p )\leq \dist (x,p')+\dist (p',p) < \dist (x,p')+\dist (p',y)=\dist (x,y
 )$. This contradicts the fact that $\dist
(x,[y,z])=\dist (x,y)$. \endproof

\begin{theorem}[\cite{DrutuSapir:TreeGraded}, Theorem 4.1 and Remark 4.2, (3)]\label{tgi}
Let $(X,\dist)$ be a geodesic metric space and let $\aaa=\{ A_i
\mid i\in I \}$ be a collection of subsets of $X$. The metric
space $X$ is asymptotically tree-graded with respect to $\aaa$ if
and only if the following properties are satisfied:
\begin{itemize}
  \item[$(\alpha_1)$] For every $\delta >0$ the diameters
of the intersections $\nn_{\delta}(A_i)\cap \nn_{\delta}(A_j)$ are
uniformly bounded for all $i\ne j$.

\item[$(\alpha_2)$] There exists $\varepsilon$ in $\left[ 0, \frac{1}{2}
\right)$ and $M>0$ such that for every geodesic $\g$ of length
$\ell$ and every $A\in \aaa$ with $\g(0),\g(\ell)\in
\nn_{\varepsilon \ell}(A)$ we have that $\g([0, \ell ])\cap
\nn_{M}(A)\neq \emptyset$.

\item[$(\alpha_3)$] For every $k\ge 2$ there exist $\theta>0$,
$\nu \geq 8 $ and $\chi
>0$ such that every $k$-gon $P$ in $X$ with
geodesic edges which is $(\theta , \nu )$--fat satisfies $P
\subset \nn_{\chi}(A)$ for some $A\in \aaa$.
\end{itemize}
\end{theorem}

\begin{rmks}[\cite{DrutuSapir:TreeGraded}, Theorem 4.1 and Remark 4.2]\label{rstr}
\begin{itemize}
    \item[(1)] Property $(\alpha_2)$ from Theorem \ref{tgi} is a slight modification of
     the similar property appearing in Theorem 4.1 in
    \cite{DrutuSapir:TreeGraded}. Nevertheless it
    implies property $\left( \alpha^\varepsilon_2 \right)$ from \cite[Remark 4.2, (3)]{DrutuSapir:TreeGraded},
    which accounts for the accuracy of the modified statement.
    \item[(2)] As a necessary condition, $(\alpha_2)$ can be
    strengthened to ``for every $\varepsilon$ from
    $\left[0,\frac{1}{2}\right)$ there exists $M>0$ such that
    etc.''
\end{itemize}
\end{rmks}

\begin{n}\label{ndiam}
We denote by $\diam_\delta$ an uniform bound provided by property
$(\alpha_1)$ for an arbitrary $\delta \geq 0$.
\end{n}

\begin{rmks}[on the condition that pieces cover the space]\label{ratg}
\begin{itemize}
    \item[(1)] If in property ($T_2$) of Definition \ref{deftgr}
    of tree-graded spaces we allow for trivial
geodesic triangles, that is if we ask that pieces cover a
tree-graded space, then in Theorem \ref{tgi} the following
condition has to be added:
\begin{itemize}
\item[$(\alpha_0)$] there exists $\tau\geq 0$ such that $X=\bigcup_{A\in
\aaa}\nn_\tau(A)$.
\end{itemize}

If $(X,\aaa)$ satisfy only the conditions $(\alpha_1),\,
(\alpha_2),\, (\alpha_3)$ but not $(\alpha_0)$ then it suffices to
add some singletons to $\aaa$ in order to ensure $(\alpha_0)$.
Indeed, for some $\tau >0$ consider in $X\setminus \bigcup_{A\in
\aaa}\nn_\tau(A)$ a maximal subset $\wp$ with the property that
$\dist (p,p')\geq \tau$ for every $p,p'\in \wp$. The space $X$
coincides with $\bigcup_{A\in \aaa}\nn_\tau(A) \cup \bigcup_{p\in
\wp}\nn_\tau(\{p\})$. Properties $(\alpha_1)$ and $(\alpha_2)$ are
obviously satisfied by singletons, whence $X$ is ATG \wrt
$\aaa'=\aaa \cup \left\{ \left\{ p\right\} \; ;\; p\in
\wp\right\}$; moreover $\aaa'$ also satisfies $(\alpha_0)$.

\item[(2)] Let $\hip^3$ be the $3$-dimensional real hyperbolic
space and let $(Hbo_n)_{n\in \N}$ be a countable collection of
pairwise disjoint open horoballs. The complementary set $X_0=
X\setminus \bigsqcup_{n\in \N} Hbo_n$ and the collection of
boundary horospheres $\aaa = \left\{ \partial Hbo_n \; ;\; n\in \N
\right\}$ is the typical example one has in mind when trying to
define relative hyperbolicity for metric spaces. The pair $(X_0,
\aaa)$ does not in general satisfy $(\alpha_0)$, one has to add
singletons to $\aaa$ to ensure that property. In order to remove
this inconvenient, we give up the condition of pieces covering the
space in Definition \ref{deftgr} of tree-graded spaces.
\end{itemize}
\end{rmks}

\begin{rmk}\label{ratg2}
If  $X$ is a metric space ATG with respect to $\aaa$, and a group
$G$ acts $\calk$-transitively (in the sense of Definition
    \ref{tk}
    , with $\calk\geq 0$) by isometries on $X$, $G$
    permuting the subsets in $\aaa$, then property $(\alpha_0)$ is satisfied with $\tau=\calk$.

    It is for instance the case when $X$ is itself a group and
    $\aaa$ is the collection of left cosets of a family of
    subgroups.
\end{rmk}

\subsection{Property $(T_2)$ and polygons with limit
edges}\label{pt2}

Property $(\alpha_3)$ in the definition of a metric space $X$ ATG
\wrt a collection $\aaa$ is used to prove property $(T_2)$ in an
arbitrary asymptotic cone of $X$ \wrt the collection of limit sets
$\aau$. If $X$ is such that any geodesic in an asymptotic cone of
it is a limit geodesic (for instance if $X$ is a CAT(0) metric
space) then it suffices to have $(\alpha_3)$ for $k=6$, that is:
\begin{itemize}
         \item[$(\beta_3)$] there exists $\theta >0$, $\nu \geq 8$ and $\chi>0$
        such that any geodesic hexagon $(\theta , \nu)$--fat is contained in $\nn_\chi (A)$, for some $A\in
        \aaa$.
\end{itemize}

This is due to the following general fact.

\begin{proposition}\label{limtr}
Let $(X,\dist)$ be a geodesic metric space and let $\theta >0$ and
$\nu \geq 8$ be two arbitrary constants. In any asymptotic cone
$\co{X;e,d}$, any simple non-trivial triangle whose edges are
limit geodesics is the limit set $\lio{H_n}$ of a sequence $(H_n)$
of geodesic hexagons that are $(\theta , \nu)$--fat \uas.
\end{proposition}

\proof Consider a non-trivial simple geodesic triangle $\Delta$ in
an asymptotic cone $\co{X;e,d}$, whose edges $[a,b]$, $[b,c]$ and
$[c,a]$ appear as limit sets of sequences $[a_n,b_n']$,
$[b_n,c_n']$ and $[c_n,a_n']$ of geodesics in $X$. We have that
$\dist (a_n,a_n'), \dist (b_n,b_n')$ and $\dist (c_n,c_n')$ are of
order $o(d_n)$ \uas.

\medskip

Let $d_n^A$ be the maximum between $\dist ([a_n,b_n'],
[a_n',c_n])$ and $\nu\theta$. Note that $d_n^A
>0$ and that $d_n^A=_\omega o(d_n)$. Take $a_n^1$ to be the
farthest from $a_n$ point on $[a_n,b_n']$ at distance $d_n^A$ from
$[a_n',c_n]$. Consider then $a_n^2$ the farthest from $a_n'$ point
on $[a_n',c_n]$ at distance $d_n^A$ from $a_n^1$. Obviously $\dist
(a_n^1,a_n^2)=d_n^A$.

The pairs of points $(b_n^1,b_n^2)$ in $[b_n,c_n']\times
[b_n',a_n]$, and respectively $(c_n^1,c_n^2)$ in $[c_n,a_n']\times
[c_n',b_n]$ are chosen similarly. Since the limit triangle
$\Delta$ is simple, it follows that the sets
$\{a_n,a_n',a_n^1,a_n^2\}$, $\{b_n,b_n',b_n^1,b_n^2\}$ and
$\{c_n,c_n',c_n^1,c_n^2\}$ have \uass diameters of order $o(d_n)$.
Hence the sequence of geodesic hexagons $H_n$ of vertices
$a_n^1,b_n^2, b_n^1, c_n^2, c_n^1, a_n^2$ with edges
$[a_n^1,b_n^2]\subset [a_n,b_n']$, $[b_n^1, c_n^2]\subset
[b_n,c_n']$, $[c_n^1, a_n^2]\subset [c_n,a_n']$, has the property
that $\lio{H_n}$ is $\Delta$. It remains to prove that $H_n$ is
\uass $(\theta ,\nu)$--fat.

$\mathbf{(F_1)}$ The fact that the edge $[a_n^1,a_n^2]$ is at
distance $O(d_n)$ from $[b_n^2, b_n^1]\cup [b_n^1, c_n^2]\cup
[c_n^2, c_n^1]$ and Lemma \ref{simplif} imply that $[a_n^1,a_n^2]$
satisfies property $(F_1)$.

In the same manner it can be shown that the edges $[b_n^1,b_n^2]$
and $[c_n^1,c_n^2]$ satisfy $(F_1)$.

The edge $[a_n^1,b_n^2]$ is at distance $O(d_n)$ from
$[c_n^1,c_n^2]$. The choice of $a_n^1$ and of the pair
$(b_n^1,b_n^2)$ implies that $[a_n^1,b_n^2]$ is at distance at
least $\nu\theta$ from $[b_n^1, c_n^2]\cup [c_n^1, a_n^2]$. Lemma
\ref{simplif} allows to conclude that $[a_n^1,b_n^2]$ satisfies
$(F_1)$.

 Similar arguments show that the edges $[b_n^1, c_n^2]$ and $[c_n^1,
 a_n^2]$ satisfy $(F_1)$.

$\mathbf{(F_2)}$ The distance from $a_n^1$ to $[a_n^2, c_n^1]$ is
at least $\nu\theta$ by the choice if $a_n^1$, while the distance
to $[b_n^2, b_n^1]\cup [b_n^1, c_n^2]\cup [c_n^2, c_n^1]$ is
$O(d_n)$. The same kind of argument shows that $(F_2)$ is
satisfied \uass  by all the vertices of $H_n$.
\endproof

In general not every geodesic in an asymptotic cone is a limit
geodesic (see the example in the end of Section \ref{prelac}).
Thus, in order to ensure property $(T_2)$ in every asymptotic cone
\wrt the collection of limit sets $\aau$, in
\cite{DrutuSapir:TreeGraded} property $(\alpha_3)$ in full
generality is used, together with the fact that limit sets are
closed, and the following result.

\begin{lemma}[\cite{DrutuSapir:TreeGraded}, Proposition 3.34]\label{lime}

Let $\Delta$ be an arbitrary simple geodesic triangle in $\co{X;
e,d}$. For every $\varepsilon >0$ sufficiently small there exists
$k_0=k_0(\varepsilon )$ and a simple geodesic triangle
$\Delta_\varepsilon$ with the following properties:
\begin{itemize}
  \item[(a)] $\dist_H\left( \Delta\, ,\, \Delta_\varepsilon \right) \leq \varepsilon$;
  \item[(b)] $\Delta_\varepsilon $ contains the midpoints
  of the edges of $\Delta $;
  \item[(c)] for every $\theta >0$ and $\nu \geq 8$ the triangle
  $\Delta_\varepsilon$ is the limit set $\lio{P^{\varepsilon}_n}$
  of a sequence $(P^{\varepsilon}_n)$ of geodesic $k$-gons in $X$,
for some $k\leq k_0$, that are $(\theta , \nu)$--fat \uas.
\end{itemize}
\end{lemma}

\begin{rmk}\label{nontr}
If $\Delta$ is non-trivial then the set of midpoints of edges of
$\Delta$ has cardinal $3$, hence the triangles
$\Delta_\varepsilon$ are also non-trivial.
\end{rmk}

In this section we prove that if in every asymptotic cone property
$(T_1)$ holds for the collection of limit sets $\aau$, then
property $(\beta_3)$ for the collection $\aaa$ suffices to deduce
$(T_2)$ for $\aau$, again in every asymptotic cone. To this
purpose, we define the following property in an asymptotic cone
$\co{X;e,d}\, $:

\begin{itemize}
\item[$(\Pi_k)$] \emph{ every simple non-trivial $k$-gon with edges
limit geodesics is contained in a subset from $\aaa_\omega$}.
\end{itemize}

\begin{cor}\label{corlimp}
Assume that in an asymptotic cone $\co{X;e,d}$ a collection
$\aaa_\omega$ of closed subsets satisfies properties $(T_1)$ and
$(\Pi_k)$ for every $k\in \N ,\, k\geq 3$. Then $\aaa_\omega$
satisfies $(T_2)$.
\end{cor}

\proof Consider a simple non-trivial geodesic triangle $\Delta$ in
$\co{X; e,d}$. By Lemma \ref{lime} for every large enough $k\in
\N$ there exists a simple non-trivial geodesic triangle $\Delta_k$
at Hausdorff distance at most $\frac{1}{k}$ from $\Delta$,
containing the midpoints of the edges of $\Delta$, moreover
$\Delta_k = \lio{P_n^{(k)}}$, where $P_n^{(k)}$ is $n$-\uass a
geodesic $m$-gon, $m=m(k)$. By property $(\Pi_m)$ the triangle
$\Delta_k$ is contained in some $A_k\in \aaa_\omega$. All $A_k$
contain the midpoints of the edges of $\Delta$. Property $(T_1)$
implies that there exists $A\in \aaa_\omega$ such that $A_k=A$ for
all $k$. All $\Delta_k$ are in $A$, $\Delta$ is the limit of
$\Delta_k$ in the Hausdorff distance, and $A$ is closed, therefore
$\Delta \subset A$.\endproof

\me

In view of Corollary \ref{corlimp} it suffices to prove that
$\aau$ satisfies $(\Pi_k)$ for all $k\geq 3$ to deduce that $\aau$
satisfies property $(T_2)$.

Obviously $(\Pi_k)$ implies $(\Pi_i)$ for every $i<k$. It turns
out that with the additional assumption that $(T_1)$ is satisfied,
the converse implication also holds.

\begin{lemma}\label{t2}
Assume that in an asymptotic cone $\co{X;e,d}$, the collection of
subsets $\aaa_\omega$ satisfies the properties $(T_1)$ and
$(\Pi_3)$. Then $\aau$ satisfies property $(\Pi_k)$ for every
$k\geq 3$.
\end{lemma}

\proof We prove property $(\Pi_k)$ by induction on $k$. The cases
$k=2$ and $k=3$ hold by hypothesis. Assume that the statement is
true for every $k\leq m-1$ and consider a simple non-trivial
geodesic $m$-gon $P$ in $\co{X;e,d}$, $m\geq 4$, with edges limit
geodesics.

Let $[x,y]$ and $[y,z]$ be two consecutive edges of $P$, in
clockwise order. Denote by $\call_1$ the union of the two edges
$[x,y]\cup [y,z]$ of $P$, and by $\call$ the union of the other
$m-2$ edges of $P$, in clockwise order.

Consider a limit geodesic $\g$ joining $x$ and $z$. If $\g$
coincides with $\call_1$ or with $\call$ then $P$ is a simple
geodesic polygon with at most $m-1$ edges, all of them limit
geodesics. By the inductive hypothesis $P$ is contained in a
subset $A$ in $\aau$.

Assume that $\g$ does not coincide either with $\call_1$ or with
$\call$.

\begin{figure}
\centering
\unitlength 1mm 
\linethickness{0.4pt}
\ifx\plotpoint\undefined\newsavebox{\plotpoint}\fi 
\begin{picture}(98.5,53.75)(0,0)
\thicklines \put(20.25,34.25){\line(1,0){65.75}}
\put(86,34.25){\line(0,1){.25}}
\qbezier(20.25,34.25)(35.63,53.75)(46.5,34.25)
\qbezier(46.5,34.25)(54.5,45.5)(60.5,34.75)
\qbezier(60.5,34.75)(74.13,52.38)(86.25,34.5)
\qbezier(20.5,34.25)(29.88,27.13)(34.75,34.5)
\qbezier(34.75,34.5)(44.5,21.13)(51.25,34.25)
\qbezier(51.25,34.25)(54.88,29.25)(58,34.25)
\qbezier(73.5,34.25)(81,24.13)(86.5,34.5)
\qbezier(58,34.25)(67.75,16.75)(73.5,34.25)
\put(15.75,35.25){\makebox(0,0)[cc]{$x$}}
\put(90.5,34.75){\makebox(0,0)[cc]{$z$}}
\put(46.75,40.5){\makebox(0,0)[cc]{$a$}}
\put(53,42){\makebox(0,0)[cc]{$\alpha$}}
\put(60.75,40.5){\makebox(0,0)[cc]{$b$}}
\put(54.25,36.75){\makebox(0,0)[cc]{$A$}}
\put(31.5,37.5){\makebox(0,0)[cc]{$A_1$}}
\put(72,37.5){\makebox(0,0)[cc]{$A_2$}}
\put(43.25,30.75){\makebox(0,0)[cc]{$A_1'$}}
\put(65.5,30.75){\makebox(0,0)[cc]{$A_2'$}}
\put(98,43.5){\makebox(0,0)[cc]{$\call$}}
\put(98.5,35){\makebox(0,0)[cc]{$\g$}}
\put(98.25,28.25){\makebox(0,0)[cc]{$\call_1$}}
\end{picture}
\caption{Step 1 in the proof of Lemma \ref{t2}} \label{lt2}
\end{figure}
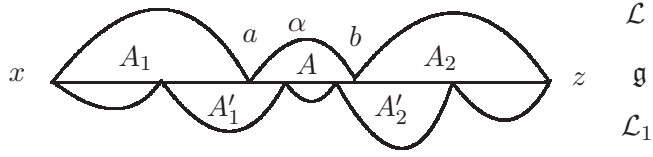

\textsc{Step 1.}\quad We prove that $\g \cup \call$ is contained
in some $A\in \aau$.

Let $\alpha \in \call \setminus \g$. Lemma \ref{lbigon} implies
that $\alpha$ is in the interior of a simple $\calt$-bigon formed
by $\call$ and $\g$, of endpoints $a$ and $b$, with $a$ closer to
$x$ than $b$ on $\call$. This $\calt$-bigon is a geodesic polygon
with at most $m-1$ edges which are limit geodesics, therefore by
the inductive hypothesis it is contained in a subset $A\in \aau$.

If $a=x$ and $b=z$ then $\g \cup \call $ is a simple $\calt$-bigon
and it is contained in some $A\in \aau$ by the inductive
hypothesis. Assume therefore that $(a,b)\neq (x,z)$. Without loss
of generality we may assume that $a\neq x$.

We apply Proposition \ref{braid} to $\call_1$, $\g_1=\g$, and
$\g_2 =\call_2$ the sub-arc of $\call$ in between $x$ and
$\alpha$. Property (3) is satisfied by the hypothesis of the
induction. It follows that the $\calt$-bigon formed by $\g$ and
$\call $ of endpoints $x$ and $a$ is contained in some $A_1\in
\aau$.

The point $a$ is in $\call \setminus \call_1$, hence it is in $\g
\setminus \call_1$. By Lemma \ref{lbigon}, $a$ is in the interior
of some simple $\calt$-bigon formed by $\call_1$ and $\g$, and by
$(\Pi_3)$ this $\calt$-bigon is in a subset $A_1'\in \aau$. Since
$A_1'\cap A$ and $A_1'\cap A_1$ contain non-trivial sub-arcs of
$\g$ property $(T_1)$ implies that $A=A_1'=A_1$.

If moreover $b\neq z$, a similar argument gives that the
$\calt$-bigon formed by $\g$ and $\call $ of endpoints $z$ and $b$
is contained in $A$ (see Figure \ref{lt2}).

We conclude that $\g \cup \call$ is contained in $A$.

\me

\textsc{Step 2.}\quad We prove that $\call_1$ is also contained in
$A$. Property $(\Pi_3)$ implies that any non-trivial simple
$\calt$-bigon formed
 by $\call_1$ and $\g$ is contained in a subset in $\aau$.
 We apply Lemma \ref{convbigon} to $\call_1$ and $\g \subset A$
  and we conclude that $\call_1\subset A$. Consequently $A$ contains
  $\call_1 \cup \call =P$.\endproof

\medskip

\begin{cor}\label{pi3}
Assume that in an asymptotic cone $\co{X;e,d}$, the collection of
closed subsets $\aaa_\omega$ satisfies properties $(T_1)$ and
$(\Pi_3)$. Then $\aau$ satisfies property $(T_2)$.
\end{cor}

\begin{cor}\label{cbeta3}
Let $X$ be a geodesic metric space and $\aaa$ a collection of
subsets in $X$, such that $(\beta_3)$ is satisfied and such that
in any asymptotic cone $\co{X;e,d}$, the collection of limit
subsets $\aaa_\omega$ satisfies property $(T_1)$. Then $\aau$
satisfies property $(T_2)$.
\end{cor}

Note that the only thing missing in Corollary \ref{cbeta3} to
conclude that $X$ is ATG \wrt $\aaa$ is that $\aau$ is composed of
geodesic subsets.

\me

Another useful consequence of Proposition \ref{limtr} is the
following.

\begin{cor}\label{chyp}
Let $(X,\dist )$ be a geodesic metric space. Assume that for some
$\theta >0$ and $\nu \geq 8$ the set of $(\theta , \nu)$--fat
geodesic hexagons is either empty or composed of hexagons of
uniformly bounded diameter. Then $X$ is hyperbolic.
\end{cor}

\proof Proposition \ref{limtr} implies that in any asymptotic cone
of $X$ any simple triangle with edges limit geodesics is trivial.
This statement can be extended by induction to all polygons.
Indeed, suppose that in any asymptotic cone of $X$ for all $3\leq
k\leq m-1$ all simple $k$-gons with edges limit geodesics are
trivial. Consider $P$ a simple $m$-gon with edges limit geodesics
in some $\co{X;e,d}$. Let $[x,y]$ and $[y,z]$ be two consecutive
edges of $P$ and let $\g$ be a limit geodesic joining $x$ and $z$.
All simple $\calt$-bigons formed by $[x,y]\cup [y,z]$ and $\g$
must be trivial by the inductive hypothesis, thus $\g=[x,y]\cup
[y,z]$. It follows that $P$ is a simple $(m-1)$--gon with edges
limit geodesics, hence by the inductive hypothesis it is trivial.

Lemma \ref{lime} and Remark \ref{nontr} imply that in any
$\co{X;e,d}$ any simple geodesic triangle must be trivial. It
follows that $\co{X;e,d}$ is a real tree, and since this holds for
all asymptotic cones we conclude that $X$ is hyperbolic (\cite[$\S
2.A$]{Gromov:Asymptotic}, see also \cite[$\S
3$]{Drutu:survey}).\endproof

\subsection{New definitions, useful for the rigidity of relatively hyperbolic
groups}\label{ndr}

In this section new versions of the definition of an ATG metric
space are stated and proved. They will play an important part in
the proof of the quasi-isometric invariance of relative
hyperbolicity.

\begin{theorem}\label{newdef}
In Theorem \ref{tgi} the following modifications in the list of
properties defining an \atgs metric space can be made:

\begin{itemize}
    \item[$\mathbf{(M_1)}$] property $(\alpha_3)$ can
be replaced by property $(\beta_3)$;

    \item[$\mathbf{(M_2)}$] property $(\alpha_2)$ can be either maintained or replaced by one
of the following two properties:
\begin{itemize}
    \item[$(\beta_2)$] there exists $\epsilon >0$ and $M\geq 0$
        such that for any geodesic $\g$ of length $\ell$ and any
        $A\in \aaa$ satisfying $\g(0),\g(\ell) \in \nn_{\epsilon\ell}
        (A)$, the middle third
        $\g \left( \left[ \frac{\ell}{3}\, ,\, \frac{2\ell}{3} \right] \right)$
         is contained in $\nn_M (A)$;
        \item[$(\fq)$]\textbf{(uniform quasi-convexity of pieces)}
        there exists $t>0$ and $K_0\geq 0$ such that for every $A\in\aaa$, $K\geq K_0$ and
$x,y\in \nn_K(A)$, every geodesic joining $x$ and $y$ in $X$ is
contained in $\nn_{tK}(A)$.
\end{itemize}
\end{itemize}
\end{theorem}

\proof Assume that $X$ is ATG with respect to $\aaa$. The uniform
quasi-convexity of pieces $(\fq)$ is satisfied by \cite[Lemma
4.3]{DrutuSapir:TreeGraded}. Property $(\beta_2)$ can be obtained
for any $\epsilon < \frac{1}{6t}$, where $t$ is the constant from
$(\fq)$, as follows. Consider a geodesic $\g$ of length $\ell$ and
$A\in \aaa$ as in $(\beta_2)$. We may assume that $\epsilon \ell
\geq K_0$, otherwise $\g$ would be contained in
$\nn_{\frac{K_0}{2\epsilon}}$. By $(\fq)$ the geodesic $\g$ is
then contained in $\nn_{t\epsilon\ell}(A)$. If $\theta=t\epsilon
<\frac{1}{6}$ then by Theorem \ref{tgi} and Remark \ref{rstr}
there exists $M=M(\theta)$ such that $\g \left( \left[ 0\, ,\,
\frac{\ell}{3} \right] \right)$ and $\g \left(
\left[\frac{2\ell}{3}\, ,\, \ell \right] \right)$ intersect $\nn_M
(A)$. Uniform convexity implies that $\g \left( \left[
\frac{\ell}{3}\, ,\, \frac{2\ell}{3} \right] \right)$ is contained
in $\nn_{tM'} (A)$, where $M'=\max (M, D_0)$.

Property $(\beta_3)$ is a particular case of property
$(\alpha_3)$.

\medskip

It remains to prove the converse statements: any of the triples of
properties $(\alpha_1)\& (\alpha_2) \& (\beta_3)$, $(\alpha_1)\&
(\beta_2) \& (\beta_3)$ or $(\alpha_1)\& (\fq) \& (\beta_3)$
implies that $X$ is ATG with respect to $\aaa$.

\begin{lemma}[\cite{DrutuSapir:TreeGraded}, Lemma
4.3]\label{impl0}
Properties $(\alpha_1)$ and $(\alpha_2)$ imply $(\fq)$.
\end{lemma}

\begin{lemma}\label{impl1}
Properties $(\alpha_1)$ and $(\beta_2)$ imply $(\fq)$ with $K_0$
equal to the constant $M$ in $(\beta_2)$.
\end{lemma}

\proof Suppose by contradiction that for every $n\in \N^*$ there
exists $A_n \in \aaa$, $K_n\geq M$ and $x_n,y_n\in \nn_{K_n}
(A_n)$ such that a geodesic $[x_n,y_n]$ is not contained in
$\nn_{n K_n} (A_n)$. For each $n\in \N^*$ we define $D_n$ to be
the infimum over the distances $\dist (x_n,y_n)$ between pairs of
points satisfying the properties above for some set in $\aaa$. In
what follows we assume that we chose $x_n,y_n$ at distance
$\delta_n\leq D_n+1$ of each other. Since $[x_n,y_n]$ is in
$\nn_{\delta_n/2}( \{x_n,y_n\})\subset \nn_{\delta_n/2+K_n}( A_n)$
it follows that $\frac{1}{2n-2}\delta_n \geq K_n$. In particular
for $n$ large enough $K_n < \epsilon \delta_n$, where $\epsilon
>0$ is the constant in $(\beta_2)$. It follows that the middle
third $[a_n,b_n]$ of $[x_n,y_n]$ is contained in $\nn_M (A_n)$.
Since $K_n \geq M$, the fact that $[x_n,y_n]\not \subset
\nn_{nK_n}(A_n)$ implies that either $[x_n,a_n]$ or $[b_n,y_n]$ is
not contained in $\nn_{nK_n}(A_n)$. It follows that $D_n\leq
\frac{\delta_n}{3}\leq \frac{D_n+1}{3}$, hence that the sequence
$(D_n)$ is uniformly bounded. This contradicts the fact that $D_n
\geq (2n-2)M$.\endproof

\begin{lemma}\label{newv}
Let $P$ be a geodesic $k$-gon with two consecutive edges $[x,y]$
and $[y,z]$, such that $\dist (x,[y,z]) =\dist (x,y)$. If $P$ is
$(\theta , \nu )$--fat then the $(k+1)$--gon $P'$ obtained from
$P$ by adding as a vertex the point $v\in [x,y]$ with $\dist (v,
y)=\frac{\nu\theta}{2}$ is $\left( \theta \,
,\,\frac{\nu}{2}\right)$--fat.
\end{lemma}

\proof

Property $(F_1)$ for $P'$ follows easily from property $(F_1)$ for
$P$.

Property $(F_2)$ holds for all the vertices different from
$x,v,y$, by property $(F_2)$ in $P$.

The polygonal line $\oo_x(P')=\oo_x(P) \cup [v,y]$ is in the
$\frac{\nu\theta}{2}$--tubular neighborhood of $\oo_x(P)$, hence
at distance at least $\frac{\nu\theta}{2}$ from $x$.

The polygonal line $\oo_y(P')$ is equal to $\oo_y(P) \cup [x,v]$.
The line $\oo_y(P)$ is at distance $\geq \nu\theta$ from $y$ and
$[x,v]$ is at distance $\frac{\nu\theta}{2}$ from $y$.

Finally, $\oo_v(P')= \oo_y(P)\cup [y,z]$. Since $\dist
(v,y)=\frac{\nu\theta}{2}$ it follows that $\dist \left( v,\oo_y
\right)\geq \frac{\nu\theta}{2}$. If there exists $p\in [y,z]$
such that $\dist (v,p)<\frac{\nu\theta}{2}$ then $\dist (x,p)\leq
\dist(x,v)+\dist(v,p)<\dist(x,v)+ \frac{\nu\theta}{2} = \dist
(x,y)$. This contradicts the hypothesis that $\dist (x,[y,z])
=\dist (x,y)$.\endproof

\begin{lemma}\label{impl2}
Properties $(\alpha_1)$, $(\fq)$ and $(\beta_3)$ imply
$(\alpha_2)$ for small enough $\varepsilon
>0$, and~$(\beta_2)$.
\end{lemma}

\proof Assume that $(\fq)$ and $(\beta_3)$ are satisfied. Let $\g
:[0,\ell ]\to X$ be a geodesic with endpoints $x=\g (0)$ and $y=\g
(\ell )$ contained in $\nn_{\varepsilon \ell} (A)$ for some $A\in
\aaa$. We shall prove that for a fixed positive constant $D$, the
geodesic $\g$ intersects $\nn_D (A)$.

According to $(\fq)$, the geodesic $\g$ is contained in
$\nn_{t\varepsilon\ell} (A)$.

\medskip

\begin{n}\label{const}
We denote $t\varepsilon$ by $\epsilon$ and we assume in what
follows that $\epsilon <\frac{1}{8}$. We denote by $D$ the maximum
between $tK_0 + 4\nu\theta$, $\chi$ and $\diam_\delta$ (with the
notation \ref{ndiam}) for $\delta = \max (\chi , tK_0)$. Here $t$
and $K_0$ are the constants appearing in $(\fq)$, while $\nu ,
\theta , \chi$ are the constants appearing in $(\beta_3)$.
\end{n}

\medskip

Suppose by contradiction that $\g$ does not intersect $\nn_D (A)$.
Note that since $\epsilon \ell \geq D$ it follows that $\ell >
8D$.

Consider $x'$ and $y'$ points in $A$ such that $\dist (x,x')$ and
$\dist (y,y')$ are at most $\varepsilon \ell$. By $(\fq)$, a
geodesic $\g'$ joining $x'$ and $y'$ is contained in $\nn_{tK_0}
(A)$.

Let $c\in \g$ and $c'\in \g'$ be two points such that $\dist
(c,c')=\dist (\g ,\g')$. Without loss of generality we may suppose
that $\dist (x,c) \geq \frac{\ell}{2}\, $. We may also suppose
that $\dist (x,x')= \dist (x, \g')$. In order to transform the
$4$-gon of vertices $x,x',c,c',$ into a fat polygon we make the
following choices. Let $x_1$ be the point on $\g$ between $x$ and
$c$ which is farthest from $x$ and at distance at most
$2\nu\theta$ from $[x,x']$. Let $x_2$ be the farthest from $x$
point on $[x,x']$ which is at distance $2\nu\theta$ from $x_1$.

We prove in the sequel that the geodesic pentagon of vertices
$x_1,x_2,x',c',c$ is $(\theta ,2\nu)$--fat. To simplify we shall
denote its edges by $[v,w]$ if $v,w$ are two consecutive vertices,
keeping in mind that $[x_1,c] \subset \g$ and that $[x',c']\subset
\g'$.

\medskip

\noindent$\mathbf{(F_1)}$\quad A point in $[c,c']\setminus
\nn_{2\theta} (\{c,c'\})$ is at distance at least $\left(
\frac{1}{2}-2\epsilon\right)\ell$ from $[x,x']$, hence at distance
at least $\left( \frac{1}{2}-2\epsilon\right)\ell - 2\nu\theta$ of
$[x_1,x_2]$. Since $\ell > 8D > 32 \nu\theta$, it follows that
$[c,c']\setminus \nn_{2\theta}
 (\{c,c'\})$ is at distance at least $\theta$ from $[x_1,x_2]\cup [x_2,
 x']$.

The choice of $c,c'$ implies that all points in
 $[c,c']\setminus \nn_{2\theta} (\{c,c'\})$ are at distance at
 least $2\theta$ from $\g$ and from $\g'$.

 The points in $[x_1,c]\setminus \nn_{2\theta} (\{x_1,c \})$ are
 at distance at least $D-tK_0$ from $\g'$, and at distance
  at least $2\nu\theta $ from $[x_2,x']$. Lemma \ref{simplif} allows to conclude that $[x_1,c]$ satisfies
 property $(F_1)$.

 The distance between $[x_1,x_2]$ and $[c,c']$ is at least $\left(
\frac{1}{2}-2\epsilon\right)\ell - 2\nu\theta$, and the one
between $[x_1,x_2]$ and $[x',c']$ is at least $D -tK_0-
2\nu\theta$. Thus, it suffices to verify that the distance between
$[x_1,x_2]\setminus \nn_{2\theta}
 (\{x_1,x_2\})$ and $[c,x_1]\cup
 [x_2,x']$ is at least $\theta$. According to the choices of
 $x_1,x_2$ this distance is $2\theta$. This, and Lemma
 \ref{simplif} also imply that $[x_2,x']$ satisfies $(F_1)$.

The fact that the edge $[x',c']$ is at distance at least $D-tK_0 -
2\nu\theta$ from $\g \cup [x_1,x_2]$, together with Lemma
\ref{simplif}, imply that $[x',c']$ satisfies $(F_1)$.

\medskip

\noindent$\mathbf{(F_2)}$\quad The vertex $c'$ is at distance at
least $D-tK_0-2\nu\theta$ from $[c,x_1]\cup [x_1,x_2]$ and at
distance at least $\left( \frac{1}{2}-2\epsilon\right)\ell$ from
$[x_2,x']$.

The vertex $c$ is at distance at least $D-tK_0$ from $[x',c']$ and
at distance at least $\left( \frac{1}{2}-2\epsilon\right)\ell
-2\nu\theta$ from $[x_1,x_2] \cup [x_2, x']$.

We have chosen $x_1$ at distance $2\nu\theta$ from $[x_2,x']$. The
same vertex is at distance at least $D-tK_0$ from $[x',c']$, and
at least $\left( \frac{1}{2}-2\epsilon\right)\ell -2\nu\theta$
from $[c,c']$.

Similarly, $x_2$ is at distance $2\nu\theta$ from $[x_1,c]$, at
distance at least $\left( \frac{1}{2}-2\epsilon\right)\ell
-2\nu\theta$ from $[c,c']$ and at least $D-tK_0-2\nu\theta$ from
$[x',c']$.

The vertex $x'$ is at distance at least $D-tK_0-2\nu\theta$ from
$[c,x_1]\cup [x_1,x_2]$ and at least $\left(
\frac{1}{2}-2\epsilon\right)\ell$ from $[c,c']$.

The pentagon of vertices $x_1,x_2,x', c',c$ is $(\theta ,
2\nu)$--fat. Lemma \ref{newv} and the fact that $\dist(c,c')=
\dist (c, [c',x'])$ implies that by adding a vertex on $[c,c']$
this pentagon becomes a hexagon $(\theta , \nu)$--fat. Therefore
by $(\beta_3)$ it is contained in $\nn_\chi (A')$ for some $A'\in
\aaa$. In particular the edge $[x',c']$ is contained in $\nn_\chi
(A') \cap \nn_{tK_0}(A)$. This edge has length at least
$\frac{\ell}{4}
> 2D$ and $D$ is at least $\diam_\delta$ for $\delta =\max (\chi , tK_0)$.
It follows that $A=A'$ and that
$D<\chi$, which is a contradiction.

We conclude that property $(\alpha_2)$ is satisfied for
$\varepsilon < \frac{1}{8t}$ and for $D$ chosen above.

Property $(\beta_2)$ is obtained as follows. If a geodesic $\g :
[0,\ell] \to X$ joins two points in $\nn_{\delta \ell } (A)$ then
it is contained in $\nn_{t\delta \ell } (A)$ by $(\fq)$. If
$\delta <\frac{\varepsilon}{3t}$ then by $(\alpha_2)$ the
sub-geodesics $g\left(\left[ 0, \frac{\ell}{3} \right]\right)$ and
$g\left(\left[\frac{2\ell}{3} , \ell \right]\right)$ intersect
$\nn_M (A)$. Then by $(\fq)$, $g\left(\left[\frac{\ell}{3} ,
\frac{2\ell}{3} \right]\right)$ is contained in $\nn_{tM} (A)$.
\endproof

In view of Lemmata \ref{impl0}, \ref{impl1} and \ref{impl2}, in
order to finish the proof of Theorem \ref{newdef} it suffices to
prove the following.

\begin{lemma}\label{rec}
Assume that $(X,\dist)$ is a geodesic space and $\aaa$ a
collection of subsets, satisfying the properties $(\alpha_1)$,
$(\alpha_2)$, $(\fq)$ and $(\beta_3)$. Then $X$ is ATG with
respect to $\aaa$.
\end{lemma}

\proof In an asymptotic cone $\co{X;e,d}$, the limit sets in
$\aaa_\omega$ are closed. Property $(\fq)$ easily implies that all
subsets in $\aau$ are geodesic.



 Property $(T_1)$ for $\aaa_\omega$ is deduced from
$(\alpha_1)$ and $(\alpha_2)$ as in \cite[Lemma
4.5]{DrutuSapir:TreeGraded}. Properties $(T_1)\& (\beta_3)$ imply
$(T_2)$ by Corollary \ref{cbeta3}. \hspace*{\fill}$\square$

\me

\begin{proposition}\label{b3big}
For any $\eta >0$ property $(\beta_3)$ can be replaced by the
following:
\begin{itemize}
\item[$(\beta_3^\eta )$] there exists $\theta >0$, $\nu \geq 8$ and $\chi>0$
        such that any geodesic hexagon $(\theta , \nu)$--fat of diameter at least $\eta $
        is contained in $\nn_\chi (A)$, for some $A\in
        \aaa$.
\end{itemize}
\end{proposition}

\proof Indeed, as a sufficient condition $(\beta_3^\eta )$ is used
to prove $(\Pi_3)$ in any asymptotic cone, by means of Proposition
\ref{limtr}. Given the sequence of hexagons $H_n$ in Proposition
\ref{limtr}, $H_n$ has \uass diameter of order $O(d_n)$. Property
$(\beta_3^\eta )$ suffices therefore to obtain property $(\Pi_3)$.

Property $(\beta_3)$ is also used in Lemma \ref{impl2} to prove
$(\alpha_2)$. It suffices to take in that proof the constant $D$
larger than $tK_0 +\eta$ to obtain that the geodesic pentagon with
vertices $x_1,x_2,x',c',c$ has diameter at least $\dist (c,c')\geq
D-tK_0 > \eta$. That pentagon is $(\theta , 2\nu)$--fat, hence by
Lemma \ref{newv} it can be made into a hexagon $(\theta ,
\nu)$--fat of diameter larger than $\eta$; therefore it is
contained in $\nn_\chi (A')$ for some $A'\in \aaa$. The rest of
the argument is carried out similarly.\endproof

\begin{cor}\label{reda}
Let $\aaa_{red}^\eta $ be the set of $A\in \aaa$ such that
$\nn_\chi (A)$ contains a $(\theta , \nu)$--fat geodesic hexagon
of diameter at least $\eta$. Then the space $X$ is ATG \wrt
$\aaa_{red}^\eta $.
\end{cor}

\proof Since $\aaa_{red}^\eta \subset \aaa$, properties
$(\alpha_1)$ and $(\beta_2)$ are satisfied. Property
$(\beta_3^\eta )$ is also satisfied by $\aaa_{red}^\eta$, hence by
Proposition \ref{b3big}, $X$ is ATG with respect to
$\aaa_{red}^\eta$.\endproof

\begin{cor}\label{redd}
For every $\lambda >0$ the space $X$ is also ATG \wrt the subset
$\aaa_\lambda$ in $\aaa$ composed of all the subsets of diameter
at least $\lambda$ in $\aaa$.
\end{cor}

\proof Indeed $\aaa_\lambda \subset \aaa$ implies that properties
$(\alpha_1)$ and $(\beta_2)$ are still satisfied.

Let $\eta = \lambda +2\chi $. Then $\aaa_{red}^\eta \subset
\aaa_\lambda$, which implies that property $(\beta_3^\eta )$ is
satisfied by $\aaa_\lambda$. By Proposition \ref{b3big}, $X$ is
ATG with respect to $\aaa_\lambda$.\endproof

\subsection{New definition, closer to the definition of
hyperbolicity}\label{ndh}

In \cite{DrutuSapir:RD} a version for groups of the following
notion has been introduced.

\begin{definition}\label{star}

Let $X$ be a geodesic metric space and let $\aaa$ be a collection
of subsets of $X$. We say that $X$ is
$(*)$--{\textit{asymptotically tree-graded with respect to}}
$\aaa$ if for every $C\geq 0$ there exist two constants $\sigma$
and $\delta$ such that every triangle $xyz$ with $(1,C)$--almost
geodesic edges is in one of the following two cases:
\begin{itemize}
  \item[\bc] there exists $a\in X$ such that $\overline{B}(a,\sigma )$ intersects each
 of the sides of the triangle;
  \item[\bp] there exists $A\in \aaa$ such that $\overline{\nn}_\sigma(A)$ intersects
each of the sides of the triangle, and the entrance (resp. exit)
points $x_1, y_1, z_1$ (resp. $y_2, z_2, x_2$) of the sides
$[x,y], [y,z], [z,x]$ in (from) $\overline{\nn}_\sigma(A)$ satisfy
$$\dist(x_1,x_2)<\delta,\, \dist(y_1, y_2)<\delta,\, \dist(z_1,
z_2)<\delta\, .$$

See Figure \ref{fig2}, taken from \cite{DrutuSapir:TreeGraded}.
\end{itemize}
\end{definition}

\begin{figure}
\centering
\unitlength .7mm 
\linethickness{0.4pt}
\ifx\plotpoint\undefined\newsavebox{\plotpoint}\fi 
\begin{picture}(105.5,87)(0,0)

\qbezier(14.5,9.25)(42.38,38.38)(53.75,84)
\qbezier(53.75,84)(68.38,31.38)(101.5,10.25)
\qbezier(101.5,10.25)(59.5,24.25)(14.5,9.25)
\put(86.65,37.5){\line(0,1){1.26}}
\put(86.63,38.76){\line(0,1){1.258}}
\multiput(86.55,40.02)(-.0308,.3135){4}{\line(0,1){.3135}}
\multiput(86.43,41.27)(-.02872,.20804){6}{\line(0,1){.20804}}
\multiput(86.26,42.52)(-.03159,.17722){7}{\line(0,1){.17722}}
\multiput(86.03,43.76)(-.03369,.15386){8}{\line(0,1){.15386}}
\multiput(85.77,44.99)(-.03176,.12194){10}{\line(0,1){.12194}}
\multiput(85.45,46.21)(-.03319,.10964){11}{\line(0,1){.10964}}
\multiput(85.08,47.42)(-.0317,.0916){13}{\line(0,1){.0916}}
\multiput(84.67,48.61)(-.032746,.083838){14}{\line(0,1){.083838}}
\multiput(84.21,49.78)(-.033606,.076991){15}{\line(0,1){.076991}}
\multiput(83.71,50.94)(-.032291,.066719){17}{\line(0,1){.066719}}
\multiput(83.16,52.07)(-.032943,.061769){18}{\line(0,1){.061769}}
\multiput(82.57,53.18)(-.033479,.05725){19}{\line(0,1){.05725}}
\multiput(81.93,54.27)(-.032297,.050571){21}{\line(0,1){.050571}}
\multiput(81.25,55.33)(-.032697,.047027){22}{\line(0,1){.047027}}
\multiput(80.53,56.37)(-.033014,.043722){23}{\line(0,1){.043722}}
\multiput(79.77,57.37)(-.033256,.040629){24}{\line(0,1){.040629}}
\multiput(78.97,58.35)(-.033429,.037722){25}{\line(0,1){.037722}}
\multiput(78.14,59.29)(-.03354,.034984){26}{\line(0,1){.034984}}
\multiput(77.27,60.2)(-.034885,.033643){26}{\line(-1,0){.034885}}
\multiput(76.36,61.07)(-.037624,.03354){25}{\line(-1,0){.037624}}
\multiput(75.42,61.91)(-.040531,.033375){24}{\line(-1,0){.040531}}
\multiput(74.45,62.71)(-.043625,.033142){23}{\line(-1,0){.043625}}
\multiput(73.44,63.48)(-.046931,.032835){22}{\line(-1,0){.046931}}
\multiput(72.41,64.2)(-.050476,.032445){21}{\line(-1,0){.050476}}
\multiput(71.35,64.88)(-.057151,.033647){19}{\line(-1,0){.057151}}
\multiput(70.26,65.52)(-.061672,.033125){18}{\line(-1,0){.061672}}
\multiput(69.15,66.12)(-.066624,.032487){17}{\line(-1,0){.066624}}
\multiput(68.02,66.67)(-.072086,.031718){16}{\line(-1,0){.072086}}
\multiput(66.87,67.18)(-.083741,.032992){14}{\line(-1,0){.083741}}
\multiput(65.7,67.64)(-.091506,.031969){13}{\line(-1,0){.091506}}
\multiput(64.51,68.05)(-.10954,.03351){11}{\line(-1,0){.10954}}
\multiput(63.3,68.42)(-.12185,.03212){10}{\line(-1,0){.12185}}
\multiput(62.08,68.74)(-.13668,.03035){9}{\line(-1,0){.13668}}
\multiput(60.85,69.02)(-.17713,.03211){7}{\line(-1,0){.17713}}
\multiput(59.61,69.24)(-.20796,.02933){6}{\line(-1,0){.20796}}
\multiput(58.37,69.42)(-.3134,.0317){4}{\line(-1,0){.3134}}
\put(57.11,69.54){\line(-1,0){1.258}}
\put(55.85,69.62){\line(-1,0){2.52}}
\put(53.33,69.63){\line(-1,0){1.258}}
\put(52.08,69.56){\line(-1,0){1.254}}
\multiput(50.82,69.44)(-.24975,-.03373){5}{\line(-1,0){.24975}}
\multiput(49.57,69.27)(-.17731,-.03107){7}{\line(-1,0){.17731}}
\multiput(48.33,69.05)(-.15396,-.03324){8}{\line(-1,0){.15396}}
\multiput(47.1,68.79)(-.12203,-.0314){10}{\line(-1,0){.12203}}
\multiput(45.88,68.47)(-.10974,-.03287){11}{\line(-1,0){.10974}}
\multiput(44.67,68.11)(-.091692,-.03143){13}{\line(-1,0){.091692}}
\multiput(43.48,67.7)(-.083934,-.032499){14}{\line(-1,0){.083934}}
\multiput(42.31,67.25)(-.077089,-.033379){15}{\line(-1,0){.077089}}
\multiput(41.15,66.75)(-.066814,-.032095){17}{\line(-1,0){.066814}}
\multiput(40.01,66.2)(-.061866,-.032761){18}{\line(-1,0){.061866}}
\multiput(38.9,65.61)(-.057348,-.03331){19}{\line(-1,0){.057348}}
\multiput(37.81,64.98)(-.050666,-.032148){21}{\line(-1,0){.050666}}
\multiput(36.75,64.3)(-.047123,-.032558){22}{\line(-1,0){.047123}}
\multiput(35.71,63.59)(-.043819,-.032885){23}{\line(-1,0){.043819}}
\multiput(34.7,62.83)(-.040726,-.033136){24}{\line(-1,0){.040726}}
\multiput(33.72,62.04)(-.037821,-.033318){25}{\line(-1,0){.037821}}
\multiput(32.78,61.2)(-.035083,-.033437){26}{\line(-1,0){.035083}}
\multiput(31.87,60.33)(-.032496,-.033498){27}{\line(0,-1){.033498}}
\multiput(30.99,59.43)(-.033651,-.037525){25}{\line(0,-1){.037525}}
\multiput(30.15,58.49)(-.033494,-.040432){24}{\line(0,-1){.040432}}
\multiput(29.34,57.52)(-.03327,-.043527){23}{\line(0,-1){.043527}}
\multiput(28.58,56.52)(-.032972,-.046834){22}{\line(0,-1){.046834}}
\multiput(27.85,55.49)(-.032593,-.05038){21}{\line(0,-1){.05038}}
\multiput(27.17,54.43)(-.032124,-.0542){20}{\line(0,-1){.0542}}
\multiput(26.53,53.35)(-.033306,-.061574){18}{\line(0,-1){.061574}}
\multiput(25.93,52.24)(-.032683,-.066528){17}{\line(0,-1){.066528}}
\multiput(25.37,51.11)(-.031929,-.071992){16}{\line(0,-1){.071992}}
\multiput(24.86,49.96)(-.033238,-.083644){14}{\line(0,-1){.083644}}
\multiput(24.4,48.78)(-.032238,-.091412){13}{\line(0,-1){.091412}}
\multiput(23.98,47.6)(-.03102,-.10032){12}{\line(0,-1){.10032}}
\multiput(23.6,46.39)(-.03247,-.12175){10}{\line(0,-1){.12175}}
\multiput(23.28,45.18)(-.03075,-.13659){9}{\line(0,-1){.13659}}
\multiput(23,43.95)(-.03263,-.17703){7}{\line(0,-1){.17703}}
\multiput(22.77,42.71)(-.02994,-.20787){6}{\line(0,-1){.20787}}
\multiput(22.59,41.46)(-.0327,-.3133){4}{\line(0,-1){.3133}}
\put(22.46,40.21){\line(0,-1){1.257}}
\put(22.38,38.95){\line(0,-1){1.26}}
\put(22.35,37.69){\line(0,-1){1.26}}
\put(22.37,36.43){\line(0,-1){1.258}}
\put(22.43,35.17){\line(0,-1){1.255}}
\multiput(22.55,33.92)(.033,-.24985){5}{\line(0,-1){.24985}}
\multiput(22.72,32.67)(.03054,-.1774){7}{\line(0,-1){.1774}}
\multiput(22.93,31.42)(.03279,-.15406){8}{\line(0,-1){.15406}}
\multiput(23.19,30.19)(.03104,-.12213){10}{\line(0,-1){.12213}}
\multiput(23.5,28.97)(.03255,-.10983){11}{\line(0,-1){.10983}}
\multiput(23.86,27.76)(.031161,-.091784){13}{\line(0,-1){.091784}}
\multiput(24.26,26.57)(.032252,-.084029){14}{\line(0,-1){.084029}}
\multiput(24.72,25.39)(.033152,-.077187){15}{\line(0,-1){.077187}}
\multiput(25.21,24.24)(.031898,-.066908){17}{\line(0,-1){.066908}}
\multiput(25.76,23.1)(.032579,-.061962){18}{\line(0,-1){.061962}}
\multiput(26.34,21.98)(.033141,-.057446){19}{\line(0,-1){.057446}}
\multiput(26.97,20.89)(.033599,-.053298){20}{\line(0,-1){.053298}}
\multiput(27.64,19.83)(.032419,-.047219){22}{\line(0,-1){.047219}}
\multiput(28.36,18.79)(.032756,-.043916){23}{\line(0,-1){.043916}}
\multiput(29.11,17.78)(.033016,-.040823){24}{\line(0,-1){.040823}}
\multiput(29.9,16.8)(.033207,-.037918){25}{\line(0,-1){.037918}}
\multiput(30.73,15.85)(.033334,-.035181){26}{\line(0,-1){.035181}}
\multiput(31.6,14.93)(.033402,-.032594){27}{\line(1,0){.033402}}
\multiput(32.5,14.05)(.035987,-.032462){26}{\line(1,0){.035987}}
\multiput(33.44,13.21)(.040334,-.033613){24}{\line(1,0){.040334}}
\multiput(34.41,12.4)(.043429,-.033398){23}{\line(1,0){.043429}}
\multiput(35.4,11.64)(.046737,-.03311){22}{\line(1,0){.046737}}
\multiput(36.43,10.91)(.050284,-.032741){21}{\line(1,0){.050284}}
\multiput(37.49,10.22)(.054105,-.032283){20}{\line(1,0){.054105}}
\multiput(38.57,9.57)(.061476,-.033487){18}{\line(1,0){.061476}}
\multiput(39.68,8.97)(.066432,-.032879){17}{\line(1,0){.066432}}
\multiput(40.81,8.41)(.071898,-.032141){16}{\line(1,0){.071898}}
\multiput(41.96,7.9)(.083546,-.033484){14}{\line(1,0){.083546}}
\multiput(43.13,7.43)(.091316,-.032507){13}{\line(1,0){.091316}}
\multiput(44.31,7.01)(.10023,-.03131){12}{\line(1,0){.10023}}
\multiput(45.52,6.63)(.12166,-.03283){10}{\line(1,0){.12166}}
\multiput(46.73,6.3)(.1365,-.03116){9}{\line(1,0){.1365}}
\multiput(47.96,6.02)(.17693,-.03315){7}{\line(1,0){.17693}}
\multiput(49.2,5.79)(.20778,-.03056){6}{\line(1,0){.20778}}
\multiput(50.45,5.61)(.3132,-.0336){4}{\line(1,0){.3132}}
\put(51.7,5.47){\line(1,0){1.257}}
\put(52.96,5.39){\line(1,0){1.26}}
\put(54.22,5.35){\line(1,0){1.26}}
\put(55.48,5.36){\line(1,0){1.259}}
\put(56.73,5.43){\line(1,0){1.255}}
\multiput(57.99,5.54)(.24994,.03226){5}{\line(1,0){.24994}}
\multiput(59.24,5.7)(.17749,.03002){7}{\line(1,0){.17749}}
\multiput(60.48,5.91)(.15416,.03234){8}{\line(1,0){.15416}}
\multiput(61.72,6.17)(.12222,.03068){10}{\line(1,0){.12222}}
\multiput(62.94,6.48)(.10993,.03222){11}{\line(1,0){.10993}}
\multiput(64.15,6.83)(.09953,.03346){12}{\line(1,0){.09953}}
\multiput(65.34,7.23)(.084123,.032005){14}{\line(1,0){.084123}}
\multiput(66.52,7.68)(.077284,.032925){15}{\line(1,0){.077284}}
\multiput(67.68,8.17)(.071189,.033683){16}{\line(1,0){.071189}}
\multiput(68.82,8.71)(.062057,.032397){18}{\line(1,0){.062057}}
\multiput(69.93,9.3)(.057543,.032972){19}{\line(1,0){.057543}}
\multiput(71.03,9.92)(.053397,.033442){20}{\line(1,0){.053397}}
\multiput(72.1,10.59)(.047314,.03228){22}{\line(1,0){.047314}}
\multiput(73.14,11.3)(.044012,.032627){23}{\line(1,0){.044012}}
\multiput(74.15,12.05)(.04092,.032896){24}{\line(1,0){.04092}}
\multiput(75.13,12.84)(.038016,.033095){25}{\line(1,0){.038016}}
\multiput(76.08,13.67)(.035279,.03323){26}{\line(1,0){.035279}}
\multiput(77,14.53)(.032692,.033306){27}{\line(0,1){.033306}}
\multiput(77.88,15.43)(.032568,.035891){26}{\line(0,1){.035891}}
\multiput(78.73,16.37)(.033731,.040235){24}{\line(0,1){.040235}}
\multiput(79.54,17.33)(.033526,.043331){23}{\line(0,1){.043331}}
\multiput(80.31,18.33)(.033247,.046639){22}{\line(0,1){.046639}}
\multiput(81.04,19.35)(.032889,.050188){21}{\line(0,1){.050188}}
\multiput(81.73,20.41)(.032442,.05401){20}{\line(0,1){.05401}}
\multiput(82.38,21.49)(.033667,.061377){18}{\line(0,1){.061377}}
\multiput(82.99,22.59)(.033074,.066335){17}{\line(0,1){.066335}}
\multiput(83.55,23.72)(.032352,.071803){16}{\line(0,1){.071803}}
\multiput(84.07,24.87)(.03373,.083447){14}{\line(0,1){.083447}}
\multiput(84.54,26.04)(.032775,.09122){13}{\line(0,1){.09122}}
\multiput(84.96,27.22)(.03161,.10014){12}{\line(0,1){.10014}}
\multiput(85.34,28.43)(.03319,.12156){10}{\line(0,1){.12156}}
\multiput(85.67,29.64)(.03156,.13641){9}{\line(0,1){.13641}}
\multiput(85.96,30.87)(.03367,.17684){7}{\line(0,1){.17684}}
\multiput(86.19,32.11)(.03117,.20769){6}{\line(0,1){.20769}}
\multiput(86.38,33.35)(.0276,.2505){5}{\line(0,1){.2505}}
\put(86.52,34.61){\line(0,1){1.257}}
\put(86.61,35.86){\line(0,1){1.638}}
\put(12,7.75){\makebox(0,0)[cc]{$x$}}
\put(54.25,87){\makebox(0,0)[cc]{$y$}}
\put(105.5,8.25){\makebox(0,0)[cc]{$z$}}
\put(22.5,22.25){\makebox(0,0)[cc]{$x_1$}}
\put(61,73.5){\makebox(0,0)[cc]{$y_1$}}
\put(78.5,10){\makebox(0,0)[cc]{$z_1$}}
\put(46.25,73.5){\makebox(0,0)[cc]{$y_2$}}
\put(89.5,25.75){\makebox(0,0)[cc]{$z_2$}}
\put(31.25,10){\makebox(0,0)[cc]{$x_2$}}
\put(49.5,69.25){\circle*{1.41}} \put(58.5,69.5){\circle*{1.8}}
\put(84.25,24.75){\circle*{1.12}}
\put(77.75,15.75){\circle*{1.41}} \put(32.75,14){\circle*{1.5}}
\put(26,22.5){\circle*{1}}
\put(55,49.6){\makebox(0,0)[cc]{$\overline{\nn}_\sigma(A)$}}
\end{picture}
\centering \caption{Case \bp.} \label{fig2}
\end{figure}

\begin{rmk}\label{starhip}
If $X$ is a geodesic metric space in which for some constant
$\sigma >0$ every
    geodesic triangle satisfies property \bc, then $X$ is a
    hyperbolic space. Conversely, in a hyperbolic geodesic metric space for
    every $L\geq 1$ and $C\geq 0$ there exists $\sigma >0$ such
    that every triangle with $(L,C)$--quasi-geodesic edges
    satisfies property \bc.
\end{rmk}

\begin{rmks}\label{starrh}
\begin{itemize}
    \item[(1)] If a metric space $X$ is ATG
     with respect to a collection of subsets $\aaa$ then $X$ is
      $(*)$--ATG with respect to $\aaa$, by \cite[Corollary 8.14 and Lemma
8.19]{DrutuSapir:TreeGraded}.

Moreover, according to \cite[Corollary 8.14]{DrutuSapir:TreeGraded}
if a geodesic triangle is in case \bps then for every $\sigma'\geq
\sigma$ there exists $\delta'$ such that the pairs of entrance
points  in $\onn_{\sigma'} (A)$ are at distance at most $\delta'$.

   \item[(2)] The notion of $(*)$--ATG space is weaker
than the one of ATG space. For instance if $X$ is a geodesic
hyperbolic space and if $\aaa$ is any collection of subsets
covering $X$, then $X$ is $(*)$--ATG with respect to $\aaa$, and
the collection $\aaa$ needs not satisfy properties $(\alpha_1)$ or
$(\fq)$, for instance.
\end{itemize}
\end{rmks}

It turns out nevertheless that one can formulate an equivalent
definition of ATG metric spaces using the $(*)$--property.

\begin{theorem}\label{tstar}
Let $(X,\dist)$ be a geodesic metric space and let $\aaa$ be a
collection of subsets of $X$. The metric space $X$ is \atgs with
respect to $\aaa$ if and only if $(X,\aaa )$ satisfy properties
$(\alpha_1)$ and $(\alpha_2)$, and moreover $X$ is $(*)$--ATG with
respect to $\aaa$.
\end{theorem}

\begin{cvn}\label{sigmaM}
In order to simplify some technical arguments of the equivalence
we make the assumption that for all $C>0$ the constant $\sigma$ in
the $(*)$--property is larger than the constant $M$ appearing in
property $(\alpha_2)$. By Remark \ref{starrh}, (1), if $X$ is ATG
then such a choice of $\sigma$ is possible.
\end{cvn}

\proof The direct implication has already been discussed, we now
prove the converse statement. As in Section \ref{ndr}, from
$(\alpha_1)$ and $(\alpha_2)$ can be deduced property $(\fq)$.
This property implies that in any asymptotic cone $\co{X;e,d}$ the
collection $\aau$ is composed of closed geodesic subsets.

Again $(\alpha_1)$ and $(\alpha_2)$ imply property $(T_1)$ for
$\aau$. According to Corollary \ref{pi3}, it remains to prove
property $(\Pi_3)$.

\begin{lemma}\label{doie}
Let $(X,\dist)$ be a geodesic metric space and let $\aaa$ be a
collection of subsets of $X$ satisfying property $(\alpha_2)$ for
some $\varepsilon \in \left[ 0,1/2 \right)$ and $M>0$.

Let $\mu \geq \nu \geq M$, let $\g$ be a geodesic and $A$ a subset
in $\aaa$ such that $\g$ intersects $\onn_\nu (A)$. If $e_\mu$ and
$e_\nu$ are the entrance points of $\g$ in $\onn_\mu (A)$ and
respectively  $\onn_\nu (A)$ then $\dist (e_\mu, e_\nu) \leq
\frac{\mu}{\varepsilon}$.
\end{lemma}

\proof If $\varepsilon \dist (e_\mu, e_\nu) > \mu$ then by
$(\alpha_2)$ the sub-arc of $\g$ between $e_\mu$ and $ e_\nu$
intersects $\nn_M (A)\subset \nn_\nu (A)$, which contradicts the
definition of $e_\nu$.\endproof

\begin{lemma}\label{distc}
If $(X,\dist)$ is a geodesic metric space $(*)$--ATG with respect
to $\aaa$, and if moreover $\aaa$ satisfies property $(\alpha_2)$,
then for every $C \geq 0$ there exist $\kappa \geq 0$ and $\lambda
\geq 0$ such that the following holds. For any two geodesics $\g ,
\g'$ with $\g_-=\g_-'$ and $\dist (\g_+ , \g_+')\leq C$, any point
$z$ on $\g'$ is either contained in $\onn_\kappa (\g)$ or it is
contained in $\onn_\kappa (A)$ for some $A\in \aaa$ such that
$\onn_\kappa (A)$ intersects $\g$. Moreover in the latter case, if
$e,f$ and $e',f'$ are the entrance and exit points from
$\onn_\kappa (A)$ of $\g$ and respectively $\g'$, then $\dist
(e,e'),\, \dist (f,f')\leq \lambda$.
\end{lemma}

\proof Let $\pgot$ be the path $\g \sqcup \left[\g_+ , \g_+'
\right]$, where $\left[\g_+ , \g_+' \right]$ is a geodesic segment
joining $\g_+$ and $\g_+'$. It is a $(1,2C)$--almost geodesic. Let
$\sigma$ and $\delta$ be the constants of property $(*)$ for $2C$,
and let $z$ be an arbitrary point on $\g'$, dividing $\g'$ into
two sub-arcs, $\g_1$ and $\g_2$. The triangle $\Delta$ of edges
$\g_1$, $\g_2$ and $\pgot$ is either in case \bcs or in case \bp.

If it is in case \bcs then there exist $a_1\in \g_1$,  $a_2\in
\g_2$ and $b\in \pgot$ such that the set $\{a_1,a_2,b\}$ has
diameter at most $2\sigma$. The point $z$ is on a geodesic joining
$a_1$ and $a_2$, hence it is at distance at most $3\sigma$ from
$b$, thus it is contained in $\onn_{3\sigma +C} (\g) $.

If $\Delta$ is in case \bps then there exists $A\in \aaa$ with
$\onn_\sigma (A)$ intersecting $\g_1$, $\g_2$ and $\pgot$. Let
$x_1,z_1$, $z_2,y_1$ and $x_2,y_2$ be the entrance and exit points
from $\onn_\sigma (A)$ of $\g_1$, $\g_2$ and $\pgot$ respectively.
Then $\dist (x_1,x_2), \dist (y_1,y_2)$ and $\dist (z_1,z_2)$ are
all at most $\delta$. Since $z$ is on a geodesic joining $z_1$ and
$z_2$, $z\in \onn_{\sigma + \delta/2} (A)$. Note that $\onn_\sigma
(A)$ intersects $\pgot$, therefore $\onn_{\sigma+C} (A)$
intersects $\g$.

Take $\kappa =\max \left( 3\sigma +C\, ,\, \sigma
+\frac{\delta}{2}\, ,\, \sigma +C \right)$.

The points $x_1$ and $y_1$ are the entrance and respectively the
exit point of $\g'$ from $\onn_{\sigma} (A)$. If we consider $e'$
and $f'$ the entrance and exit points of $\g'$ from $\onn_{\kappa}
(A)$, Lemma \ref{doie} implies that $\dist (x_1,e')$ and $\dist
(y_1,f')$ are at most $\frac{\kappa}{\varepsilon}$. Hence $\dist
(e', x_2)$ and $\dist (f', y_2)$ are at most
$\frac{\kappa}{\varepsilon} + \delta$.

Let $e$ and $f$ be the entrance and exit points of $\g$ from
$\onn_{\kappa} (A)$. If either $x_2$ or both $x_2$ and $y_2$ are
in $\g$ then they are the entrance and respectively the exit point
of $\g$ from $\onn_{\sigma} (A)$. Lemma \ref{doie} implies that
either $\dist (x_2,e)$ or both $\dist (x_2,e)$ and $\dist (y_2,f)$
are at most $\frac{\kappa}{\varepsilon}$, hence that either $\dist
(e,e')$ or both $\dist (e,e'),\dist (f,f')$ are $O(1)$.

Assume that $y_2\in [\g_+, \g_+']$. Then $\g_+$ is in
$\onn_{\kappa} (A)$, hence $\g_+=f$. It follows that $\dist
(f,y_2)\leq C$ and that $\dist (f,f')\leq
C+\frac{\kappa}{\varepsilon} + \delta$.

Assume that $x_2\in [\g_+, \g_+']$. The point $\g_+$ is in
$\onn_{\kappa} (A)$, and if $\varepsilon \dist (e, \g_+)> \kappa$
then $\g$ intersects $\nn_M(A)$ between $e$ and $\g_+$. According
to convention \ref{sigmaM}, $\nn_M(A)\subset \nn_\sigma (A)$,
hence $\g$ intersects $\nn_\sigma (A)$ between $e$ and $\g_+$.
This contradicts the fact that $x_2$ is the entrance point of
$\pgot$ into $\onn_\sigma (A)$. Thus $\dist (e, \g_+) \leq
\frac{\kappa}{\varepsilon}$ and $\dist (e, x_2)\leq
\frac{\kappa}{\varepsilon} +C$, whence $\dist (e, e')\leq
2\frac{\kappa}{\varepsilon} +C +\delta$.\endproof

\begin{lemma}\label{4fat}
Let $(X,\dist)$ be a geodesic metric space $(*)$--ATG with respect
to a collection of subsets $\aaa$. Assume moreover that $(X,\aaa
)$ satisfy properties $(\alpha_1)$, $(\alpha_2)$ and $(\fq)$. Then
there exist $\theta >0$, $\nu \geq 8$ and $\chi >0$ such that any
geodesic quadrilateral which is $(\theta , \nu)$--fat is contained
in $\nn_\chi (A)$ for some $A\in \aaa$.
\end{lemma}

\proof Let $P$ be a $(\theta , \nu)$--fat geodesic quadrilateral
with vertices $x,y,z,w$ in counterclockwise order. Let $[x,z]$ be
a geodesic joining  the opposite vertices $x$ and $z$.

\me

\textsc{Case 1.} Assume that both geodesic triangles $xyz$ and
$xzw$ are in case \bc. Then there exists $a_1\in [x,y]\, ,\,
a_2\in [y,z]$ and $a_3\in [x,z]$ such that the set
$\{a_1,a_2,a_2\}$ has diameter at most $2\sigma$. Likewise there
exists $b_1\in [z,w]\, ,\, b_2\in [w,x]$ and $b_3\in [z,x]$ such
that $\{b_1,b_2,b_3\}$ has diameter at most $2\sigma$. If $\theta
>2\sigma$ then $a_1,a_2\in B(y, 2\theta)$ and $b_1,b_2\in B(w,
2\theta)$.

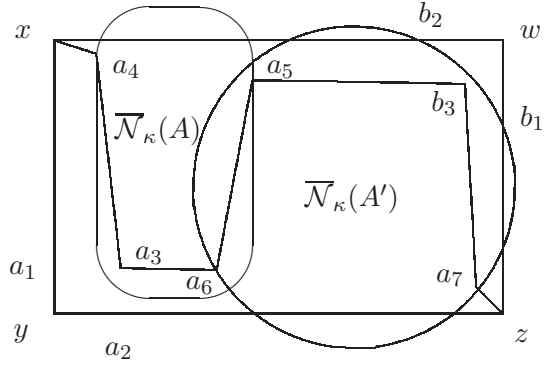
\begin{figure}[!ht]
\centering
\unitlength 1mm 
\linethickness{0.4pt}
\ifx\plotpoint\undefined\newsavebox{\plotpoint}\fi 
\begin{picture}(89.5,66.25)(0,0)
\put(26,25.5){\framebox(59.5,36.25)[cc]{}}
\put(41.88,46.88){\oval(20.75,38.75)[]}
\put(21.75,63){\makebox(0,0)[cc]{$x$}}
\put(21.5,22.75){\makebox(0,0)[cc]{$y$}}
\put(88,22.5){\makebox(0,0)[cc]{$z$}}
\put(89.25,63){\makebox(0,0)[cc]{$w$}}
\put(21.75,31.25){\makebox(0,0)[cc]{$a_1$}}
\put(38.5,33.25){\makebox(0,0)[cc]{$a_3$}}
\put(87.15,42.25){\line(0,1){.968}}
\put(87.12,43.22){\line(0,1){.966}}
\put(87.06,44.18){\line(0,1){.962}}
\multiput(86.95,45.15)(-.03056,.19128){5}{\line(0,1){.19128}}
\multiput(86.8,46.1)(-.03266,.15809){6}{\line(0,1){.15809}}
\multiput(86.6,47.05)(-.02983,.11733){8}{\line(0,1){.11733}}
\multiput(86.36,47.99)(-.03121,.10299){9}{\line(0,1){.10299}}
\multiput(86.08,48.92)(-.03226,.09133){10}{\line(0,1){.09133}}
\multiput(85.76,49.83)(-.03305,.08161){11}{\line(0,1){.08161}}
\multiput(85.4,50.73)(-.03365,.07336){12}{\line(0,1){.07336}}
\multiput(84.99,51.61)(-.03166,.061512){14}{\line(0,1){.061512}}
\multiput(84.55,52.47)(-.032117,.056015){15}{\line(0,1){.056015}}
\multiput(84.07,53.31)(-.032456,.051098){16}{\line(0,1){.051098}}
\multiput(83.55,54.13)(-.032692,.04666){17}{\line(0,1){.04666}}
\multiput(82.99,54.92)(-.032838,.042626){18}{\line(0,1){.042626}}
\multiput(82.4,55.69)(-.032905,.038933){19}{\line(0,1){.038933}}
\multiput(81.77,56.43)(-.032902,.035534){20}{\line(0,1){.035534}}
\multiput(81.12,57.14)(-.032834,.032389){21}{\line(-1,0){.032834}}
\multiput(80.43,57.82)(-.03598,.032413){20}{\line(-1,0){.03598}}
\multiput(79.71,58.47)(-.039379,.03237){19}{\line(-1,0){.039379}}
\multiput(78.96,59.08)(-.04307,.032253){18}{\line(-1,0){.04307}}
\multiput(78.18,59.66)(-.047103,.032051){17}{\line(-1,0){.047103}}
\multiput(77.38,60.21)(-.051536,.031755){16}{\line(-1,0){.051536}}
\multiput(76.56,60.72)(-.060481,.033588){14}{\line(-1,0){.060481}}
\multiput(75.71,61.19)(-.066703,.033188){13}{\line(-1,0){.066703}}
\multiput(74.85,61.62)(-.07381,.03265){12}{\line(-1,0){.07381}}
\multiput(73.96,62.01)(-.08205,.03193){11}{\line(-1,0){.08205}}
\multiput(73.06,62.36)(-.09176,.03101){10}{\line(-1,0){.09176}}
\multiput(72.14,62.67)(-.11633,.03353){8}{\line(-1,0){.11633}}
\multiput(71.21,62.94)(-.13455,.03226){7}{\line(-1,0){.13455}}
\multiput(70.27,63.16)(-.15852,.0305){6}{\line(-1,0){.15852}}
\multiput(69.32,63.35)(-.19168,.02795){5}{\line(-1,0){.19168}}
\put(68.36,63.49){\line(-1,0){.964}}
\put(67.39,63.58){\line(-1,0){.967}}
\put(66.43,63.64){\line(-1,0){.969}}
\put(65.46,63.64){\line(-1,0){.968}}
\put(64.49,63.61){\line(-1,0){.965}}
\multiput(63.52,63.53)(-.2402,-.0306){4}{\line(-1,0){.2402}}
\multiput(62.56,63.41)(-.19085,-.03317){5}{\line(-1,0){.19085}}
\multiput(61.61,63.24)(-.13511,-.02984){7}{\line(-1,0){.13511}}
\multiput(60.66,63.03)(-.11692,-.03143){8}{\line(-1,0){.11692}}
\multiput(59.73,62.78)(-.10255,-.03262){9}{\line(-1,0){.10255}}
\multiput(58.81,62.49)(-.09088,-.0335){10}{\line(-1,0){.09088}}
\multiput(57.9,62.15)(-.07439,-.03132){12}{\line(-1,0){.07439}}
\multiput(57,61.78)(-.067288,-.031985){13}{\line(-1,0){.067288}}
\multiput(56.13,61.36)(-.061074,-.032497){14}{\line(-1,0){.061074}}
\multiput(55.27,60.91)(-.055571,-.032879){15}{\line(-1,0){.055571}}
\multiput(54.44,60.41)(-.05065,-.033151){16}{\line(-1,0){.05065}}
\multiput(53.63,59.88)(-.04621,-.033326){17}{\line(-1,0){.04621}}
\multiput(52.84,59.32)(-.042173,-.033417){18}{\line(-1,0){.042173}}
\multiput(52.09,58.72)(-.03848,-.033434){19}{\line(-1,0){.03848}}
\multiput(51.35,58.08)(-.035081,-.033384){20}{\line(-1,0){.035081}}
\multiput(50.65,57.41)(-.033534,-.034938){20}{\line(0,-1){.034938}}
\multiput(49.98,56.71)(-.033599,-.038336){19}{\line(0,-1){.038336}}
\multiput(49.34,55.98)(-.033598,-.042029){18}{\line(0,-1){.042029}}
\multiput(48.74,55.23)(-.033524,-.046066){17}{\line(0,-1){.046066}}
\multiput(48.17,54.45)(-.033368,-.050507){16}{\line(0,-1){.050507}}
\multiput(47.64,53.64)(-.033118,-.05543){15}{\line(0,-1){.05543}}
\multiput(47.14,52.81)(-.032759,-.060934){14}{\line(0,-1){.060934}}
\multiput(46.68,51.95)(-.032273,-.06715){13}{\line(0,-1){.06715}}
\multiput(46.26,51.08)(-.03163,-.07425){12}{\line(0,-1){.07425}}
\multiput(45.88,50.19)(-.03081,-.08248){11}{\line(0,-1){.08248}}
\multiput(45.54,49.28)(-.03306,-.10241){9}{\line(0,-1){.10241}}
\multiput(45.24,48.36)(-.03194,-.11678){8}{\line(0,-1){.11678}}
\multiput(44.99,47.43)(-.03042,-.13498){7}{\line(0,-1){.13498}}
\multiput(44.78,46.48)(-.02833,-.15892){6}{\line(0,-1){.15892}}
\multiput(44.61,45.53)(-.0317,-.2401){4}{\line(0,-1){.2401}}
\put(44.48,44.57){\line(0,-1){.965}}
\put(44.4,43.6){\line(0,-1){1.936}}
\put(44.36,41.67){\line(0,-1){.967}}
\put(44.41,40.7){\line(0,-1){.964}}
\multiput(44.5,39.73)(.02713,-.1918){5}{\line(0,-1){.1918}}
\multiput(44.64,38.78)(.02981,-.15865){6}{\line(0,-1){.15865}}
\multiput(44.82,37.82)(.03168,-.13469){7}{\line(0,-1){.13469}}
\multiput(45.04,36.88)(.03303,-.11648){8}{\line(0,-1){.11648}}
\multiput(45.3,35.95)(.03061,-.09189){10}{\line(0,-1){.09189}}
\multiput(45.61,35.03)(.03158,-.08219){11}{\line(0,-1){.08219}}
\multiput(45.96,34.13)(.03233,-.07395){12}{\line(0,-1){.07395}}
\multiput(46.34,33.24)(.032901,-.066845){13}{\line(0,-1){.066845}}
\multiput(46.77,32.37)(.033328,-.060624){14}{\line(0,-1){.060624}}
\multiput(47.24,31.52)(.033635,-.055117){15}{\line(0,-1){.055117}}
\multiput(47.74,30.69)(.031849,-.04724){17}{\line(0,-1){.04724}}
\multiput(48.28,29.89)(.032068,-.043208){18}{\line(0,-1){.043208}}
\multiput(48.86,29.11)(.032201,-.039517){19}{\line(0,-1){.039517}}
\multiput(49.47,28.36)(.032259,-.036119){20}{\line(0,-1){.036119}}
\multiput(50.12,27.64)(.032248,-.032973){21}{\line(0,-1){.032973}}
\multiput(50.8,26.95)(.035392,-.033054){20}{\line(1,0){.035392}}
\multiput(51.5,26.29)(.038791,-.033072){19}{\line(1,0){.038791}}
\multiput(52.24,25.66)(.042484,-.033021){18}{\line(1,0){.042484}}
\multiput(53.01,25.06)(.04652,-.032892){17}{\line(1,0){.04652}}
\multiput(53.8,24.5)(.050958,-.032675){16}{\line(1,0){.050958}}
\multiput(54.61,23.98)(.055877,-.032357){15}{\line(1,0){.055877}}
\multiput(55.45,23.5)(.061376,-.031924){14}{\line(1,0){.061376}}
\multiput(56.31,23.05)(.067585,-.031353){13}{\line(1,0){.067585}}
\multiput(57.19,22.64)(.08147,-.0334){11}{\line(1,0){.08147}}
\multiput(58.08,22.27)(.09119,-.03265){10}{\line(1,0){.09119}}
\multiput(59,21.95)(.10286,-.03165){9}{\line(1,0){.10286}}
\multiput(59.92,21.66)(.11721,-.03034){8}{\line(1,0){.11721}}
\multiput(60.86,21.42)(.15794,-.03334){6}{\line(1,0){.15794}}
\multiput(61.81,21.22)(.19115,-.03138){5}{\line(1,0){.19115}}
\put(62.76,21.06){\line(1,0){.962}}
\put(63.72,20.95){\line(1,0){.966}}
\put(64.69,20.88){\line(1,0){1.937}}
\put(66.63,20.87){\line(1,0){.967}}
\put(67.59,20.93){\line(1,0){.963}}
\multiput(68.56,21.04)(.19141,.02974){5}{\line(1,0){.19141}}
\multiput(69.51,21.19)(.15822,.03198){6}{\line(1,0){.15822}}
\multiput(70.46,21.38)(.13424,.03352){7}{\line(1,0){.13424}}
\multiput(71.4,21.61)(.10312,.03077){9}{\line(1,0){.10312}}
\multiput(72.33,21.89)(.09146,.03186){10}{\line(1,0){.09146}}
\multiput(73.24,22.21)(.08175,.0327){11}{\line(1,0){.08175}}
\multiput(74.14,22.57)(.07351,.03334){12}{\line(1,0){.07351}}
\multiput(75.03,22.97)(.061647,.031396){14}{\line(1,0){.061647}}
\multiput(75.89,23.41)(.056153,.031877){15}{\line(1,0){.056153}}
\multiput(76.73,23.89)(.051237,.032236){16}{\line(1,0){.051237}}
\multiput(77.55,24.4)(.0468,.032491){17}{\line(1,0){.0468}}
\multiput(78.35,24.95)(.042766,.032655){18}{\line(1,0){.042766}}
\multiput(79.12,25.54)(.039074,.032738){19}{\line(1,0){.039074}}
\multiput(79.86,26.16)(.035675,.032749){20}{\line(1,0){.035675}}
\multiput(80.57,26.82)(.03253,.032695){21}{\line(0,1){.032695}}
\multiput(81.26,27.51)(.032568,.03584){20}{\line(0,1){.03584}}
\multiput(81.91,28.22)(.032539,.03924){19}{\line(0,1){.03924}}
\multiput(82.53,28.97)(.032437,.042931){18}{\line(0,1){.042931}}
\multiput(83.11,29.74)(.032253,.046965){17}{\line(0,1){.046965}}
\multiput(83.66,30.54)(.031976,.0514){16}{\line(0,1){.0514}}
\multiput(84.17,31.36)(.031591,.056314){15}{\line(0,1){.056314}}
\multiput(84.64,32.21)(.033474,.06656){13}{\line(0,1){.06656}}
\multiput(85.08,33.07)(.03296,.07367){12}{\line(0,1){.07367}}
\multiput(85.47,33.96)(.03229,.08192){11}{\line(0,1){.08192}}
\multiput(85.83,34.86)(.0314,.09162){10}{\line(0,1){.09162}}
\multiput(86.14,35.77)(.03025,.10328){9}{\line(0,1){.10328}}
\multiput(86.41,36.7)(.03284,.13441){7}{\line(0,1){.13441}}
\multiput(86.64,37.64)(.03118,.15838){6}{\line(0,1){.15838}}
\multiput(86.83,38.59)(.02877,.19156){5}{\line(0,1){.19156}}
\put(86.98,39.55){\line(0,1){.963}}
\put(87.08,40.52){\line(0,1){1.735}}
\multiput(26,61.75)(.1057692,-.0336538){52}{\line(1,0){.1057692}}
\multiput(31.5,60)(.03350515,-.29381443){97}{\line(0,-1){.29381443}}
\multiput(34.75,31.5)(1.59375,-.03125){8}{\line(1,0){1.59375}}
\multiput(47.5,31.25)(.03368794,.17907801){141}{\line(0,1){.17907801}}
\multiput(52.25,56.5)(1.883333,-.033333){15}{\line(1,0){1.883333}}
\multiput(80.5,56)(.0333333,-.6){45}{\line(0,-1){.6}}
\multiput(82,29)(.03350515,-.03350515){97}{\line(0,-1){.03350515}}
\put(65.5,41){\makebox(0,0)[cc]{$\onn_\kappa (A')$}}
\put(39.75,50){\makebox(0,0)[cc]{$\onn_\kappa (A)$}}
\put(34.5,20.5){\makebox(0,0)[cc]{$a_2$}}
\put(35.5,58){\makebox(0,0)[cc]{$a_4$}}
\put(56,58){\makebox(0,0)[cc]{$a_5$}}
\put(45.25,29.25){\makebox(0,0)[cc]{$a_6$}}
\put(77.75,53.5){\makebox(0,0)[cc]{$b_3$}}
\put(76,65.25){\makebox(0,0)[cc]{$b_2$}}
\put(89.5,51.5){\makebox(0,0)[cc]{$b_1$}}
\put(78.5,30.25){\makebox(0,0)[cc]{$a_7$}}
\end{picture}

\caption{Case 1 in proof of Lemma \ref{4fat}.} \label{caz1l4fat}
\end{figure}

Without loss of generality we may assume that $a_3\in [x, b_3]$.

\begin{n}
For $C=\max \left( 2\sigma\, ,\, \delta \right)$ let $\kappa$ and
$\lambda$ be the constants given by Lemma \ref{distc}.
\end{n}

Lemma \ref{distc} applied to $[x,b_2]$ and to $a_3\in [x, b_3]$
implies that either $a_3\in \onn_\kappa ([x,b_2])$, or $a_3\in
\onn_\kappa (A)$ such that the entrance respectively exit point,
$a_4,a_5$, of $[x,b_3]$ from $\onn_\kappa (A)$ are at distance at
most $\lambda$ from $[x,b_2]$. If $\theta > 2\sigma +\kappa$ then
by Lemma \ref{latop} the first case cannot occur.

In the second case we have that $\dist \left( a_3, \{ a_4, a_5 \}
\right)$ is at least $\theta - 2\sigma - \lambda$. Lemma
\ref{distc} applied to $[z,a_2]$ and to $a_5\in [z,a_3]$ implies
that either $a_5\in \onn_\kappa ([z,a_2])$ or that $a_5\in
\onn_\kappa (A')$ such that the entrance and the exit point
$a_6,a_7$, of $[a_3, z]$ from $\onn_\kappa (A')$ are at distance
at most $\lambda$ from $[z,a_2]$. The first case cannot occur if
$\theta
>\lambda + \kappa$. In the second case
$\dist \left( a_5\, ,\, \{ a_6, a_7\} \right)\geq \theta
-2\lambda$. The intersection $[a_3,a_5]\cap [a_6,a_5]$ has length
at least $\theta -2\sigma - 2 \lambda$. By property $(\fq)$ this
intersection is contained in $\nn_{t\kappa+1} (A) \cap
\nn_{t\kappa+1} (A')$. If $\theta > 2\sigma + 2\lambda +
\diam_{t\kappa +1}+1$ then $A=A'$.

The point $a_4$ is the entrance point of $[x,b_3]$ in $\onn_\kappa
(A)$ while $a_6$ is the entrance point of $[a_3,z]$ in
$\onn_\kappa (A)$. If $a_4\in [a_3,b_3]$ then $a_4=a_6$, and this
point is at distance at most $\lambda$ from both $[x,w]$ and
$[y,z]$. If $\theta > 2 \lambda$ then this cannot occur. Thus we
may assume that $a_4\in [x,a_3]$. Likewise we have that $a_7\in
[b_3, z]$ (see Figure \ref{caz1l4fat}).

We apply Lemma \ref{distc} to $[x,a_1]$ and to $a_4\in [x,a_3]$.
If we are in the second case of the conclusion then $a_4\in
\onn_\kappa (A'')$, and the entrance and exit point, $a_4',a_4''$,
of $[x,a_3]$ from $\onn_\kappa (A'')$ are at distance at most
$\lambda$ from $[x,a_1]$. If $\dist (a_4,a_4'')\geq \diam_{t\kappa
+1}+1$ then $A''=A$ and $a_4=a_4'$. Thus, in all cases $a_4$ is at
distance $O(1)$ from $[x, a_1]$. Recall that $a_4$ is at distance
at most $\lambda$ from $[x,w]$. It follows that if $\theta$ is
large enough then $a_4\in B(x, 2\theta +\lambda)$.

A similar argument gives that $a_7\in B(z, 2\theta +\lambda)$.

We have thus that $\{x,z\}\subset \onn_{2\theta + \lambda
+\kappa}(A)$. Also, since $\{a_3,b_3\}\subset [a_4,a_7]\subset
\onn_{t\kappa }(A)$ it follows that $\{y,w\}\subset \onn_{2\theta
+ 2\sigma + t\kappa} (A)$. By property $(\fq)$, $P\subset
\onn_\chi (A)$ where $\chi = t(2\theta +\lambda +2\sigma +t
\kappa)$.

\me

\textsc{Case 2.}\quad Assume that the triangle $xyz$ is in case
\bps while $xzw$ is in case \bc. Then there exists $A\in \aaa$
such that $\onn_\sigma (A)$ intersects all the edges of $xyz$.
Moreover if $x_2,y_1$ are the entrance and exit point of $[x,y]$
in $\onn_\sigma (A)$, while $y_2,z_1$ and $z_2,x_1$ are the
entrance and exit points of $[y,z]$ and respectively $[z,x]$ in
$\onn_\sigma (A)$ then $\dist (x_1,x_2), \dist (y_1,y_2)$ and
$\dist (z_1,z_2)$ are at most $\delta$.

Let also $b_1\in [z,w]$, $b_2\in [w,x]$ and $b_3\in [x,z]$ be such
that $\{b_1, b_2,b_3\}$ has diameter at most $2\sigma$. If $\theta
>2\sigma$ then property $(F_1)$ implies that $\{b_1, b_2 \}
\subset B(w, 2\theta)$.

\me

\textsc{Case 2.a.}\quad Assume that $b_3\in [x_1,z_2]$. Note that
$\dist (b_3 , x_1)\geq \dist (w, [x,y]) - 2 \theta -2\sigma -
\delta \geq 6\theta -2\sigma - \delta$. Same for $\dist (b_3 ,
z_2)$. Thus for $\theta$ large both $\dist (b_3 , x_1)$ and $\dist
(b_3 , z_2)$ are large.

Lemma \ref{distc} applied to $x_1\in [x,b_3]$ and to $[x, b_2]$
implies that either $x_1\in \onn_\kappa ([x,b_2])$ or $x_1 \in
\onn_\kappa (A_1)$ such that the entrance and exit points
$x_1',x_1''$ of $[x,b_3]$ from $\onn_\kappa (A_1)$ are at distance
at most $\lambda$ from $[x,b_2]$. In the latter case if $\dist
(x_1,x_1'')
> \diam_\tau$ where $\tau =t\max (\sigma , \kappa) +1$ then
$A_1=A$ and $x_1=x_1'$. Thus in all cases $\dist (x_1,
[x,w])=O(1)$. For $\theta$ large enough it follows that $x_2\in
B(x, 2\theta)$, hence $x\in \onn_{2\theta + \sigma }(A)$. A
similar argument gives that $z\in \onn_{2\theta + \sigma }(A)$.

If $\theta > \delta$ then $\dist (y,y_1)<2\theta$ and $y\in
\onn_{2\theta +\sigma} (A)$.

Also $b_3\in [x_1,z_2]\subset \onn_{t\sigma } (A)$, hence $w\in
\onn_{t\sigma + 2\theta +2\sigma} (A)$. We conclude by property
$(\fq)$ that $P\subset \nn_\chi (A)$ for some $\chi = O(1)$.

\me

\textsc{Case 2.b.}\quad Assume that $b_3\not\in [x_1,z_2]$.
Without loss of generality we may assume that $[x_1,z_2]\subset
[x,b_3)$.

Lemma \ref{distc} applied to $b_3\in [z_2,z]$ and to $[z_1,z]$
implies that either $b_3\in \onn_\kappa ([z_1,z])$ or that $b_3\in
\onn_\kappa (B)$ such that $\onn_\kappa (B)$ intersects $[z_1,z]$,
and if $b_4,b_5$ are the entrance and exit point of $[z_2,z]$ from
$\onn_\kappa (B)$ then these points are at distance at most
$\lambda$ from $[z_1,z]$. For $\theta$ large enough the first case
cannot occur. In the second case $\dist (b_3, \{ b_4,b_5\})\geq
\dist (w, [y,z]) - \lambda - 2(\theta + \sigma)\geq 6 \theta
-\lambda -2\sigma$.

Applying Lemma \ref{distc} now to $b_4\in [x,b_3]$ and to
$[x,b_2]$ gives that for $\theta$ large enough $b_4\in \onn_\kappa
(B')$ such that $\onn_\kappa (B')$ intersects $[x,b_3]$, and  the
entrance and exit points $b_6,b_7$ of $[x,b_3]$ from $\onn_\kappa
(B')$ are at distance at most $\lambda$ from $[x,b_2]$ (see Figure
\ref{caz2bl4fat}). Moreover $\dist (b_4 , \{b_6, b_7\})\geq \theta
-2 \lambda$. Thus $\onn_{t\kappa}(B)\cap \onn_{t\kappa}(B')$
contains $[b_4,b_7]\cap [b_4,b_3]$ of length $\min (6 \theta
-\lambda -2\sigma \, ,\, \theta -2 \lambda)$. For $\theta$ large
we conclude that $B=B'$.

\begin{figure}[!ht]
\centering
\unitlength 1mm 
\linethickness{0.4pt}
\ifx\plotpoint\undefined\newsavebox{\plotpoint}\fi 
\begin{picture}(116,72.25)(0,0)
\put(31.25,19.25){\framebox(76.5,49)[cc]{}}
\put(71.13,46){\oval(22.75,37)[]}
\put(74.95,42.25){\line(0,1){1.053}}
\put(74.93,43.3){\line(0,1){1.051}}
\put(74.86,44.35){\line(0,1){1.047}}
\multiput(74.74,45.4)(-.03195,.20814){5}{\line(0,1){.20814}}
\multiput(74.58,46.44)(-.02926,.14754){7}{\line(0,1){.14754}}
\multiput(74.38,47.47)(-.0312,.12786){8}{\line(0,1){.12786}}
\multiput(74.13,48.5)(-.03265,.11234){9}{\line(0,1){.11234}}
\multiput(73.84,49.51)(-.03069,.09067){11}{\line(0,1){.09067}}
\multiput(73.5,50.5)(-.03172,.08181){12}{\line(0,1){.08181}}
\multiput(73.12,51.49)(-.032535,.074171){13}{\line(0,1){.074171}}
\multiput(72.69,52.45)(-.033178,.067494){14}{\line(0,1){.067494}}
\multiput(72.23,53.4)(-.033677,.061588){15}{\line(0,1){.061588}}
\multiput(71.73,54.32)(-.032051,.052998){17}{\line(0,1){.052998}}
\multiput(71.18,55.22)(-.032419,.04869){18}{\line(0,1){.04869}}
\multiput(70.6,56.1)(-.03269,.044747){19}{\line(0,1){.044747}}
\multiput(69.98,56.95)(-.032875,.041119){20}{\line(0,1){.041119}}
\multiput(69.32,57.77)(-.032983,.037762){21}{\line(0,1){.037762}}
\multiput(68.63,58.56)(-.033022,.034642){22}{\line(0,1){.034642}}
\multiput(67.9,59.32)(-.034498,.033173){22}{\line(-1,0){.034498}}
\multiput(67.14,60.05)(-.037618,.033147){21}{\line(-1,0){.037618}}
\multiput(66.35,60.75)(-.040976,.033053){20}{\line(-1,0){.040976}}
\multiput(65.53,61.41)(-.044605,.032884){19}{\line(-1,0){.044605}}
\multiput(64.68,62.04)(-.048548,.03263){18}{\line(-1,0){.048548}}
\multiput(63.81,62.62)(-.052858,.032281){17}{\line(-1,0){.052858}}
\multiput(62.91,63.17)(-.057601,.031823){16}{\line(-1,0){.057601}}
\multiput(61.99,63.68)(-.067349,.033471){14}{\line(-1,0){.067349}}
\multiput(61.05,64.15)(-.074029,.032857){13}{\line(-1,0){.074029}}
\multiput(60.08,64.58)(-.08167,.03207){12}{\line(-1,0){.08167}}
\multiput(59.1,64.96)(-.09053,.03108){11}{\line(-1,0){.09053}}
\multiput(58.11,65.3)(-.1122,.03314){9}{\line(-1,0){.1122}}
\multiput(57.1,65.6)(-.12773,.03175){8}{\line(-1,0){.12773}}
\multiput(56.08,65.86)(-.14741,.02991){7}{\line(-1,0){.14741}}
\multiput(55.04,66.07)(-.208,.03285){5}{\line(-1,0){.208}}
\put(54,66.23){\line(-1,0){1.046}}
\put(52.96,66.35){\line(-1,0){1.05}}
\put(51.91,66.42){\line(-1,0){2.105}}
\put(49.8,66.43){\line(-1,0){1.051}}
\put(48.75,66.37){\line(-1,0){1.047}}
\multiput(47.7,66.26)(-.20828,-.03104){5}{\line(-1,0){.20828}}
\multiput(46.66,66.1)(-.17228,-.03339){6}{\line(-1,0){.17228}}
\multiput(45.63,65.9)(-.128,-.03064){8}{\line(-1,0){.128}}
\multiput(44.61,65.66)(-.11248,-.03216){9}{\line(-1,0){.11248}}
\multiput(43.59,65.37)(-.09988,-.03332){10}{\line(-1,0){.09988}}
\multiput(42.59,65.03)(-.08195,-.03136){12}{\line(-1,0){.08195}}
\multiput(41.61,64.66)(-.074312,-.032212){13}{\line(-1,0){.074312}}
\multiput(40.64,64.24)(-.067637,-.032884){14}{\line(-1,0){.067637}}
\multiput(39.7,63.78)(-.061733,-.033409){15}{\line(-1,0){.061733}}
\multiput(38.77,63.28)(-.053137,-.03182){17}{\line(-1,0){.053137}}
\multiput(37.87,62.74)(-.04883,-.032207){18}{\line(-1,0){.04883}}
\multiput(36.99,62.16)(-.044889,-.032495){19}{\line(-1,0){.044889}}
\multiput(36.14,61.54)(-.041261,-.032696){20}{\line(-1,0){.041261}}
\multiput(35.31,60.89)(-.037905,-.032819){21}{\line(-1,0){.037905}}
\multiput(34.52,60.2)(-.034785,-.032872){22}{\line(-1,0){.034785}}
\multiput(33.75,59.47)(-.033322,-.034353){22}{\line(0,-1){.034353}}
\multiput(33.02,58.72)(-.03331,-.037474){21}{\line(0,-1){.037474}}
\multiput(32.32,57.93)(-.033231,-.040832){20}{\line(0,-1){.040832}}
\multiput(31.65,57.11)(-.033077,-.044462){19}{\line(0,-1){.044462}}
\multiput(31.02,56.27)(-.032841,-.048406){18}{\line(0,-1){.048406}}
\multiput(30.43,55.4)(-.03251,-.052718){17}{\line(0,-1){.052718}}
\multiput(29.88,54.5)(-.032073,-.057462){16}{\line(0,-1){.057462}}
\multiput(29.37,53.58)(-.031512,-.062723){15}{\line(0,-1){.062723}}
\multiput(28.89,52.64)(-.033178,-.073886){13}{\line(0,-1){.073886}}
\multiput(28.46,51.68)(-.03243,-.08153){12}{\line(0,-1){.08153}}
\multiput(28.07,50.7)(-.03147,-.0904){11}{\line(0,-1){.0904}}
\multiput(27.73,49.71)(-.03362,-.11205){9}{\line(0,-1){.11205}}
\multiput(27.43,48.7)(-.03231,-.12759){8}{\line(0,-1){.12759}}
\multiput(27.17,47.68)(-.03055,-.14728){7}{\line(0,-1){.14728}}
\multiput(26.95,46.65)(-.02813,-.17322){6}{\line(0,-1){.17322}}
\multiput(26.78,45.61)(-.0309,-.2614){4}{\line(0,-1){.2614}}
\put(26.66,44.56){\line(0,-1){1.05}}
\put(26.58,43.51){\line(0,-1){1.052}}
\put(26.55,42.46){\line(0,-1){1.053}}
\put(26.56,41.41){\line(0,-1){1.051}}
\put(26.62,40.36){\line(0,-1){1.048}}
\multiput(26.73,39.31)(.03014,-.20841){5}{\line(0,-1){.20841}}
\multiput(26.88,38.27)(.03264,-.17242){6}{\line(0,-1){.17242}}
\multiput(27.08,37.23)(.03009,-.12813){8}{\line(0,-1){.12813}}
\multiput(27.32,36.21)(.03167,-.11262){9}{\line(0,-1){.11262}}
\multiput(27.6,35.19)(.03289,-.10002){10}{\line(0,-1){.10002}}
\multiput(27.93,34.19)(.03101,-.08208){12}{\line(0,-1){.08208}}
\multiput(28.3,33.21)(.031889,-.074451){13}{\line(0,-1){.074451}}
\multiput(28.72,32.24)(.03259,-.06778){14}{\line(0,-1){.06778}}
\multiput(29.17,31.29)(.033141,-.061878){15}{\line(0,-1){.061878}}
\multiput(29.67,30.36)(.033563,-.056604){16}{\line(0,-1){.056604}}
\multiput(30.21,29.46)(.031994,-.048969){18}{\line(0,-1){.048969}}
\multiput(30.78,28.58)(.0323,-.04503){19}{\line(0,-1){.04503}}
\multiput(31.4,27.72)(.032517,-.041403){20}{\line(0,-1){.041403}}
\multiput(32.05,26.89)(.032654,-.038047){21}{\line(0,-1){.038047}}
\multiput(32.73,26.09)(.03272,-.034927){22}{\line(0,-1){.034927}}
\multiput(33.45,25.33)(.034208,-.033471){22}{\line(1,0){.034208}}
\multiput(34.21,24.59)(.037329,-.033473){21}{\line(1,0){.037329}}
\multiput(34.99,23.89)(.040687,-.033408){20}{\line(1,0){.040687}}
\multiput(35.8,23.22)(.044317,-.03327){19}{\line(1,0){.044317}}
\multiput(36.65,22.59)(.048263,-.033051){18}{\line(1,0){.048263}}
\multiput(37.51,21.99)(.052576,-.032739){17}{\line(1,0){.052576}}
\multiput(38.41,21.43)(.057322,-.032322){16}{\line(1,0){.057322}}
\multiput(39.33,20.92)(.062585,-.031785){15}{\line(1,0){.062585}}
\multiput(40.26,20.44)(.073741,-.033499){13}{\line(1,0){.073741}}
\multiput(41.22,20)(.08139,-.03278){12}{\line(1,0){.08139}}
\multiput(42.2,19.61)(.09026,-.03186){11}{\line(1,0){.09026}}
\multiput(43.19,19.26)(.10072,-.0307){10}{\line(1,0){.10072}}
\multiput(44.2,18.95)(.12744,-.03286){8}{\line(1,0){.12744}}
\multiput(45.22,18.69)(.14715,-.03118){7}{\line(1,0){.14715}}
\multiput(46.25,18.47)(.17309,-.02888){6}{\line(1,0){.17309}}
\multiput(47.29,18.3)(.2613,-.032){4}{\line(1,0){.2613}}
\put(48.33,18.17){\line(1,0){1.05}}
\put(49.38,18.09){\line(1,0){1.052}}
\put(50.43,18.05){\line(1,0){1.053}}
\put(51.49,18.06){\line(1,0){1.051}}
\put(52.54,18.12){\line(1,0){1.048}}
\multiput(53.59,18.22)(.20854,.02923){5}{\line(1,0){.20854}}
\multiput(54.63,18.36)(.17256,.0319){6}{\line(1,0){.17256}}
\multiput(55.67,18.55)(.12826,.02953){8}{\line(1,0){.12826}}
\multiput(56.69,18.79)(.11276,.03118){9}{\line(1,0){.11276}}
\multiput(57.71,19.07)(.10016,.03245){10}{\line(1,0){.10016}}
\multiput(58.71,19.4)(.08969,.03343){11}{\line(1,0){.08969}}
\multiput(59.69,19.76)(.074589,.031565){13}{\line(1,0){.074589}}
\multiput(60.66,20.17)(.06792,.032296){14}{\line(1,0){.06792}}
\multiput(61.61,20.63)(.062021,.032871){15}{\line(1,0){.062021}}
\multiput(62.55,21.12)(.056749,.033317){16}{\line(1,0){.056749}}
\multiput(63.45,21.65)(.051997,.033651){17}{\line(1,0){.051997}}
\multiput(64.34,22.22)(.04517,.032104){19}{\line(1,0){.04517}}
\multiput(65.2,22.83)(.041544,.032336){20}{\line(1,0){.041544}}
\multiput(66.03,23.48)(.038188,.032488){21}{\line(1,0){.038188}}
\multiput(66.83,24.16)(.035069,.032568){22}{\line(1,0){.035069}}
\multiput(67.6,24.88)(.033619,.034063){22}{\line(0,1){.034063}}
\multiput(68.34,25.63)(.033635,.037183){21}{\line(0,1){.037183}}
\multiput(69.05,26.41)(.033585,.040541){20}{\line(0,1){.040541}}
\multiput(69.72,27.22)(.033463,.044173){19}{\line(0,1){.044173}}
\multiput(70.35,28.06)(.03326,.048119){18}{\line(0,1){.048119}}
\multiput(70.95,28.93)(.032967,.052433){17}{\line(0,1){.052433}}
\multiput(71.51,29.82)(.032571,.057181){16}{\line(0,1){.057181}}
\multiput(72.03,30.73)(.032056,.062447){15}{\line(0,1){.062447}}
\multiput(72.51,31.67)(.031403,.068338){14}{\line(0,1){.068338}}
\multiput(72.95,32.63)(.03313,.08125){12}{\line(0,1){.08125}}
\multiput(73.35,33.6)(.03226,.09012){11}{\line(0,1){.09012}}
\multiput(73.71,34.59)(.03114,.10058){10}{\line(0,1){.10058}}
\multiput(74.02,35.6)(.03341,.1273){8}{\line(0,1){.1273}}
\multiput(74.28,36.62)(.03182,.14701){7}{\line(0,1){.14701}}
\multiput(74.51,37.65)(.02963,.17296){6}{\line(0,1){.17296}}
\multiput(74.69,38.68)(.0331,.2611){4}{\line(0,1){.2611}}
\put(74.82,39.73){\line(0,1){1.049}}
\put(74.9,40.78){\line(0,1){1.473}}
\put(91.88,41.25){\oval(24.75,50)[]}
\multiput(31.25,68.25)(.03361345,-.05882353){119}{\line(0,-1){.05882353}}
\put(35.25,61.25){\line(1,0){25.75}}
\multiput(61,61.25)(.033737024,-.109861592){289}{\line(0,-1){.109861592}}
\multiput(70.75,29.5)(.291667,.033333){30}{\line(1,0){.291667}}
\multiput(79.5,30.5)(.0335821,.4477612){67}{\line(0,1){.4477612}}
\put(81.75,60.5){\line(1,0){19.5}}
\multiput(101.25,60.5)(.0333333,-.5066667){75}{\line(0,-1){.5066667}}
\put(103.75,22.5){\line(4,-3){4}}
\put(36.5,58.25){\makebox(0,0)[cc]{$x_1$}}
\put(27,70.25){\makebox(0,0)[cc]{$x$}}
\put(26.25,16){\makebox(0,0)[cc]{$y$}}
\put(116,13.25){\makebox(0,0)[cc]{$z$}}
\put(112,69){\makebox(0,0)[cc]{$w$}}
\put(25.75,60.25){\makebox(0,0)[cc]{$x_2$}}
\put(27.5,26.5){\makebox(0,0)[cc]{$y_1$}}
\put(41,12.25){\makebox(0,0)[cc]{$y_2$}}
\put(64.25,14.5){\makebox(0,0)[cc]{$z_1$}}
\put(43.5,39.75){\makebox(0,0)[cc]{$\onn_\sigma ( A )$}}
\put(93,28.25){\makebox(0,0)[cc]{$\onn_\kappa ( B )$}}
\put(72.5,50){\makebox(0,0)[cc]{$\onn_\kappa ( B' ) $}}
\put(55.5,58.5){\makebox(0,0)[cc]{$b_6$}}
\put(65.75,31.25){\makebox(0,0)[cc]{$z_2$}}
\put(84.25,28.75){\makebox(0,0)[cc]{$b_4$}}
\put(84.75,61.75){\makebox(0,0)[cc]{$b_7$}}
\put(97.5,56.75){\makebox(0,0)[cc]{$b_3$}}
\put(99.75,72.25){\makebox(0,0)[cc]{$b_2$}}
\put(111.25,61){\makebox(0,0)[cc]{$b_1$}}
\put(100.5,24){\makebox(0,0)[cc]{$b_5$}}
\end{picture}

\caption{Case 2.B in proof of Lemma \ref{4fat}.} \label{caz2bl4fat}
\end{figure}

The point $b_6$ is the entrance point of $[x,b_3]$ into
$\onn_\kappa (B)$ while $b_4$ is the entrance point of $[z_2,z]$
into $\onn_\kappa (B)$. If $b_6\in [z_2,b_3]$ then $b_6=b_4$ and
$\dist ([y,z], [x,w])\leq 2\lambda$. For $\theta$ large enough
this case is impossible. Thus $b_6\in [x,z_2]$.

The intersection $[x_1,z_2]\cap [b_6,z_2]$ is in $\onn_{t\sigma}
(A) \cap \onn_{t\kappa }(B)$. Note that $\dist (b_6, z_2)\geq
\theta - \lambda - \delta$, thus for $\theta$ large we may assume
that $\dist (b_6, z_2)> \diam_\tau +1$ with $\tau = t \max (\sigma
, \kappa )+1$.

If $\dist (x_1,z_2) \leq \diam_\tau +1$ then $\{ x_2,x_1, z_1\}$
has diameter $O(1)$ and we are back in Case 1 with $a_1=x_2,
a_2=z_1$ and $a_3=x_1$ and with the constant $\sigma$ possibly
larger. We may then use the proof in Case 1 to finish the
argument.

Assume now that $\dist (x_1,z_2) > \diam_\tau +1$. Then $A=B$ and
 $b_5$, the entrance point of $[z,z_2]$ into $\onn_\kappa (B)$, is
 also the entrance point of $[z,x]$ into $\onn_\kappa (A)$. As
 $z_2$ is the entrance point of $[z,x]$ into $\onn_\sigma (A)$,
 Lemma \ref{doie} implies that $\dist (z_2, b_5)=O(1)$.

 By construction $b_3\in [z_2, b_5]$, hence $\dist (b_3, z_2)=O(1)$.
 On the other hand $\dist(b_3, z_2)\geq \dist (w,[y,z])-2\sigma -2\theta -\delta
 \geq 6\theta -2\sigma - \delta $. Thus for $\theta$ large
 enough we obtain a contradiction.

\me

\begin{figure}[!ht]
\centering
\unitlength 1mm 
\linethickness{0.4pt}
\ifx\plotpoint\undefined\newsavebox{\plotpoint}\fi 
\begin{picture}(104.5,76.25)(0,0)
\put(20.75,21.25){\framebox(79.5,51)[cc]{}}
\put(64.25,53.75){\oval(21,44.5)[]}
\put(101,62.75){\line(0,1){.575}}
\put(100.99,63.32){\line(0,1){.573}}
\put(100.94,63.9){\line(0,1){.57}}
\put(100.86,64.47){\line(0,1){.565}}
\multiput(100.76,65.03)(-.0274,.11167){5}{\line(0,1){.11167}}
\multiput(100.62,65.59)(-.03333,.11004){5}{\line(0,1){.11004}}
\multiput(100.45,66.14)(-.03264,.09009){6}{\line(0,1){.09009}}
\multiput(100.26,66.68)(-.03206,.07561){7}{\line(0,1){.07561}}
\multiput(100.03,67.21)(-.03155,.06457){8}{\line(0,1){.06457}}
\multiput(99.78,67.73)(-.03107,.05581){9}{\line(0,1){.05581}}
\multiput(99.5,68.23)(-.03061,.04866){10}{\line(0,1){.04866}}
\multiput(99.2,68.72)(-.03317,.04696){10}{\line(0,1){.04696}}
\multiput(98.86,69.19)(-.03239,.04102){11}{\line(0,1){.04102}}
\multiput(98.51,69.64)(-.03166,.03596){12}{\line(0,1){.03596}}
\multiput(98.13,70.07)(-.03353,.03422){12}{\line(0,1){.03422}}
\multiput(97.73,70.48)(-.03531,.03237){12}{\line(-1,0){.03531}}
\multiput(97.3,70.87)(-.04036,.03321){11}{\line(-1,0){.04036}}
\multiput(96.86,71.23)(-.04207,.031){11}{\line(-1,0){.04207}}
\multiput(96.4,71.57)(-.04804,.03158){10}{\line(-1,0){.04804}}
\multiput(95.91,71.89)(-.05518,.03219){9}{\line(-1,0){.05518}}
\multiput(95.42,72.18)(-.06392,.03284){8}{\line(-1,0){.06392}}
\multiput(94.91,72.44)(-.07495,.03358){7}{\line(-1,0){.07495}}
\multiput(94.38,72.68)(-.07664,.02952){7}{\line(-1,0){.07664}}
\multiput(93.85,72.88)(-.09113,.02961){6}{\line(-1,0){.09113}}
\multiput(93.3,73.06)(-.11109,.02964){5}{\line(-1,0){.11109}}
\put(92.74,73.21){\line(-1,0){.563}}
\put(92.18,73.33){\line(-1,0){.568}}
\put(91.61,73.42){\line(-1,0){.572}}
\put(91.04,73.47){\line(-1,0){.574}}
\put(90.47,73.5){\line(-1,0){.575}}
\put(89.89,73.5){\line(-1,0){.574}}
\put(89.32,73.46){\line(-1,0){.571}}
\put(88.75,73.4){\line(-1,0){.567}}
\multiput(88.18,73.3)(-.1402,-.0314){4}{\line(-1,0){.1402}}
\multiput(87.62,73.18)(-.11069,-.03111){5}{\line(-1,0){.11069}}
\multiput(87.07,73.02)(-.09073,-.03082){6}{\line(-1,0){.09073}}
\multiput(86.52,72.84)(-.07624,-.03053){7}{\line(-1,0){.07624}}
\multiput(85.99,72.62)(-.06519,-.03024){8}{\line(-1,0){.06519}}
\multiput(85.47,72.38)(-.06348,-.03368){8}{\line(-1,0){.06348}}
\multiput(84.96,72.11)(-.05475,-.03291){9}{\line(-1,0){.05475}}
\multiput(84.47,71.81)(-.04762,-.03221){10}{\line(-1,0){.04762}}
\multiput(83.99,71.49)(-.04166,-.03156){11}{\line(-1,0){.04166}}
\multiput(83.53,71.14)(-.03992,-.03374){11}{\line(-1,0){.03992}}
\multiput(83.09,70.77)(-.03488,-.03284){12}{\line(-1,0){.03488}}
\multiput(82.67,70.38)(-.03308,-.03466){12}{\line(0,-1){.03466}}
\multiput(82.28,69.96)(-.03118,-.03637){12}{\line(0,-1){.03637}}
\multiput(81.9,69.53)(-.03184,-.04144){11}{\line(0,-1){.04144}}
\multiput(81.55,69.07)(-.03254,-.04739){10}{\line(0,-1){.04739}}
\multiput(81.23,68.6)(-.03329,-.05452){9}{\line(0,-1){.05452}}
\multiput(80.93,68.11)(-.03033,-.05622){9}{\line(0,-1){.05622}}
\multiput(80.65,67.6)(-.03069,-.06498){8}{\line(0,-1){.06498}}
\multiput(80.41,67.08)(-.03106,-.07603){7}{\line(0,-1){.07603}}
\multiput(80.19,66.55)(-.03144,-.09051){6}{\line(0,-1){.09051}}
\multiput(80,66.01)(-.03187,-.11047){5}{\line(0,-1){.11047}}
\multiput(79.84,65.45)(-.0324,-.14){4}{\line(0,-1){.14}}
\put(79.71,64.89){\line(0,-1){.566}}
\put(79.61,64.33){\line(0,-1){.571}}
\put(79.54,63.76){\line(0,-1){1.148}}
\put(79.5,62.61){\line(0,-1){.574}}
\put(79.52,62.03){\line(0,-1){.572}}
\put(79.57,61.46){\line(0,-1){.569}}
\put(79.66,60.89){\line(0,-1){.563}}
\multiput(79.77,60.33)(.02887,-.1113){5}{\line(0,-1){.1113}}
\multiput(79.92,59.77)(.02898,-.09133){6}{\line(0,-1){.09133}}
\multiput(80.09,59.22)(.02899,-.07684){7}{\line(0,-1){.07684}}
\multiput(80.29,58.69)(.03306,-.07518){7}{\line(0,-1){.07518}}
\multiput(80.53,58.16)(.0324,-.06414){8}{\line(0,-1){.06414}}
\multiput(80.78,57.65)(.03181,-.0554){9}{\line(0,-1){.0554}}
\multiput(81.07,57.15)(.03125,-.04826){10}{\line(0,-1){.04826}}
\multiput(81.38,56.67)(.03071,-.04229){11}{\line(0,-1){.04229}}
\multiput(81.72,56.2)(.03293,-.04059){11}{\line(0,-1){.04059}}
\multiput(82.08,55.75)(.03213,-.03554){12}{\line(0,-1){.03554}}
\multiput(82.47,55.33)(.031369,-.031172){13}{\line(1,0){.031369}}
\multiput(82.88,54.92)(.03574,-.03191){12}{\line(1,0){.03574}}
\multiput(83.31,54.54)(.04079,-.03267){11}{\line(1,0){.04079}}
\multiput(83.75,54.18)(.04673,-.03349){10}{\line(1,0){.04673}}
\multiput(84.22,53.85)(.04845,-.03094){10}{\line(1,0){.04845}}
\multiput(84.71,53.54)(.0556,-.03146){9}{\line(1,0){.0556}}
\multiput(85.21,53.25)(.06435,-.03199){8}{\line(1,0){.06435}}
\multiput(85.72,53)(.07539,-.03258){7}{\line(1,0){.07539}}
\multiput(86.25,52.77)(.08986,-.03326){6}{\line(1,0){.08986}}
\multiput(86.79,52.57)(.09151,-.02841){6}{\line(1,0){.09151}}
\multiput(87.34,52.4)(.11148,-.02817){5}{\line(1,0){.11148}}
\put(87.89,52.26){\line(1,0){.564}}
\put(88.46,52.15){\line(1,0){.569}}
\put(89.03,52.07){\line(1,0){.573}}
\put(89.6,52.02){\line(1,0){1.149}}
\put(90.75,52.01){\line(1,0){.573}}
\put(91.32,52.05){\line(1,0){.57}}
\put(91.89,52.12){\line(1,0){.566}}
\multiput(92.46,52.23)(.1398,.0333){4}{\line(1,0){.1398}}
\multiput(93.02,52.36)(.11027,.03257){5}{\line(1,0){.11027}}
\multiput(93.57,52.52)(.09031,.03201){6}{\line(1,0){.09031}}
\multiput(94.11,52.71)(.07583,.03154){7}{\line(1,0){.07583}}
\multiput(94.64,52.94)(.06478,.0311){8}{\line(1,0){.06478}}
\multiput(95.16,53.18)(.05603,.03068){9}{\line(1,0){.05603}}
\multiput(95.67,53.46)(.05431,.03363){9}{\line(1,0){.05431}}
\multiput(96.15,53.76)(.04719,.03284){10}{\line(1,0){.04719}}
\multiput(96.63,54.09)(.04124,.0321){11}{\line(1,0){.04124}}
\multiput(97.08,54.44)(.03618,.03141){12}{\line(1,0){.03618}}
\multiput(97.51,54.82)(.03445,.0333){12}{\line(1,0){.03445}}
\multiput(97.93,55.22)(.03262,.03509){12}{\line(0,1){.03509}}
\multiput(98.32,55.64)(.03349,.04013){11}{\line(0,1){.04013}}
\multiput(98.69,56.08)(.03129,.04186){11}{\line(0,1){.04186}}
\multiput(99.03,56.54)(.03191,.04782){10}{\line(0,1){.04782}}
\multiput(99.35,57.02)(.03257,.05495){9}{\line(0,1){.05495}}
\multiput(99.64,57.52)(.03328,.06369){8}{\line(0,1){.06369}}
\multiput(99.91,58.03)(.02983,.06538){8}{\line(0,1){.06538}}
\multiput(100.15,58.55)(.03005,.07643){7}{\line(0,1){.07643}}
\multiput(100.36,59.08)(.03024,.09092){6}{\line(0,1){.09092}}
\multiput(100.54,59.63)(.03041,.11089){5}{\line(0,1){.11089}}
\multiput(100.69,60.18)(.0306,.1404){4}{\line(0,1){.1404}}
\put(100.81,60.75){\line(0,1){.567}}
\put(100.91,61.31){\line(0,1){.572}}
\put(100.97,61.89){\line(0,1){.865}}
\put(66.01,42.75){\line(0,1){1.04}}
\put(65.99,43.79){\line(0,1){1.038}}
\put(65.92,44.83){\line(0,1){1.034}}
\multiput(65.81,45.86)(-.03175,.2056){5}{\line(0,1){.2056}}
\multiput(65.65,46.89)(-.02908,.14572){7}{\line(0,1){.14572}}
\multiput(65.44,47.91)(-.031,.12627){8}{\line(0,1){.12627}}
\multiput(65.2,48.92)(-.03245,.11093){9}{\line(0,1){.11093}}
\multiput(64.9,49.92)(-.03354,.09846){10}{\line(0,1){.09846}}
\multiput(64.57,50.9)(-.03152,.08075){12}{\line(0,1){.08075}}
\multiput(64.19,51.87)(-.032326,.073192){13}{\line(0,1){.073192}}
\multiput(63.77,52.82)(-.032962,.066585){14}{\line(0,1){.066585}}
\multiput(63.31,53.76)(-.033455,.06074){15}{\line(0,1){.06074}}
\multiput(62.81,54.67)(-.031836,.052251){17}{\line(0,1){.052251}}
\multiput(62.27,55.55)(-.032199,.047985){18}{\line(0,1){.047985}}
\multiput(61.69,56.42)(-.032464,.044081){19}{\line(0,1){.044081}}
\multiput(61.07,57.26)(-.032644,.040487){20}{\line(0,1){.040487}}
\multiput(60.42,58.07)(-.032747,.037162){21}{\line(0,1){.037162}}
\multiput(59.73,58.85)(-.032781,.034071){22}{\line(0,1){.034071}}
\multiput(59.01,59.6)(-.034241,.032603){22}{\line(-1,0){.034241}}
\multiput(58.26,60.31)(-.037332,.032553){21}{\line(-1,0){.037332}}
\multiput(57.47,61)(-.040657,.032433){20}{\line(-1,0){.040657}}
\multiput(56.66,61.65)(-.04425,.032234){19}{\line(-1,0){.04425}}
\multiput(55.82,62.26)(-.048152,.031948){18}{\line(-1,0){.048152}}
\multiput(54.95,62.83)(-.055692,.033536){16}{\line(-1,0){.055692}}
\multiput(54.06,63.37)(-.060914,.033138){15}{\line(-1,0){.060914}}
\multiput(53.15,63.87)(-.066756,.032615){14}{\line(-1,0){.066756}}
\multiput(52.21,64.32)(-.07336,.031944){13}{\line(-1,0){.07336}}
\multiput(51.26,64.74)(-.08091,.03109){12}{\line(-1,0){.08091}}
\multiput(50.29,65.11)(-.09863,.03303){10}{\line(-1,0){.09863}}
\multiput(49.3,65.44)(-.11109,.03187){9}{\line(-1,0){.11109}}
\multiput(48.3,65.73)(-.12643,.03035){8}{\line(-1,0){.12643}}
\multiput(47.29,65.97)(-.17018,.03305){6}{\line(-1,0){.17018}}
\multiput(46.27,66.17)(-.20576,.03068){5}{\line(-1,0){.20576}}
\put(45.24,66.32){\line(-1,0){1.035}}
\put(44.2,66.43){\line(-1,0){1.038}}
\put(43.17,66.49){\line(-1,0){1.04}}
\put(42.13,66.51){\line(-1,0){1.04}}
\put(41.09,66.48){\line(-1,0){1.038}}
\put(40.05,66.41){\line(-1,0){1.033}}
\multiput(39.02,66.29)(-.20543,-.03282){5}{\line(-1,0){.20543}}
\multiput(37.99,66.13)(-.14557,-.02984){7}{\line(-1,0){.14557}}
\multiput(36.97,65.92)(-.12611,-.03166){8}{\line(-1,0){.12611}}
\multiput(35.96,65.66)(-.11076,-.03302){9}{\line(-1,0){.11076}}
\multiput(34.96,65.37)(-.08935,-.03096){11}{\line(-1,0){.08935}}
\multiput(33.98,65.03)(-.08058,-.03194){12}{\line(-1,0){.08058}}
\multiput(33.01,64.64)(-.073023,-.032707){13}{\line(-1,0){.073023}}
\multiput(32.06,64.22)(-.066413,-.033309){14}{\line(-1,0){.066413}}
\multiput(31.13,63.75)(-.05678,-.03166){16}{\line(-1,0){.05678}}
\multiput(30.23,63.25)(-.052085,-.032108){17}{\line(-1,0){.052085}}
\multiput(29.34,62.7)(-.047817,-.032448){18}{\line(-1,0){.047817}}
\multiput(28.48,62.12)(-.043911,-.032694){19}{\line(-1,0){.043911}}
\multiput(27.65,61.49)(-.040316,-.032855){20}{\line(-1,0){.040316}}
\multiput(26.84,60.84)(-.03699,-.032941){21}{\line(-1,0){.03699}}
\multiput(26.06,60.15)(-.033899,-.032958){22}{\line(-1,0){.033899}}
\multiput(25.32,59.42)(-.032424,-.03441){22}{\line(0,-1){.03441}}
\multiput(24.6,58.66)(-.032358,-.037501){21}{\line(0,-1){.037501}}
\multiput(23.92,57.88)(-.03222,-.040825){20}{\line(0,-1){.040825}}
\multiput(23.28,57.06)(-.032003,-.044417){19}{\line(0,-1){.044417}}
\multiput(22.67,56.22)(-.033561,-.051161){17}{\line(0,-1){.051161}}
\multiput(22.1,55.35)(-.033245,-.055867){16}{\line(0,-1){.055867}}
\multiput(21.57,54.45)(-.03282,-.061086){15}{\line(0,-1){.061086}}
\multiput(21.08,53.54)(-.032266,-.066925){14}{\line(0,-1){.066925}}
\multiput(20.63,52.6)(-.031561,-.073525){13}{\line(0,-1){.073525}}
\multiput(20.21,51.64)(-.03346,-.08844){11}{\line(0,-1){.08844}}
\multiput(19.85,50.67)(-.03251,-.0988){10}{\line(0,-1){.0988}}
\multiput(19.52,49.68)(-.03129,-.11126){9}{\line(0,-1){.11126}}
\multiput(19.24,48.68)(-.02969,-.12659){8}{\line(0,-1){.12659}}
\multiput(19,47.67)(-.03216,-.17035){6}{\line(0,-1){.17035}}
\multiput(18.81,46.65)(-.02961,-.20592){5}{\line(0,-1){.20592}}
\put(18.66,45.62){\line(0,-1){1.035}}
\put(18.56,44.58){\line(0,-1){1.039}}
\put(18.5,43.54){\line(0,-1){1.04}}
\put(18.49,42.5){\line(0,-1){1.04}}
\put(18.52,41.46){\line(0,-1){1.037}}
\multiput(18.6,40.43)(.0311,-.2582){4}{\line(0,-1){.2582}}
\multiput(18.73,39.39)(.02824,-.17104){6}{\line(0,-1){.17104}}
\multiput(18.9,38.37)(.0306,-.14541){7}{\line(0,-1){.14541}}
\multiput(19.11,37.35)(.03232,-.12594){8}{\line(0,-1){.12594}}
\multiput(19.37,36.34)(.0336,-.11058){9}{\line(0,-1){.11058}}
\multiput(19.67,35.35)(.03142,-.08919){11}{\line(0,-1){.08919}}
\multiput(20.02,34.36)(.03236,-.08042){12}{\line(0,-1){.08042}}
\multiput(20.41,33.4)(.033087,-.072851){13}{\line(0,-1){.072851}}
\multiput(20.84,32.45)(.033655,-.066238){14}{\line(0,-1){.066238}}
\multiput(21.31,31.53)(.031956,-.056614){16}{\line(0,-1){.056614}}
\multiput(21.82,30.62)(.032379,-.051916){17}{\line(0,-1){.051916}}
\multiput(22.37,29.74)(.032697,-.047647){18}{\line(0,-1){.047647}}
\multiput(22.96,28.88)(.032922,-.04374){19}{\line(0,-1){.04374}}
\multiput(23.58,28.05)(.033064,-.040145){20}{\line(0,-1){.040145}}
\multiput(24.24,27.25)(.033133,-.036818){21}{\line(0,-1){.036818}}
\multiput(24.94,26.47)(.033135,-.033727){22}{\line(0,-1){.033727}}
\multiput(25.67,25.73)(.034579,-.032244){22}{\line(1,0){.034579}}
\multiput(26.43,25.02)(.037669,-.032162){21}{\line(1,0){.037669}}
\multiput(27.22,24.35)(.04315,-.033692){19}{\line(1,0){.04315}}
\multiput(28.04,23.71)(.04706,-.033536){18}{\line(1,0){.04706}}
\multiput(28.89,23.1)(.051335,-.033294){17}{\line(1,0){.051335}}
\multiput(29.76,22.54)(.056039,-.032954){16}{\line(1,0){.056039}}
\multiput(30.66,22.01)(.061256,-.032501){15}{\line(1,0){.061256}}
\multiput(31.58,21.52)(.067093,-.031917){14}{\line(1,0){.067093}}
\multiput(32.51,21.07)(.073689,-.031177){13}{\line(1,0){.073689}}
\multiput(33.47,20.67)(.08862,-.033){11}{\line(1,0){.08862}}
\multiput(34.45,20.31)(.09897,-.032){10}{\line(1,0){.09897}}
\multiput(35.44,19.99)(.11142,-.03071){9}{\line(1,0){.11142}}
\multiput(36.44,19.71)(.14484,-.03317){7}{\line(1,0){.14484}}
\multiput(37.45,19.48)(.17052,-.03127){6}{\line(1,0){.17052}}
\multiput(38.48,19.29)(.20607,-.02853){5}{\line(1,0){.20607}}
\put(39.51,19.15){\line(1,0){1.036}}
\put(40.54,19.05){\line(1,0){1.039}}
\put(41.58,19){\line(1,0){1.04}}
\put(42.62,18.99){\line(1,0){1.039}}
\put(43.66,19.03){\line(1,0){1.037}}
\multiput(44.7,19.11)(.258,.0324){4}{\line(1,0){.258}}
\multiput(45.73,19.24)(.17089,.02914){6}{\line(1,0){.17089}}
\multiput(46.76,19.42)(.14525,.03136){7}{\line(1,0){.14525}}
\multiput(47.77,19.64)(.12577,.03298){8}{\line(1,0){.12577}}
\multiput(48.78,19.9)(.09936,.03076){10}{\line(1,0){.09936}}
\multiput(49.77,20.21)(.08902,.03189){11}{\line(1,0){.08902}}
\multiput(50.75,20.56)(.08025,.03277){12}{\line(1,0){.08025}}
\multiput(51.71,20.95)(.072678,.033466){13}{\line(1,0){.072678}}
\multiput(52.66,21.39)(.061658,.031733){15}{\line(1,0){.061658}}
\multiput(53.58,21.87)(.056447,.032251){16}{\line(1,0){.056447}}
\multiput(54.49,22.38)(.051747,.03265){17}{\line(1,0){.051747}}
\multiput(55.37,22.94)(.047476,.032945){18}{\line(1,0){.047476}}
\multiput(56.22,23.53)(.043568,.03315){19}{\line(1,0){.043568}}
\multiput(57.05,24.16)(.039972,.033273){20}{\line(1,0){.039972}}
\multiput(57.85,24.82)(.036645,.033324){21}{\line(1,0){.036645}}
\multiput(58.62,25.52)(.033554,.03331){22}{\line(1,0){.033554}}
\multiput(59.36,26.26)(.033591,.036401){21}{\line(0,1){.036401}}
\multiput(60.06,27.02)(.033564,.039728){20}{\line(0,1){.039728}}
\multiput(60.73,27.82)(.033466,.043325){19}{\line(0,1){.043325}}
\multiput(61.37,28.64)(.03329,.047234){18}{\line(0,1){.047234}}
\multiput(61.97,29.49)(.033026,.051508){17}{\line(0,1){.051508}}
\multiput(62.53,30.37)(.032661,.05621){16}{\line(0,1){.05621}}
\multiput(63.05,31.26)(.032181,.061425){15}{\line(0,1){.061425}}
\multiput(63.53,32.19)(.031567,.067258){14}{\line(0,1){.067258}}
\multiput(63.98,33.13)(.03336,.08){12}{\line(0,1){.08}}
\multiput(64.38,34.09)(.03254,.08879){11}{\line(0,1){.08879}}
\multiput(64.73,35.06)(.03148,.09914){10}{\line(0,1){.09914}}
\multiput(65.05,36.06)(.03013,.11158){9}{\line(0,1){.11158}}
\multiput(65.32,37.06)(.03242,.14502){7}{\line(0,1){.14502}}
\multiput(65.55,38.07)(.03038,.17068){6}{\line(0,1){.17068}}
\multiput(65.73,39.1)(.02746,.20621){5}{\line(0,1){.20621}}
\put(65.87,40.13){\line(0,1){1.036}}
\put(65.96,41.17){\line(0,1){1.584}}
\multiput(21,72.25)(.03368794,-.08687943){141}{\line(0,-1){.08687943}}
\multiput(25.75,60)(.6222222,.0333333){45}{\line(1,0){.6222222}}
\put(53.75,61.5){\line(0,1){0}}
\multiput(53.75,61.5)(.033717105,-.092105263){304}{\line(0,-1){.092105263}}
\multiput(64,33.5)(.03369906,.096394984){319}{\line(0,1){.096394984}}
\multiput(74.75,64.25)(.0783582,.0335821){67}{\line(1,0){.0783582}}
\put(80,66.5){\line(0,1){0}}
\multiput(80,66.5)(.03553616,-.033665835){401}{\line(1,0){.03553616}}
\multiput(94.25,53)(.03362573,-.18421053){171}{\line(0,-1){.18421053}}
\put(100,21.5){\line(0,1){0}}
\put(17.5,73.5){\makebox(0,0)[cc]{$x$}}
\put(17.75,17){\makebox(0,0)[cc]{$y$}}
\put(102.75,17.75){\makebox(0,0)[cc]{$z$}}
\put(104.5,73.75){\makebox(0,0)[cc]{$w$}}
\put(37.75,40){\makebox(0,0)[cc]{$\onn_\sigma (A_1)$}}
\put(65,68.25){\makebox(0,0)[cc]{$\onn_\kappa (A_3)$}}
\put(92.75,66.25){\makebox(0,0)[cc]{$\onn_\sigma ( A_2)$}}
\put(26.75,57){\makebox(0,0)[cc]{$x_1$}}
\put(16.25,53.75){\makebox(0,0)[cc]{$x_2$}}
\put(16.25,31.25){\makebox(0,0)[cc]{$y_1$}}
\put(31,16.5){\makebox(0,0)[cc]{$y_2$}}
\put(53.75,17.5){\makebox(0,0)[cc]{$z_1$}}
\put(49.75,58){\makebox(0,0)[cc]{$z_3$}}
\put(68.5,36.25){\makebox(0,0)[cc]{$z_2$}}
\put(83.5,76.25){\makebox(0,0)[cc]{$x_1'$}}
\put(103,66.5){\makebox(0,0)[cc]{$w_1$}}
\put(95.25,76.25){\makebox(0,0)[cc]{$w_2$}}
\put(91.5,48.25){\makebox(0,0)[cc]{$z_1'$}}
\put(103.25,56.25){\makebox(0,0)[cc]{$z_2'$}}
\put(71.25,64.5){\makebox(0,0)[cc]{$z_4$}}
\put(84,66.5){\makebox(0,0)[cc]{$x_2'$}}
\end{picture}

\caption{Case 3 in proof of Lemma \ref{4fat}.} \label{caz3l4fat}
\end{figure}
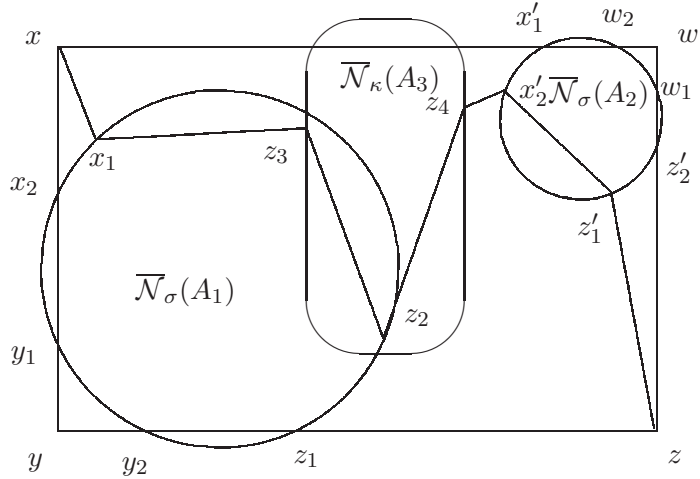

\textsc{Case 3.}\quad Assume that both geodesic triangles $xyz$
and $xzw$ are in case \bp. Then there exists $A_1$ in $\aaa$ such
that $\onn_\sigma (A_1)$ intersects all the edges of $xyz$.
Moreover the pairs of entrance points in $\onn_\sigma (A_1)$,
$(x_1,x_2)$, $(y_1,y_2)$ and $(z_1,z_2)$ are all at respective
distances less than $\delta$. Likewise there exists $A_2$ in
$\aaa$ such that $\onn_\sigma (A_2)$ intersects all the edges of
$xzw$, and the pairs of entrance points in $\onn_\sigma (A_2)$,
$(x'_1,x'_2)$, $(z_1',z_2')$ and $(w_1,w_2)$ are all at distances
less than $\delta$ (see Figure \ref{caz3l4fat}).

If $\theta > \delta$ then $y_1,y_2\in B(y, 2\theta)$ and $w_1,w_2
\in B(w, 2\theta)$.

If $A_1=A_2=A$ then $x_1=x_2'$,  $z_2=z_1'$, hence $\dist
(x_1',x_2)$ and $\dist (z_1,z_2')$ are less than $2\delta$. If
$\theta >2\delta$ it follows that $x_1',x_2 \in B(x,2\theta)$ and
that $z_1,z_2' \in B(z, 2\theta)$. Thus $x,y,z,w \in \onn_{2\theta
+\sigma} (A)$, which by $(\fq)$ implies that $P\subset \nn_{\chi}
(A)$ for $\chi> t(2\theta +\sigma)$.

Assume that $A_1\neq A_2$. Then $[x_1,z_2]$ and $[x_2',z_1']$ are
either disjoint or they intersect in a sub-geodesic of length at
most $\diam_{t\sigma +1}$.

Let $\tau = t\max (\sigma , \kappa )+1$. If  either $[x_1,z_2]$ or
$[x_2', z_1']$ is of length at most $\diam_\tau+1$ then either $\{
x_2,x_1, z_1\}$ or $\{ z_1', z_2', x_1' \}$ is of diameter at most
$\diam_{\tau} +1 +2\delta $. Therefore we find ourselves in Case
2.B, with the constant $\sigma$ possibly larger. We can then
finish the argument as in that case.

If both $[x_1,z_2]$ and $[x_2', z_1']$ have length larger than
$\diam_\tau+1$ and their intersection is non-empty then either
$x_1'$ and $z_1$ or $x_2$ and $z_2'$ are at distance at most
$2\delta +\diam_{t\sigma +1}$. If $\theta > 2\delta
+\diam_{t\sigma +1}$ then this is impossible. We may therefore
assume that $[x_1,z_2]$ and $[x_2',z_1']$ do not intersect, and
are both of length larger than $\diam_\tau+1$.

Without loss of generality we may also assume that
$[x_1,z_2]\subset [x,x_2')$. Lemma \ref{distc} applied to the
geodesic $[x,x_1']$ and the point $z_2\in [x,x_2']$ implies that
either $z_2\in \onn_\kappa ([x,w])$ or $z_2\in \onn_\kappa (A_3)$,
where $\onn_\kappa (A_3)$ intersects $[x,x_1']$, and the entrance
and exit points $z_3 , z_4$ of $[x, x_2']$ in $\onn_\kappa (A_3)$
are at distance at most $\lambda$ from $[x,x_1']$ (see Figure
\ref{caz3l4fat}). If $\theta
>\delta +\kappa$ then the first case cannot occur. In the second case
we have that $\dist \left( z_2, \{ z_3, z_4\} \right) > \theta
-\delta -\lambda$. In particular we may assume that $\dist (z_2,
z_3) > \diam_\tau+1$. We also have the assumption that $\dist
(x_1, z_2)> \diam_\tau +1$. Then $[x_1, z_2]\cap
 [z_3,z_2]$ has diameter $> \diam_\tau$ and it is contained in
 $\onn_{t\sigma} (A_1)\cap \onn_{t\kappa } (A_3)$. Therefore
 $A_1=A_3$. In particular $z_4$, the exit point of $[x, x_2']$
  from $\onn_\kappa (A_1)$, and $z_2$, the exit point of $[x, z]$
  (therefore also of $[x, x_2']$) from $\onn_\sigma (A_1)$  are at distance $O(1)$
 by Lemma \ref{doie}. It follows that $z_1$ and $[x,w]$ are at
 distance $O(1)$, and if $\theta$ is large enough this
 gives a contradiction.\endproof

\begin{cor}\label{pi2}
Let $(X, \dist )$ and $\aaa$ satisfy the hypotheses of Lemma
\ref{4fat}. Then in any asymptotic cone of $X$ property $(\Pi_2)$
is satisfied by the collection of limit sets $\aau$.
\end{cor}

\proof In an asymptotic cone $\co{X;e,d}$ consider a simple bigon
of endpoints $x,y$ whose edges are limit geodesics. Then there
exist two sequences of geodesics $[x_n,y_n]$ and $[x_n',y_n']$
such that \uass $\dist (x_n,x_n')$ and $\dist (y_n,y_n')$ are of
order $o(d_n)$, while $\dist (x_n,y_n)$ and $\dist (x_n',y_n')$
are of order $O(d_n)$. Let $m_n$ and $m_n'$ be the middlepoints of
$[x_n,y_n]$ and respectively of $[x_n',y_n']$. Let $\delta_n $ be
the maximum between $\dist ([x_n,m_n]\, ,\, [x_n',m_n'])$ and
$\nu\theta$, where $\theta$ and $\nu$ are the constants provided
by Lemma \ref{4fat}. Then $\delta_n=_\omega o(d_n)$. Similarly,
$\delta_n'=\max\left\{ \dist ([m_n,y_n]\, ,\, [m_n',y_n'])\, ,\,
\nu\theta \right\}$ satisfies $\delta_n'=_\omega o(d_n)$. Let
$x_n^1$ be the farthest from $x_n$ point on $[x_n,m_n]$ at
distance $\delta_n$ from $[x_n',m_n']$, and let $x_n^2$ be the
farthest from $x_n'$ point on $[x_n',m_n']$ at distance $\delta_n$
from $x_n^1$. We choose in a similar manner $y_n^1\in [y_n,m_n]$
and $y_n^2\in [y_n',m_n']$. Since the limit bigon is simple, it
follows that the sets $\{x_n,x_n', x_n^1,x_n^2\}$ and $\{y_n,y_n',
y_n^1,y_n^2\}$ have diameters of order $o(d_n)$ \uas.

We prove that any quadrilateral having as two opposite edges
$[x_n^1,y_n^1]\subset [x_n,y_n]$ and $[x_n^2,y_n^2]\subset
[x_n',y_n']$ is $(\theta , \nu)$--fat. This suffices to finish the
argument, by Lemma \ref{4fat}.

\medskip

$\mathbf{(F_1)}$\quad By construction $\dist ([x_n^1,y_n^1] ,
[x_n^2,y_n^2])\geq \nu \theta$, while the edges $[x_n^1,x_n^2]$
and $[y_n^1,y_n^2]$ are at distance $O(d_n)$ from each other. The
rest of the property follows by Lemma \ref{simplif}.

\medskip

Property $\mathbf{(F_2)}$ follows immediately from the fact that
$\dist ([x_n^1,y_n^1] , [x_n^2,y_n^2])\geq \nu \theta$ and that
$\dist ([x_n^1,x_n^2],[y_n^1,y_n^2])=O(d_n)$.\endproof

The following statement ends the proof of Theorem \ref{tstar}.

\begin{proposition}\label{starpi3}
Let $(X,\dist)$ be a geodesic metric space $(*)$--ATG with respect
to a collection of subsets $\aaa$. If $(X,\aaa )$ moreover satisfy
properties $(\alpha_1)$ and $(\alpha_2)$, then in any asymptotic
cone of $X$ property $(\Pi_3)$ is satisfied.
\end{proposition}

\proof Let $\co{X;e,d}$ be an arbitrary asymptotic cone of $X$ and
let $\aau$ be the collection of limit sets of sequences from
$\aaa$. Since $(\alpha_1) \& (\alpha_2)\Rightarrow (\fq)$ it
follows that the sets in $\aau$ are geodesic. Also $(\alpha_1) \&
(\alpha_2)$ imply that $\aau$ satisfies $(T_1)$.

Let $\Delta$ be a non-trivial simple geodesic triangle in
$\co{X;e,d}$, whose edges $[x,y]$, $[y,z]$ and $[z,x]$ appear as
limits of sequences $[x_n,y_n']$, $[y_n,z_n']$ and $[z_n,x_n']$ of
geodesics in $X$. Then \uass $\dist (x_n,x_n'), \dist (y_n,y_n')$
and $\dist (z_n,z_n')$ are of order $o(d_n)$, while the lengths of
$[x_n,y_n']$, $[y_n,z_n']$ and $[z_n,x_n']$ are of order $O(d_n)$.
Let $T_n$ be a geodesic triangle with vertices $x_n,y_n,z_n$. We
denote its edges by $[u,v]$, with $u,v\in \{x_n,y_n,z_n\}$. The
three limit geodesics $\g^x=\lio{[y_n,z_n]}$,
$\g^y=\lio{[x_n,z_n]}$ and $\g^z=\lio{[x_n,y_n]}$ compose the
limit triangle $T=\lio{T_n}$.

\medskip

\textsc{Case 1.}\quad Assume that \uass $T_n$ is in case \bc. Then
there exists $a_n^1\in [x_n,y_n]$, $a_n^2\in [y_n,z_n]$ and
$a_n^3\in [z_n,x_n]$ such that the set $\{ a_n^1, a_n^2,a_n^3 \}$
has \uass diameter at most $2\sigma$ for some constant $\sigma$.
It follows that $\lio{a_n^1}=\lio{a_n^2}=\lio{a_n^3}=a$. The point
$a$ is on the three edges of $T$.

Without loss of generality we may assume that $a\not\in \{x,y\}$.
The fact that $a\neq x$ implies that either $\g^z \neq [x,y]$ or
$\g^y\neq [x,z]$. Property $(\Pi_2)$ implies that we may apply
Proposition \ref{braid} to $\call_1=[x,y]$, $\call_2=[x,z]$,
$\g_1=\g^z$, $\g_2=\g^y$ and to the intersection point $a\in
\g^z\cap \g^y$. We conclude that the $\calt$-bigon formed by
$\g^z$ and $\g^y$ of endpoints $a,x$ is contained in a subset
$A_x\in \aau$.

Similarly we deduce that the $\calt$-bigon formed by $\g^z$ and
$\g^x$ of endpoints $a,y$ is contained in a subset $A_y\in \aau$.

If $a=z$ then $\g_x\subset A_y$, $\g_y\subset A_x$. Also $a\not\in
[x,y]$, which by Lemma \ref{lbigon} implies that $a$ is contained
in the interior of a simple $\calt$-bigon formed by $[x,y]$ and
$\g^z$. By property $(\Pi_2)$ this $\calt$-bigon is contained in
some $A\in \aau$. The intersections $A\cap A_x$ and $A\cap A_y$
contain non-trivial sub-arcs of $\g^z$ therefore by $(T_1)$ we
conclude that $A=A_x=A_y$. The subset $A$ contains also $\g^z$.

Property $(\Pi_2)$ allows to apply Lemma \ref{convbigon} to the
pairs of arcs $(\g^x , [y,z])$, $(\g^y , [x,z])$ and $(\g^z ,
[x,y])$ and deduce that $\Delta = [x,y]\cup [y,z] \cup [z,x]$ is
contained in $A$.

\me

If $a\neq z$ then again by Proposition \ref{braid} the
$\calt$-bigon formed by $\g^x$ and $\g^y$ of endpoints $a,z$ is
contained in a subset $A_z\in \aau$. Since $a$ is not a vertex in
$\Delta$ it is contained in at most one edge of $\Delta$. Without
loss of generality we assume that $a\not\in [x,y]\cup [y,z]$.

The fact that $a\not\in [x,y]$, Lemma \ref{lbigon}, properties
$(\Pi_2)$ and $(T_1)$ imply as above that $A_x=A_y$. Likewise from
$a\not\in [y,z]$ we deduce that $A_y=A_z$. Thus $A=A_x=A_y=A_z$
contains $T$. Property $(\Pi_2)$ and Lemma \ref{convbigon} imply
that $\Delta$ is also contained in $A$.

\medskip

\textsc{Case 2.}\quad Assume that \uass $T_n$ is in case \bp. Then
there exist $A_n$ in $\aaa$ such that $\onn_\sigma (A_n)$
intersects all the edges of $T_n$. Moreover if $(x_n^2, y_n^1)$,
$(y_n^2, z_n^1)$ and $(z_n^2, x_n^1)$ are the pairs of entrance
and exit points from $\onn_\sigma (A)$ of $[x_n,y_n]$, $[y_n,z_n]$
and $[z_n,x_n]$ respectively, then $\dist (x_n^1\, ,\, x_n^2)\,
,\, \dist (y_n^1\, ,\, y_n^2)$ and $\dist (z_n^1\, ,\, z_n^2)$ are
less than $\delta$.

Let $x'=\lio{x_n^1}=\lio{x_n^2}$, $y'=\lio{y_n^1}=\lio{y_n^2}$ and
$z'=\lio{z_n^1}=\lio{z_n^2}$.

Assume that $\{x',y',z'\}$ has cardinal at most $2$. Assume for
instance that $x'= y'$. Then the point $a=x'=y'$ is in $\g^x \cap
\g^y \cap \g^z$. With the same argument as in Case 1 we deduce
that both $T$ and $\Delta$ are contained in some $A\in\aau$.

Assume now that $\{x',y',z'\}$ has cardinal $3$. The geodesic
triangle $T'$ of vertices $x',y',z'$ and with edges contained in
the edges of $T$ is included in the piece $A=\lio{A_n}$.

Proposition \ref{braid} implies that the $\calt$-bigon of
endpoints $x,x'$ formed by $\g^z$ and $\g^y$ is either trivial or
contained in some $A_x\in \aau$. Similarly, the $\calt$-bigon of
endpoints $y,y'$ formed by $\g^z$ and $\g^x$ is either trivial or
in some $A_y$, and the $\calt$-bigon of endpoints $z,z'$ formed by
$\g^x$ and $\g^y$ is either trivial or in some $A_z$.

If $x'\neq x$ then $x'$ cannot be contained both in $[x,y]$ and in
$[x,z]$. Suppose that $x'\not\in [x,y]$. Then $x'$ is in the
interior of a non-trivial simple $\calt$-bigon formed by $\g^z$
and $[x,y]$. This $\calt$-bigon is contained in some $B_x\in \aau$
by $(\Pi_2)$, and its intersections with $A_x$ and with $A$
contain a non-trivial sub-arc of $\g^z$. Hence $A_x=B_x=A$. Thus
the $\calt$-bigon of endpoints $x,x'$ is contained in $A$.

In the same way we obtain that the $\calt$-bigons of endpoints
$y,y'$ and $z,z'$ are contained in $A$. Thus in all cases
$T\subset A$, which by Lemma \ref{convbigon} implies that $\Delta
\subset A$.\endproof

\section{Quasi-isometric rigidity of relatively hyperbolic
groups}\label{sqirrh}

In this section we prove the following.

\begin{theorem}\label{atgrh}
Let $(G ,\dist )$ be an infinite finitely generated group endowed
with a word metric, which is asymptotically tree-graded with
respect to a collection $\aaa$ of subsets of $G$. Let $\varkappa $
be the maximum between the constant $M$ in property $(\beta_2)$
and the constant $\chi$ in property $(\beta_3)$ of $(G,\aaa)$.

Then the group $G$ is either hyperbolic or relatively hyperbolic
with respect to a family of subgroups $\{H_1,...,H_m\}$, such that
each $H_i$ is contained in $\nn_\varkappa (A_i)$ for some $A_i\in
\aaa$.
\end{theorem}

\begin{rmk}\label{hypapart}
If $G$ is hyperbolic then it is hyperbolic  relative to $H=\{1\}$.
Still, in this case one cannot state that $H$ is contained in some
$\nn_\varkappa (A)$ with $A\in \aaa$, because in the definition
that we adopt of asymptotically tree-graded metric spaces the
finite radius tubular neighborhoods of sets $A\in \aaa$ do not
cover the whole space (see Remark \ref{ratg}).
\end{rmk}

\proof The pair $(G,\aaa )$ satisfies properties $(\alpha_1)$,
$(\beta_2)$ and $(\beta_3)$.

By Corollary \ref{chyp}, if for $\theta>0$ and $\nu \geq 8$ from
$(\beta_3)$ either there exists no $(\theta , \nu)$--fat geodesic
hexagon in the Cayley graph of $G$, or the $(\theta , \nu)$--fat
geodesic hexagons have uniformly bounded diameter, then $G$ is
hyperbolic.

Assume from now on that for every $\eta >0$ there exists a
$(\theta , \nu)$--fat geodesic hexagon of diameter at least
$\eta$.

For $\varkappa$ as in the theorem and
 $\diam_\varkappa$ given by property $(\alpha_1)$ of $\aaa$,
consider the set
$$
\Phi = \left\{P \: \mbox{geodesic hexagon} \; ;\; P \: \mbox{is
}(\theta , \nu)\mbox{--fat}\, , \; \mbox{diameter}(P)\geq
\diam_\varkappa +1\right\}\, .
$$

Let $g\in G$. The metric space$(G,\dist)$ is asymptotically
tree-graded with respect to the collection of subsets $g\aaa =\{
gA\; ;\; A\in \aaa \}$, moreover the constants in the properties
$(\alpha_1)$, $(\beta_2)$ and $(\beta_3)$ are the same as for
$\aaa$.

Let $P\in \Phi$. Then $P$ is contained in $\nn_\varkappa (gA)$ for
some $A\in \aaa$. If $P$ is also contained in $\nn_\varkappa
(gA')$ for $A'\in \aaa$ then $\nn_\varkappa (A) \cap \nn_\varkappa
(A')$ has diameter at least the diameter of $P$, hence at least
$\diam_\varkappa +1$, consequently $A=A'$. Thus $P$ defines a map
$$
\cali_P : G \to \aaa \, ,\, \cali_P (g)=A\mbox{ such that
}P\subset \nn_\varkappa (gA) \, .
$$

We may then define
$$
\cali : \Phi \to Map\,(G,\aaa)\, ,\, \cali (P) = \cali_P \, ,
$$ where $Map\, (G,\aaa )$ is the set of maps from $G$ to $\aaa$.
Consider the equivalence relation on $\Phi$ induced by $\cali$,
that is
$$
P\sim P' \Leftrightarrow \cali(P)=\cali(P')\Leftrightarrow \forall
g\in G\, ,\, P\mbox{ and } P'\mbox{ are in the same }\nn_\varkappa
(g A)\, .
$$

Let $[P]$ be the equivalence class of a hexagon $P$ in $\Phi $. To
it we associate the set
$$
B[P]=\bigcap_{g\in G} \nn_\varkappa \left( g\, \cali_P (g) \right)
$$

\begin{proposition}\label{atgb}
The metric space $(G , \dist)$ is ATG \wrt
$$
\calb = \left\{ B[P]\; ;\; [P]\in \Phi /\sim  \right\}\, .
$$
\end{proposition}

\proof According to Proposition \ref{b3big} it suffices to prove
$(\alpha_1), (\beta_2)$ and $(\beta_3^\eta )$ for some $\eta
>0$. The proof relies on the simple remark that for every $r>0$,
$$
\nn_r \left( B[P] \right) \subset \bigcap_{g\in G} \nn_{r+
\varkappa} \left( g\, \cali_P (g) \right).
$$

$(\alpha_1)$\quad Let $[P]\neq [P']$, which is equivalent to the
fact that there exists $g_0\in G$ such that $P\subset \nk{g_0A}$
and $P'\subset \nk{g_0A'}$ with $A\neq A'$. For every $\delta >
0$, $$\nn_\delta \left( B[P] \right)\cap \nn_\delta \left(
B[P']\right) \subset \nn_{\delta +\varkappa} \left(g_0A\right)\cap
\nn_{\delta +\varkappa} \left(g_0A'\right)= g_0\left[ \nn_{\delta
+\varkappa} \left(A\right)\cap \nn_{\delta +\varkappa}
\left(A'\right) \right]\, .$$

Property $(\alpha_1)$ for $\aaa$ implies that the diameter of
$\nn_\delta \left( B[P] \right)\cap \nn_\delta \left(
B[P']\right)$ is uniformly bounded.

\me

$(\beta_2)$\quad Let $\epsilon$ be the constant appearing in
$(\beta_2)$ for $\aaa$. Take $\epsilon'=\frac{\epsilon}{2}$ and
take $M'= \frac{\epsilon +1}{\epsilon}\varkappa$. We prove that
$(\beta_2)$ holds for $\calb$ with the constants $\epsilon'$ and
$M'$.

Let $\g$ be a geodesic of length $\ell$ and let $[P]\in
\Phi/\sim\, $ be such that $\g (0)$ and $\g (\ell )$ are in
$\nn_{\epsilon' \ell} \left( B[P]\right)$. It follows that for
every $g\in G$, $\g (0)$ and $\g (\ell )$ are in $\nn_{\epsilon'
\ell + \varkappa} \left( g\, \cali_P (g) \right)$.

If $\varkappa \geq \frac{\epsilon}{2}\ell \Leftrightarrow \ell
\leq \frac{2\varkappa }{\epsilon }$ then $\g \subset
\onn_{\frac{\varkappa}{\epsilon }} \left( \{ \g (0),\g (\ell ) \}
\right)\subset \nn_{\frac{\varkappa}{\epsilon } +\varkappa }\left(
B[P]\right)=\nn_{M'}\left( B[P]\right)$.

Assume that $\varkappa < \frac{\epsilon}{2}\ell$. Then for every
$g\in G$ the geodesic $g\iv \g$ of length $\ell$ has its endpoints
in $\nn_{\epsilon \ell } \left(\cali_P (g) \right)$. Property
$(\beta_2)$ implies that $g\iv \g \left( \left[ \frac{\ell}{3}\,
,\, \frac{2\ell}{3} \right] \right)$ is contained in $\nn_M
\left(\cali_P (g) \right) \subset \nn_\varkappa \left(\cali_P (g)
\right)$.

We have thus obtained that for every $g\in G$, $\g \left( \left[
\frac{\ell}{3}\, ,\, \frac{2\ell}{3} \right] \right)$ is contained
in $\nn_\varkappa \left(g \cali_P (g) \right)$. It follows that
$\g \left( \left[ \frac{\ell}{3}\, ,\, \frac{2\ell}{3} \right]
\right)$ is contained in $B[P]$.

Property $(\beta_3^\eta)$ holds for the constants $\theta$ and
$\nu$ same as in $(\beta_3)$ for $\aaa$, for the constant $\chi$
equal to $0$, and $\eta =\diam_\varkappa +1$. Indeed every $P\in
\Phi$ is contained in $B[P]$.\endproof

\begin{lemma}[the group permutes the pieces]\label{perm}
\begin{itemize}
    \item[(1)] If $P\sim P'$ and $g\in G$ then $gP\sim gP'$.

    Consequently
    $G$ acts on the left on $\Phi/\sim \,$.
    \item[(2)] For every $g\in G$ and $P\in \Phi$, $gB[P]=B[gP]$.
\end{itemize}
\end{lemma}

\proof \textbf{(1)} The set $\cali_P(\gamma)$ is defined by
$P\subset \nk{\gamma \cali_P(\gamma)}$. For every $g\in G$,
$gP\subset \nk{g\gamma \cali_P(\gamma)}$, hence
$\cali_{gP}(g\gamma)=\cali_P(\gamma)$. From this can be deduced
that $P\sim P'\Rightarrow gP \sim gP'$.

\me

\textbf{(2)} The translate $gB[P]=\bigcap_{\gamma\in G}
\nn_\varkappa \left( g\gamma\, \cali_P (\gamma) \right)$ is equal
to $$\bigcap_{\gamma\in G} \nn_\varkappa \left( g\gamma\,
\cali_{gP}(g\gamma) \right)= \bigcap_{\gamma'\in G} \nn_\varkappa
\left( \gamma'\, \cali_{gP}(\gamma') \right)=B[gP]\, .$$
\endproof

The following statement finishes the proof of Theorem \ref{atgrh}.

\begin{proposition}[equivariant ATG structure implies relative
hyperbolicity]\label{equivatg}

Consider a finitely generated group endowed with a word metric $(G
,\dist )$, which is ATG with respect to a collection of subsets
$\calb$, such that $G$ permutes the subsets in $\calb$.

Then $G$ is either hyperbolic, or hyperbolic relative to a family
of subgroups $\{H_1,..., H_m\}$ such that for each $H_i$ there
exists a unique $B_i\in \calb$ satisfying $H_i\subset B_i\subset
\nn_\mathcal{K} (H_i)$ , where $\mathcal{K}$ is a constant
depending only on $(G,\dist )$ and $\calb$.
\end{proposition}



The proof is done in several steps.

\begin{lemma}\label{fin}
Finitely many subsets in $\calb$ contain $1$.
\end{lemma}

\proof By property $(\fq )$ of $\calb$ there exists $\tau>0$ such
that for any $x,y$ in some $B\in \calb$ any geodesic $[x,y]$ is
contained in $\nn_\tau \left( B\right)$. Property $(\alpha_1)$ for
$\calb$ implies that there exists $D_\tau$ such that for $B\neq
B'$, $\nn_\tau \left( B\right)\cap \nn_\tau \left( B'\right) $ has
diameter at most $D_\tau$.

Assume that $B\in \calb$ contains $1$ and has diameter at most
$3D_\tau$. Then $B\subset \overline{B}(1,3D_\tau)$. As
$\overline{B}(1,3D_\tau)$ is finite, only finitely many $B\in
\calb$ can be in this case.

Assume that $B$ contains $1$ and has diameter larger than
$3D_\tau$. Then $B$ contains some point $x$ with $\dist
(1,x)>3D_\tau$. The geodesic $[1,x]$ is contained in $\nn_\tau
\left( B\right)$ and it intersects the sphere around $1$ of radius
$2D_\tau$, $S(1, 2D_\tau)$. We define a map from the set $\left\{
B\in \calb\; ;\; 1\in B\, ,\, \diam \, B
>3D_\tau \right\}$ to the set of subsets of $S(1, 2D_\tau)$,
associating to each $B$ the non-empty intersection $\nn_\tau
\left( B\right) \cap S(1, 2D_\tau)$. By $(\alpha_1)$ and the
choice of $D_\tau$,  two distinct subsets $B,B'$ have disjoint
images by the above map, in particular the map is injective. Since
the set of subsets of $S(1, 2D_\tau)$ is finite, so is the
considered subset of $\calb$.\endproof

Let $\calf=\{ B_1,B_2,...,B_k\}$ be the set of $B\in \calb$
containing $1$. For every $i\in \{1,2,...,k\}$ let
$$\mathcal{I}_i=\{j\in \{1,2,...,k\}\mid \exists g\in G\mbox{ such
that }gB_i=B_j \}\, .$$ For every $j\in \mathcal{I}_i$ we fix
$g_j\in G$ such that $g_j B_i=B_j$.

\begin{n}\label{k}
Define the constants $\mathcal{K}_i=\max_{j\in \mathcal{I}_i}\dist
(1,g_j)$ and $\mathcal{K}=\max_{1\leq i\leq k}\mathcal{K}_i$.
\end{n}

\begin{lemma}\label{pstab}
For every $B\in \calb \,$ the stabilizer $\stab \left( B
\right)=\left\{ g\in G \mid gB=B \right\}$ is a subgroup of $G$
acting $\mathcal{K}$-transitively on $B$ (in the sense of
Definition \ref{tk}).
\end{lemma}

\proof Let $x$ and $b$ be arbitrary points in $B$. Both subsets
$b\iv B$ and $x\iv B$ contain $1$ and are in $\calb$. It follows
that $b\iv B=B_i$ and $x\iv B=B_j$ for some $i,j\in
\{1,2,...,k\}$. Since $b\iv xB_j= B_i$ it follows that $j\in
\mathcal{I}_i$ and that $B_j=g_jB_i$. The last equality can be
re-written as $x\iv B = g_j b\iv B$ which implies that $xg_jb\iv
\in \stab (B)$, hence that $x$ is at distance at most $\dist
(1,g_j)$ from $\stab (B) b$.\endproof

\begin{cor}\label{dinc}
For $i\in \{1,2,...,k\}$, $\stab (B_i)\subset B_i\subset
\nn_{\mathcal{K}} \left(\stab (B_i) \right)$.
\end{cor}

Let $\diam_{2\calk}$ be the uniform bound given by property
$(\alpha_1)$ for $(G,\calb)$ and $\delta = 2\calk $.

If all the subsets in $\calb$ have diameter at most
$\diam_{2\calk}+1$ then $G$ is hyperbolic by Corollary \ref{chyp}.
Thus, in what follows we may assume that $\calb$ contains subsets
of diameter larger than $\diam_{2\calk}+1$.

Denote by $\calb'$ the set of $B\in \calb $ of diameter larger
than $\diam_{2\calk}+1$. Proposition \ref{atgb} and Corollary
\ref{redd} imply that $G$ is ATG \wrt $\calb'$. Obviously $G$ also
permutes the subsets in $\calb'$.


Let $\calf'=\calf \cap \calb'$. Let $\calf_0$ be a subset of
$\calf'$ such that for every $B\in \calf'$, its orbit $G\cdot B$
intersects $\calf_0$ in a unique element. Such a subset can be
obtained for instance by considering one by one the elements $B_i$
in $ \calf'$, and deleting from $\calf'$ all $B_j$ with $j\in
\mathcal{I}_i,\, j\neq i$.

It follows that for every $B\in \calb'$, the orbit $G\cdot B$
intersects $\calf_0$ in only one element.

Let $\bar{B}_1,\dots, \bar{B}_m$ be the elements of $\calf_0$.

\begin{lemma}\label{bij}
For every $B\in \calb'$ there exists a unique $j\in \{1,2,...,m\}$
and a unique left coset $g\stab \left( \bB_j \right)$ such that
\begin{equation}\label{edinc}
g\stab \left( \bB_j \right)\subset B\subset \nn_{\mathcal{K}}
\left(g\stab \left( \bB_j \right) \right)\, .
\end{equation}
\end{lemma}

\proof \textit{Existence.} Let $g\in B$. Then $g\iv B\in \calb'$
and $1\in g\iv B$. Therefore $g\iv B = \bB_j$ for some $j\in
\{1,2,...,m\}$. Corollary \ref{dinc} implies the double inclusion
(\ref{edinc}).

\me

\textit{Unicity.} Assume that $g\stab \left( \bB_j \right)$ and
$g'\stab \left( \bB_l \right)$ both satisfy (\ref{edinc}), for
$j,l\in \{1,2,...,m\}$. Then
$$
g\bB_j \subset \nn_\calk \left( g\stab \left( \bB_j \right)
\right)\subset \nn_\calk \left( B \right)\subset \nn_{2\calk}
\left( g'\stab \left( \bB_l \right) \right)\subset \nn_{2\calk}
\left( g'\bB_l \right)\, .
$$

Both $g\bB_j$ and $g'\bB_l$ are in $\calb'$, in particular
$g\bB_j$ has diameter at least $\diam_{2\calk}+1$. Property
$(\alpha_1)$ implies that $g\bB_j = g'\bB_l$. According to the
definition of $\calf_0$ this can only happen if $j=l$. Then $g\iv
g'$ is in $\stab \left( \bB_j \right)$, and  $g'\stab \left( \bB_l
\right)$ coincides with $g\stab \left( \bB_j \right)$.\endproof

\begin{lemma}\label{cend}
The group $G$ is hyperbolic relative to $\left\{ H_1,...,H_m
\right\}$, where $H_j=\stab \left( \bB_j \right)$.
\end{lemma}

\proof The fact that $G$ is ATG \wrt $\calb'$, Lemma \ref{bij} and
Remark \ref{tna}, (2), imply that $G$ is ATG \wrt $\left\{ gH_j\;
;\; g\in G/H_j\, ,\, j\in \{1,2,...,m\}\right\}\, . $ In
particular by $(\fq)$ each $H_j$ is quasi-convex in $G$, hence
each $H_j$ is finitely generated. Theorem \ref{dsrh} implies that
$G$ is hyperbolic relative to $H_1,..., H_m$.

If $G= H_j= \stab \left( \bB_j \right)$ then Corollary \ref{dinc}
implies that $G=\bB_j$.\hspace*{\fill}$\square$

\bigskip

A consequence of Theorem \ref{atgrh} is the following.

\begin{theorem}\label{qir}
Let $G$ be a group hyperbolic relative to a family of subgroups
$\hh=\{H_1,...,H_n\}$. If a group $G'$ is $(L,C)$-quasi-isometric
to $G$ then $G'$ is hyperbolic relative to
$\hh'=\{H_1',...,H_m'\}$, where each $H_i'$ can be embedded
$(\lambda , \kappa)$-quasi-isometrically in $H_j$ for some
$j=j(i)\in \{1,2,...,n\}$, where $(\lambda , \kappa )$ depend on
$(L,C)$ and on $(G,\hh )$.
\end{theorem}

\proof If the group $G$ is finite then the group $G'$ is also
finite. We assume henceforth that both groups are infinite.

Let $\q$ be an $(L,C)$-quasi-isometry from $G$ to $G'$, and let
$\bar{\q}$ be its quasi-converse, such that $\dist (\q \circ
\bar{\q} , \mathrm{id}_{G'})\leq D$ and $\dist (\bar{\q}\circ \q ,
\mathrm{id}_{G})\leq D$, where $D=D(L,C)$. By Theorem \ref{dsrh},
$G$ is ATG \wrt the collection of left cosets $\aaa = \{ gH_i\;
;\; g\in G/H_i\, ,\, i\in \{1,2,...,n\}\}$. Theorem 5.1 in
\cite{DrutuSapir:TreeGraded} implies that $G'$ is ATG \wrt $\q
(\aaa )= \{ \q (A) \; ;\; A\in \aaa \}$. Moreover all constants
appearing in the properties $(\alpha_i)$, $i=1,2,3$, $(\beta_j)$,
$j=2,3,$ and $(\fq)$ for $(G', \q (\aaa ))$ can be expressed as
functions of $(L,C)$ and of the constants in the similar
properties for $(G,\aaa )$.

Theorem \ref{atgrh} implies that $G'$ is either hyperbolic or
relatively hyperbolic with respect to a family of subgroups
$\{H_1',...,H_m'\}$; moreover each $H_i'$ is contained in
$\nn_\varkappa \left( \q \left( A_i \right) \right)$ for some
$A_i\in \aaa$, where $\varkappa$ is a constant depending on
$(L,C)$, on the constant $M$ in $(\beta_2)$ for $(G,\aaa )$, and
on the constant $\chi $ in $(\beta_3)$ for $(G,\aaa )$.

Let $\pi_1 : \nn_\varkappa \left( \q \left( A_i \right) \right)
\to \q \left( A_i \right)$ be a map such that $\dist (x,
\pi_1(x))\leq \varkappa$. Then $\pi_1$ is a
$(1,2\varkappa)$--quasi-isometric embedding. Let $\pi_2 :
\nn_D(A_i)\to A_i$ be a $(1,2D)$--quasi-isometric embedding
constructed similarly. The restriction to $H_i'$ of $\pi_2 \circ
\bar{\q} \circ \pi_1$ is a $(\lambda , \kappa )$-quasi-isometric
embedding of $H_i'$ into $A_i=gH_j$, for some $j\in
\{1,2,...,n\}$, with $(\lambda , \kappa)$ depending on $(L,C)$,
$\varkappa$ and $D$.

If $G'$ is hyperbolic then $G'$ is relatively hyperbolic with
respect to $\{1\}\neq \{G'\}$ and all the statements in the
theorem hold.

If $G'=H_i'$ then $G'=\nn_\varkappa \left( \q \left( A_i \right)
\right)$, which implies that $G\subset \nn_C (\bar{\q}
(G'))\subset \nn_{L\varkappa +2C+D} (A_i)$. By Theorem \ref{dsrh},
this contradicts the fact that $G$ is (properly) hyperbolic
relative to $\hh$.\endproof



\begin{minipage}[t]{2.9 in}
\noindent Cornelia Dru\c tu\\ UFR de Math\'ematiques\\
Universit\'e de Lille 1\\ F-59655 Villeneuve d'Ascq \\FRANCE \\ Cornelia.Drutu@math.univ-lille1.fr\\
\end{minipage}

\end{document}